\documentclass[12pt]{article}
\usepackage[T1]{fontenc}
\usepackage[utf8]{inputenc}
\usepackage{lmodern}
\usepackage{microtype}
\usepackage{amsthm,amsmath,amssymb,bm}
\usepackage[authoryear]{natbib}
\usepackage{graphicx,booktabs}
\usepackage{needspace}
\usepackage{xcolor}
\usepackage{tikz}
\usetikzlibrary{arrows.meta,decorations.pathreplacing}
\usepackage[hidelinks]{hyperref}
\usepackage{bookmark}
\usepackage{orcidlink}
\def\spacingset#1{\renewcommand{\baselinestretch}{#1}\small\normalsize
}

\providecommand{\anon}{1}

\hypersetup{
    pdftitle={Species sampling priors through exact augmentation},
    pdfkeywords={Bayesian nonparametrics, Species sampling process,
    Finite-index representation, Prior construction,
    Size-biased invariance}
}

\ifnum\anon=1
\hypersetup{
    pdfauthor={Ramses H. Mena; Christos Merkatas;
    Theodoros Nicoleris; Carlos E. Rodriguez}
}
\else
\hypersetup{pdfauthor={}}
\fi

\newtheorem{theorem}{Theorem}[section]
\newtheorem{proposition}[theorem]{Proposition}
\newtheorem{corollary}[theorem]{Corollary}

\theoremstyle{definition}
\newtheorem{example}[theorem]{Example}

\newcommand{\N}{\mathbb N}
\newcommand{\Pp}{\mathbb P}
\newcommand{\E}{\mathbb E}
\newcommand{\dd}{\,\mathrm d}
\newcommand{\TV}{d_{\mathrm{TV}}}

\newcommand{\xib}{\bm\xi}
\newcommand{\wb}{\bm w}
\newcommand{\thetab}{\bm\theta}
\newcommand{\zb}{\bm z}
\newcommand{\kb}{\bm k}

\DeclareMathOperator{\Geom}{Geom}
\DeclareMathOperator{\Beta}{Beta}
\DeclareMathOperator{\NegBin}{NegBin}

\begin{document}
\spacingset{1}

\ifnum\anon=1

\title{\bfseries
Species sampling priors through exact augmentation}

\author{
Rams\'es H. Mena$^{1}$\,\orcidlink{0000-0001-9608-8059},
Christos Merkatas$^{2}$\,\orcidlink{0000-0001-7880-7443},\\
Theodoros Nicoleris$^{3}$,
Carlos E. Rodr\'iguez$^{1}$\,\orcidlink{0000-0002-4796-0784}\thanks{Corresponding author:
\href{mailto:carloserwin@sigma.iimas.unam.mx}
{carloserwin@sigma.iimas.unam.mx}}
}
\date{}
\maketitle

\vspace{-1.5em}

\begin{center}
\small
$^{1}$ Department of Probability and Statistics,
IIMAS, Universidad Nacional Aut\'onoma de M\'exico\\[0.3em]
$^{2}$ Department of Statistics and
Actuarial--Financial Mathematics,
University of the Aegean\\[0.3em]
$^{3}$ Department of Economics,
National and Kapodistrian University of Athens
\end{center}

\else

\bigskip
\bigskip
\bigskip
\begin{center}
{\LARGE\bfseries
Species sampling priors through exact augmentation
}
\end{center}
\medskip

\fi

\bigskip

\begin{abstract}
Discrete random probability measures are central to Bayesian nonparametric modeling. In practice, their use often favors a few tractable priors with full weak support. Full support alone, however, does not specify how prior mass reflects beliefs about clustering, recurrence,  discovery and other sample features. We develop a novel and general construction of species sampling process (SSP) priors from an index distribution and finite allocation probabilities, linked to these sample characteristics. Every proper SSP admits an exact two-stage representation through an almost surely finite index and a finite conditional measure. Averaging over the index recovers each realized SSP on every measurable set. Reversing this identity constructs priors whose tie probabilities, expected occupancy and discovery rates can be assessed without evaluating the infinite weights. We characterize compatibility among index laws, schedules and weights, and construct light-tailed, heavy-tailed and multimodal families with full weak support under suitable hyperpriors. A prior with paired weights combines full weak support with invariance under size-biased permutation in its defining order, an explicit finite index law and partition prediction. Exact simulation, recursive reweighting and finite posterior augmentation follow. Examples distinguish prior elicitation, exact averaging and fixed-level approximation.
\end{abstract}

\vfill
\newpage

\spacingset{1.8}

\section{Introduction}
Random probability measures express uncertainty about an unknown distribution, whether it describes observations directly or latent quantities in a hierarchical model. Their representations connect prior specification, prediction,  computation and generalizations. For example, the canonical Dirichlet process (DP) of \citet{Ferguson1973} can be obtained by normalizing a gamma process, a construction that motivated normalized completely random measures (CRMs), \citep[see, for example,][]{JamesLijoiPrunster2009}. Its P\'olya-urn and stick-breaking descriptions provide complementary ways to study and use the same prior \citep{BlackwellMacQueen1973,Sethuraman1994}, and subsequently were used to develop more general priors. 

Among the reasons to adopt a particular prior, 
mathematical tractability and full weak support are the most important. However, while full support assigns positive probability to every weak neighborhood of an admissible distribution, it does not specify how strongly the prior favors one pattern of recurrence or diversity over another. Prior information may instead concern whether two sampled individuals share a latent subgroup, how many distinct groups a sample of a given size is likely to contain, or whether a further observation will reveal an unobserved species with certain probability. For example, \citet{Lee2013} discuss clinical clustering where groups of comparable prevalence may be more plausible than a few dominant groups accompanied by many small ones. Encoding such beliefs directly in a law on infinitely many weights is difficult. Finite sampling ingredients offer a more direct connection between these beliefs and Bayesian prior elicitation, assessment and updating \citep[see, {\it e.g.},][]{NavarreteQuintanaMueller2008}.

Proper species sampling processes (SSPs) provide
a common framework for these purposes. They encompass the DP and
Pitman--Yor process (PY) \citep{Pitman1996,PitmanYor1997},
the broad class of homogeneous normalized CRMs
\citep{FavaroEtAl2016}, and geometric, beta-binomial,
exchangeable-fraction and Markov stick-breaking priors
\citep{FuentesGarcia2010,GilLeyvaMenaNicoleris2020,
GilLeyvaMena2023,GilLeyvaLijoiMenaPrunster2026}.
Their partition and predictive distributions connect these
weight laws to recurrence and discovery and, when computable,
provide a basis for comparing and eliciting priors
\citep{GnedinPitman2006,DeBlasi2015,LijoiPrunster2010}.
However, although these distributions are always well defined, they
need not admit tractable expressions.

Let $\mathsf X$ be the space of atom locations, assumed Polish and equipped with its Borel $\sigma$-field $\mathcal B$. Write $\wb=(w_j)_{j\geq1}$ for a random sequence of nonnegative weights and $\thetab=(\theta_j)_{j\geq1}$ for a random sequence of atom locations, drawn independently from a diffuse base probability measure $G_0$ and independently of $\wb$. A proper SSP has the representation
\[
 G=\sum_{j\geq1}w_j\delta_{\theta_j},\qquad
 \sum_{j\geq1}w_j=1\quad\text{almost surely},
 \qquad \theta_j\stackrel{\mathrm{iid}}\sim G_0,
\]
where $\delta_\theta$ is unit mass at $\theta$. The joint law of the weights, together with $G_0$, specifies the prior on $G$ \citep{Pitman1996,GilLeyvaLijoiMenaPrunster2026}. An \emph{enumeration} labels the weight--atom pairs $(w_j,\theta_j)$ by
$j=1,2,\ldots$, without requiring the weights to be ranked or the atom space to be ordered. We keep this enumeration fixed, distinguishing the ordered description from the measure $G$, which is unchanged by simultaneous relabeling of its weight--atom pairs.

When $G$ is used as a mixing distribution, observations $x_1,\ldots,x_n$ satisfy
\begin{equation}\label{eq:ssp-mixture}
 x_i\mid G\stackrel{\mathrm{iid}}\sim f_G,
 \qquad f_G(x)=\int f(x\mid\theta)G(\dd\theta)
             =\sum_{j\geq1}w_j f(x\mid\theta_j),
\end{equation}where $f(\cdot\mid\theta)$ is the sampling kernel.

Equivalently, each observation $x_i$ is generated by selecting an atom $\theta_j$ with probability $w_j$ and then sampling from $f(\cdot\mid\theta_j)$. Writing $z_i$ for the label of the selected atom, we have
\[
 \Pp(z_i=j\mid\wb)=w_j,\qquad
 x_i\mid z_i,\thetab\sim f(\cdot\mid\theta_{z_i}),
 \qquad i=1,\ldots,n.
\]
As the labels are conditionally i.i.d. given the shared weights, their equality pattern defines an exchangeable partition. Its exchangeable partition probability function (EPPF) is symmetric in the block sizes \citep{Pitman1996,Pitman2006}. Note that this partition law requires neither size-biased invariance of
the ordered weights nor a closed-form EPPF. Its exchangeability
and consistency follow from conditionally i.i.d.\ sampling from
a shared SSP, whereas predictive rules specified directly must
satisfy conditions ensuring that they define an exchangeable
partition law \citep{Lee2013}.

A representation also determines which posterior methods are available. Marginal methods integrate out the random measure through an EPPF or predictive rule, whereas conditional methods retain its weights and atoms. Slice sampling \citep{Walker2007,Kalli2011} and retrospective simulation \citep{PapaspiliopoulosRoberts2008} handle the infinite posterior through almost surely finite updates. Error-controlled truncations instead produce approximations at a prescribed accuracy \citep{Arbel2019}.

We develop an exact decomposition that makes the finite conditional measures and their mixing probabilities explicit. Let $\xib=(\xi_j)_{j\geq1}$ be a strictly decreasing sequence of positive numbers tending to zero. We call it a \emph{visibility schedule} because, at an auxiliary level $u>0$, label $j$ is visible when $u\leq\xi_j$. The eligible labels consequently form a finite initial set. \citet{Kalli2011} introduced this schedule mechanism in a generalization of the slice sampler of \citet{Walker2007}. Here it identifies a discrete two-stage sampling law whose inverse provides a method of prior construction.

As we will see, the construction we propose first draws an almost surely finite index $K$ and then selects among the first $K$ atoms
with probabilities proportional to $w_j/\xi_j$. Averaging the
resulting finite conditional measures over $K$ recovers the
realized SSP on every measurable set. For a fixed weight law,
the schedule can therefore be varied to obtain different
representations while preserving the prior and its marginal
posterior.

Reversing the identity turns the finite sampling specification
into a method of prior construction. An index law and a finite
allocation rule determine a weight sequence, with random weights
obtained by drawing their parameters once per random measure.
Finite conditional calculations then allow the induced prior
to be assessed through tie probabilities, expected occupancy,
discovery probabilities and other sample characteristics.
With uniform allocation, for example, pairwise coincidence is
evaluated by averaging over the two observation-specific indices,
without calculating the infinite weight sequence.
Different allocation rules can accompany the same index law
and induce different priors, whereas a compatible conditional
index law and schedule determine the realized weights uniquely.

The representation covers every proper SSP in any fixed enumeration, regardless of monotonicity or independence of the stick-breaking fractions. Equal allocations recover earlier decreasing-weight constructions \citep{FuentesGarcia2010,DeBlasiMartinezMenaPrunster2020,HatjispyrosMerkatasWalker2023}, while unequal allocations yield new weight structures. We characterize inverse compatibility and give hyperprior conditions for full weak support, including examples with finite expected index. We also construct weights invariant under size-biased permutation in their defining order. In particular, splitting each Dirichlet-process weight equally between two species with independent locations yields a paired prior with full weak support, explicit partition and predictive laws, and a finite marginal index generator. Combining this splitting operation within the fragmentation framework of \citet[Section~5.4]{Pitman2006} with a defining order generated through size-biased sampling
\citep{PitmanISBP1996,PermanPitmanYor1992} gives an invariant
finite-index specification outside the class of homogeneous
normalized CRMs.

Here, the finite index $K$ determines how many labels are
available for an individual allocation, rather than the total
number of components in a finite-mixture model. Observations
share the same weights and atoms, while each has its own index.
As in finite and infinite mixture models, the nonempty sets
$\{i:z_i=j\}$ form the sample partition, whose number of blocks
is $c_n=|\{z_1,\ldots,z_n\}|$ \citep{Nobile2004,McCullaghYang2008}. Accordingly, prior elicitation can focus on the partition and
predictive laws induced jointly by the index distribution and
allocation rule, without interpreting $K$ as the number of
subpopulations.

We first establish the representation and then reverse it to construct priors, examine their support and relate their finite sampling ingredients to predictive beliefs. The subsequent analysis develops simulation and posterior updates, with numerical comparisons showing how prior specification differs from computational representation. Supporting proofs and further examples appear in the Supplementary Material \citep{MenaSSPSupp2026}, to which we use lettered sections as references.

\section{Exact two-stage representation}
To obtain an exact finite sampling rule, we first condition on the weights and atoms of a proper SSP. Without loss of generality, we fix such a realization and a positive, strictly decreasing schedule $\xib$ tending to zero. The schedule may be specified independently of the weights (exogenous) or determined by their realization (endogenous). Write $R_k=\sum_{j>k}w_j$ for the remaining original mass beyond label $k$, with $R_0=1$. Define
\begin{equation}\label{eq:ak-sk-intro}
 a_j=w_j/\xi_j,\qquad s_k=\sum_{j=1}^k a_j.
\end{equation}
Here $a_j$ is the reweighted mass assigned to label $j$, and $s_k$ is the total reweighted mass of the first $k$ labels.

\begin{theorem}[Setwise representation]\label{thm:main}
Let $G=\sum_{j\geq1}w_j\delta_{\theta_j}$ be a proper SSP with
weight sequence $\wb$ and atom sequence $\thetab$ in a fixed
enumeration. Given a visibility schedule $\xib$, let $a_j$ and
$s_k$ be as in \eqref{eq:ak-sk-intro} and define
\begin{equation}\label{eq:pk}
p_k(\wb,\xib)=\Pp(K=k\mid\wb,\xib)
 =(\xi_k-\xi_{k+1})s_k,\qquad k\geq1.
\end{equation}
We write $p_k$ when the conditioning sequences $(\wb,\xib)$
are understood. For each $k$ with $s_k>0$, define
\begin{equation*}
G_k^\star(A)=\sum_{j=1}^k
\widetilde w_j^{(k)}\delta_{\theta_j}(A),
\qquad
\widetilde w_j^{(k)}=\frac{a_j}{s_k}
=\frac{w_j/\xi_j}{\sum_{h=1}^k w_h/\xi_h}.
\end{equation*}
When $s_k=0$, set $G_k^\star=\delta_{\theta_1}$. Since $p_k=0$,
this convention has no effect on the sampling law.
Then, almost surely, $(p_k)_{k\geq1}$ is a probability mass
function and
\begin{equation}\label{eq:barycenter}
G(A)=\sum_{k\geq1}p_k(\wb,\xib)G_k^\star(A),
\qquad A\in\mathcal B.
\end{equation}
The corresponding two-stage draw first generates
$K\mid(\wb,\xib)\sim(p_k)$. Given $K=k$ with $p_k>0$, it then
chooses a label $J\in\{1,\ldots,k\}$ with probabilities
$\widetilde w_1^{(k)},\ldots,\widetilde w_k^{(k)}$ and returns
the atom $\theta_J$. The selected label has marginal probabilities
\begin{equation*}
\Pp(J=j\mid\wb,\xib)=w_j,
\end{equation*}
so $\theta_J$ is an exact draw from the realized measure $G$.
\end{theorem}

\begin{proof}
Fix a realization of the weights and atoms of $G$, together
with an admissible schedule $\xib$. Since $\xi_m\to0$,
\[
 \sum_{k\geq j}(\xi_k-\xi_{k+1})=\xi_j.
\]
Consider the candidate joint probabilities
\begin{equation}\label{eq:jointJK}
 \Pp(J=j,K=k\mid\wb,\xib)
 =a_j(\xi_k-\xi_{k+1})\mathbf1\{j\leq k\}.
\end{equation}
The right-hand side is nonnegative, and summing over $k\geq j$
gives $a_j\xi_j=w_j$. Since $\sum_jw_j=1$,
\eqref{eq:jointJK} defines a joint probability.  Its $K$-marginal is
$(\xi_k-\xi_{k+1})s_k=p_k$, while conditioning on $K=k$ with
$p_k>0$ gives the allocation probabilities $a_j/s_k$ for
$j\leq k$. Levels with $s_k=0$ have $p_k=0$ and add
nothing. Finally, Tonelli's theorem gives, for every
$A\in\mathcal B$,
\[
 \sum_{k\geq1}p_kG_k^\star(A)
 =\sum_{j\geq1}a_j\delta_{\theta_j}(A)
       \sum_{k\geq j}(\xi_k-\xi_{k+1})
 =\sum_{j\geq1}w_j\delta_{\theta_j}(A)=G(A).
\]
The argument holds for every admissible realization and
therefore almost surely.
\end{proof}
In Theorem~\ref{thm:main}, $K$ determines the available labels $1,\ldots,K$, whose relative selection probabilities are
specified by the coefficients $a_j=w_j/\xi_j$. For fixed
weights and schedule, each label retains its coefficient as
$k$ increases, so adding labels can change the normalizing
constant $s_k$ while preserving the probability ratios among
previously included labels with positive weight. The finite
representation depends on the chosen enumeration, although
relabeling the weight--atom pairs leaves $G$ unchanged.

Note that the identity holds for every realisation and measurable set, even with coincident atom locations when $G_0$ has atoms. We focus on pairwise distinct locations, which occur almost surely for diffuse $G_0$ and are required for some of the results below. Note also that a common positive rescaling of the schedule leaves $p_k$ and $G_k^\star$ unchanged.

Independent draws from $G$, including the latent parameters in \eqref{eq:ssp-mixture}, use Theorem~\ref{thm:main} separately conditional on the shared $(\wb,\thetab,\xib)$. Draw independent indices $K_i$ with probabilities $(p_k)$, then independently select $z_i\leq K_i$ with probabilities $\widetilde w_j^{(K_i)}$. For measurable sets $A_1,\ldots,A_n\in\mathcal B$, the joint sampling probability is
\begin{equation}\label{eq:independent-indices-product}
\Pp(\theta_{z_i}\in A_i\text{ for all }i=1,\ldots,n
            \mid\wb,\thetab,\xib)=\sum_{k_1,\ldots,k_n}\prod_{i=1}^n
            p_{k_i}G_{k_i}^\star(A_i)
       =\prod_{i=1}^nG(A_i).
\end{equation}
Drawing $x_i$ from $f(\cdot\mid\theta_{z_i})$ recovers the
original mixture model. Using a common $K$ for the entire
sample, as if it were the number of components in a
finite-mixture model, would instead generate all latent
parameters from the same $G_K^\star$. Averaging over this
index would generally introduce additional dependence.
Supplement~B.1 derives this distinction.

For strictly decreasing positive weights, the endogenous choice $\xi_j=w_j$ gives
\begin{equation*}
 p_k=k(w_k-w_{k+1}),\qquad
 G_k^\star=\frac1k\sum_{j=1}^k\delta_{\theta_j},\qquad
 w_j=\sum_{k\geq j}\frac{p_k}{k}.
\end{equation*}
Averaging the equal finite allocations recovers the unequal weights, as in \citet{FuentesGarcia2010,DeBlasiMartinezMenaPrunster2020,HatjispyrosMerkatasWalker2023}. General schedules require no monotonicity and also apply to familiar stick-breaking SSPs. If $w_j=v_jR_{j-1}$, with $R_{j-1}=\prod_{\ell<j}(1-v_\ell)$  and $0<v_j<1$ almost surely, choosing $\xi_j=R_{j-1}$ gives $a_j=v_j$ and $p_k=w_k\sum_{j\leq k}v_j$. Thus Theorem~\ref{thm:main} covers, in particular, the DP, the two-parameter PY process and exchangeable stick-breaking (ESB) priors through their respective stick-breaking fractions \citep{Sethuraman1994,PitmanYor1997,GilLeyvaMena2023}.

\section{SSP priors from a prescribed index distribution}
To construct a prior, we specify an index law
$\boldsymbol\pi=(\pi_k)_{k\geq1}$ with $\pi_k>0$ and
$\sum_k\pi_k=1$, together with a finite allocation rule,
and obtain the weights by averaging. The notation $\pi_k$
distinguishes these prescribed probabilities from
$p_k(\wb,\xib)$ in \eqref{eq:pk}. Compatibility requires
\begin{equation*}
 \pi_k=(\xi_k-\xi_{k+1})\sum_{j\leq k}w_j/\xi_j.
\end{equation*}
Any compatible pair determines the third ingredient, with the schedule unique up to scale. For random weights, prescribing $\Pp(K=k\mid\wb,\xib)=\pi_k$ fixes the law at every realization, whereas prescribing only $\Pp(K=k)=\pi_k$ fixes its prior average. We develop the direct construction before considering marginal calibration for an existing SSP.

\subsection{Construction through finite allocation probabilities}
Let $\boldsymbol a=(a_j)_{j\geq1}$ be a sequence of positive coefficients whose values are chosen deterministically or
drawn jointly with $\boldsymbol\pi$ from a prior independent
of the atom sequence. Conditional on $\boldsymbol a$ and
$K=k$, the allocation rule selects label $j\leq k$ with
probability $a_j/s_k$, where $s_k=\sum_{h\leq k}a_h$.
Thus equal coefficients give uniform allocation, whereas
increasing coefficients favor later labels and alternating
coefficients permit repeated local preferences. Normalizing
the same coefficients over the first $k$ labels preserves
the probability ratios $a_j/a_h$ whenever both labels are
available and remains valid even if
$\sum_{j\geq1}a_j=\infty$.

\begin{theorem}[SSP construction from an index law]\label{thm:count-prior}
Let $\boldsymbol\pi$ be a strictly positive probability mass function on $\N$ and let $a_j\in(0,\infty)$ for every $j$. Put $s_k=\sum_{j=1}^k a_j$ and define
\begin{equation}\label{eq:count-prior-map}
\xi_j=\sum_{k=j}^{\infty}\frac{\pi_k}{s_k},
\qquad
w_j=a_j\xi_j,
\qquad
\widetilde w_j^{(k)}=\frac{a_j}{s_k},\quad j\leq k.
\end{equation}
Then $\xi_1>\xi_2>\cdots\downarrow0$, $w_j>0$, $\sum_jw_j=1$, and
\begin{equation}\label{eq:count-prior-forward}
(\xi_k-\xi_{k+1})\sum_{j=1}^k\frac{w_j}{\xi_j}=\pi_k.
\end{equation}
For any atom sequence in $\mathsf X$,
\begin{equation}\label{eq:count-prior-barycenter}
G=\sum_{j\geq1}w_j\delta_{\theta_j}
  =\sum_{k\geq1}\pi_kG_k^\star,
\qquad \mbox{with}\qquad
G_k^\star=\sum_{j=1}^k\frac{a_j}{s_k}\delta_{\theta_j}.
\end{equation}
Consequently, any joint prior on $(\boldsymbol\pi,\boldsymbol a)$ satisfying these assumptions almost surely, together with atoms drawn independently and identically from a diffuse $G_0$ independently of that pair, defines a proper SSP. Conversely, every positive-weight representation in Theorem~\ref{thm:main} is obtained from \eqref{eq:count-prior-map} by taking $a_j=w_j/\xi_j$ and $\pi_k=p_k(\wb,\xib)$.
\end{theorem}
The proof is given in Supplement~A.1.

An appealing feature of Theorem~\ref{thm:count-prior} is that
it generates observations through $\boldsymbol\pi$ without
evaluating the infinite sums defining $w_j$. For fixed
$\boldsymbol\pi$, rescaling $\boldsymbol a\mapsto c\boldsymbol a$,
$c>0$, preserves both $\widetilde w_j^{(k)}=a_j/s_k$ and the
induced weights.

When $\pi_k=\pi_k(\eta)$, the values of $(\eta,\boldsymbol a)$
may be fixed or drawn jointly from a prior independent of the
atoms. A single realization of $(\eta,\boldsymbol a,\thetab)$
is shared by the sample, with observations generated
conditionally independently through
\begin{equation*}
 k_i\mid\eta\sim\boldsymbol\pi(\eta),\qquad
 \Pp(z_i=j\mid k_i,\boldsymbol a)
 =\frac{a_j}{s_{k_i}}\mathbf1\{j\leq k_i\},
 \qquad
 x_i\mid z_i,\thetab\sim f(\cdot\mid\theta_{z_i}),
\end{equation*}
for $i=1,\ldots,n$. Writing $\boldsymbol x=(x_1,\ldots,x_n)$,
$\kb=(k_1,\ldots,k_n)$ and $\zb=(z_1,\ldots,z_n)$, the augmented
likelihood is
\begin{equation}\label{eq:count-prior-likelihood}
 \mathcal L_{\rm ind}(\boldsymbol x,\zb,\kb\mid\eta,\boldsymbol a,\thetab)
 =\prod_{i=1}^n\pi_{k_i}(\eta)\frac{a_{z_i}}{s_{k_i}}
 f(x_i\mid\theta_{z_i})\mathbf1\{z_i\leq k_i\}.
\end{equation}

Evaluating $\mathcal L_{\rm ind}$ requires only
$a_1,\ldots,a_{\max_i k_i}$. Since all observations share
$(\eta,\boldsymbol a)$, the joint prior on these quantities
is updated using the likelihood of the full sample.

\par\medskip\noindent\textit{Prior-predictive specification.}\enspace
We can calibrate $(\boldsymbol\pi,\boldsymbol a)$ to prior
beliefs about pairwise coincidence, the expected number of
observed species and the probability of discovering a new species.
Conditional on $(\boldsymbol\pi,\boldsymbol a)$, let $K_1,K_2$ be independent with law $\boldsymbol\pi$, and draw $J_i$ with probabilities $a_j/s_{K_i}$, for $j\leq K_i$. Their coincidence probability is
\begin{equation}
 \begin{aligned}
 T_2(\boldsymbol\pi,\boldsymbol a)
 &=\Pp(J_1=J_2\mid\boldsymbol\pi,\boldsymbol a)
   =\sum_{j\geq1}w_j^2\\
 &=\E\left\{\frac{\sum_{j\leq\min(K_1,K_2)}a_j^2}
                       {s_{K_1}s_{K_2}}\,\middle|\,\boldsymbol\pi,\boldsymbol a\right\}.
 \end{aligned}
\end{equation}
Under uniform allocation,
$\Pp(J_1=J_2\mid K_1=k,K_2=\ell,\boldsymbol a=\boldsymbol1)
=1/\max(k,\ell)$, so
\begin{equation}
 \Pp(J_1=J_2\mid\boldsymbol\pi,\boldsymbol a=\boldsymbol1)
 =T_2(\boldsymbol\pi,\boldsymbol1)
 =\E\{1/\max(K_1,K_2)\mid\boldsymbol\pi\}.
\end{equation}
Thus, under uniform allocation, the index law directly determines
a fundamental clustering probability without requiring the infinite
weights to be calculated. The expectation still averages over an
unbounded index law, but each simulated evaluation is finite and exact.

Notice that for arbitrary positive coefficients $\boldsymbol a$,
$J_i\mid(\boldsymbol\pi,\boldsymbol a)
\stackrel{\mathrm{iid}}{\sim}(w_j)$. The expected
occupancy and the probability of discovering a new label
therefore satisfy
\begin{equation}\label{eq:prior-occupancy-discovery}
 \begin{aligned}
 \E(c_n\mid\boldsymbol\pi,\boldsymbol a)
   &=\sum_{j\geq1}\{1-(1-w_j)^n\},\\
 \Pp(J_{n+1}\notin\{J_1,\ldots,J_n\}
      \mid\boldsymbol\pi,\boldsymbol a)
   &=\sum_{j\geq1}w_j(1-w_j)^n.
 \end{aligned}
\end{equation}

These calculations describe sample characteristics under the
prior, before conditioning on observed data.
Supplement~A.6 gives equivalent
finite-index formulas and derivations, including expected
new species, and illustrates posterior updating of an
elicited coincidence probability through the augmented
likelihood. Averaging the prior summaries over a hyperprior
on the shared parameters can guide elicitation.

It is worth noting that, fixing the conditional index law also restricts the support of the induced prior, since $w_1\geq\pi_1>0$. However, an hyperprior on its parameters allows the index probabilities and weights to vary across random measures, whereas calibrating a marginal index law for an existing SSP preserves its weight prior.

\subsection{Compatibility and the inverse relations}
For specified weights let $W_k=\sum_{j\leq k}w_j$, $Q_k=\sum_{j\leq k}\pi_j$ and $W_0=Q_0=0$. The inclusion $\{K\leq k\}\subseteq\{J\leq k\}$ gives $Q_k\leq W_k$. This condition is also sufficient when the finite allocation probabilities can be chosen freely, as proved in Lemma~A.1 of Supplement~A.1. The strictly decreasing schedule formulation has the following sharper characterization.

\begin{proposition}[Inverse determination of the schedule]\label{prop:count-inverse-schedule}
Let $w_j>0$, $\sum_jw_j=1$, and $\pi_k>0$, $\sum_k\pi_k=1$. A positive, strictly decreasing schedule tending to zero for which $p_k(\wb,\xib)=\pi_k$ exists if and only if
\begin{equation}\label{eq:count-strict-order}
Q_k<W_k\qquad\text{for every finite }k.
\end{equation}
It is unique up to a common positive factor. With $\xi_1=1$, it is determined recursively by
\begin{equation}\label{eq:count-inverse-schedule}
\xi_{k+1}=\xi_k\frac{W_k-Q_k}{W_k-Q_{k-1}},\qquad k\geq1.
\end{equation}
\end{proposition}
The cumulative difference has the probabilistic interpretation
\[
 W_k-Q_k=\Pp(J\leq k,\ K>k\mid\wb,\xib)=\xi_{k+1}s_k.
\]The proof in Supplement~A.1 derives the recursion, verifies its converse and establishes uniqueness.

For an SSP prior specified in advance, the same identities instead
determine a compatible representation. Let $\wb\sim\Pi$, with
$\overline w_j=\E_\Pi(w_j)$ and
$\overline W_k=\sum_{j\leq k}\overline w_j$. Keeping $\Pi$ fixed,
we look for a deterministic schedule yielding a prescribed
\emph{marginal} law for $K$, while allowing the conditional index
probabilities to vary with $\wb$. Since \eqref{eq:pk} is linear
in the weights for a fixed schedule, averaging over $\Pi$ reduces
the problem to the mean sequence.

\begin{corollary}[Marginal calibration of a fixed SSP]\label{cor:count-marginal}
Assume $\overline w_j>0$ for every $j$. There is a deterministic positive, strictly decreasing schedule tending to zero with $\Pp(K=k)=\pi_k$ if and only if $Q_k<\overline W_k$ for every finite $k$. Furthermore, it is unique up to scale and, with $\xi_1=1$, satisfies
\begin{equation*}
\xi_{k+1}=\xi_k
\frac{\overline W_k-Q_k}{\overline W_k-Q_{k-1}}.
\end{equation*}
The weight prior remains $\Pi$, and $K$ has conditional probabilities $p_k(\wb,\xib)$. Equivalently, if $K\sim\boldsymbol\pi$ is generated first, the conditional law of the weights is
\begin{equation}\label{eq:count-tilted-prior}
 \Pi(\dd\wb\mid K=k)
 =\frac{s_k(\wb)}{\overline s_k}\,\Pi(\dd\wb),
 \qquad
 s_k(\wb)=\sum_{j=1}^k\frac{w_j}{\xi_j},\quad
 \overline s_k=\sum_{j=1}^k\frac{\overline w_j}{\xi_j}.
\end{equation}
\end{corollary}
The proof is given in Supplement~A.1.

Thus calibration retains the original weight prior. One can draw $\wb\sim\Pi$ and then $K$ with probabilities $p_k(\wb,\xib)$, or draw $K\sim\boldsymbol\pi$ first and use \eqref{eq:count-tilted-prior}. Both orders specify the same joint law, where $K$ and $\wb$ are generally dependent. Alternatively, prescribing the index law and schedule determines the total reweighted mass at each level.

\begin{proposition}[Inverse determination of the weights]\label{prop:count-inverse-weights}
Let $\boldsymbol\pi$ be a probability mass function on $\N$, with $\pi_k>0$ and $\sum_{k\geq1}\pi_k=1$, and let $\xi_1>\xi_2>\cdots\downarrow0$. Define
\begin{equation*}
t_0=0,\qquad t_k=\frac{\pi_k}{\xi_k-\xi_{k+1}},\quad k\geq1.
\end{equation*}
There is a nonnegative probability sequence $\wb$ yielding  $p_k(\wb,\xib)=\pi_k$ if and only if $(t_k)$ is nondecreasing. In that case the weights are uniquely determined by
\begin{equation*}
w_j=\xi_j(t_j-t_{j-1}),\qquad j\geq1,
\end{equation*}
and their reweighted normalizing constant is $s_k=t_k$. All weights are positive precisely when $(t_k)_{k\geq0}$
is strictly increasing, although the weights themselves
need not be monotone in their defining order. Every compatible pair also satisfies
\begin{equation}\label{eq:index-schedule-barrier}
 \xi_j\leq \frac{H_j}{t_1},\qquad\mbox{with}\qquad
 H_j=\sum_{k\geq j}\pi_k,\quad j\geq1.
\end{equation}
\end{proposition}
The proof is given in Supplement~A.1. A compatible schedule $\xib$ and conditional probabilities
$\pi_k=\Pp(K=k\mid\wb,\xib)$ uniquely determine the realized
weights $\wb$. For a fixed deterministic schedule, prescribing
only the marginal probabilities $\pi_k=\Pp(K=k)$ instead
determines the mean weights $\E(w_j)$, while allowing different
joint prior distributions for $\wb$ with those same means.

A joint prior supported on compatible pairs $(\boldsymbol\pi,\xib)$
therefore defines an SSP, and Theorem~\ref{thm:count-prior}
provides a coefficient parametrization that guarantees
compatibility. Retaining a previously chosen SSP, however,
requires the cumulative inequalities to hold for every
$k\geq1$, so checking a finite range alone is insufficient.
The bound \eqref{eq:index-schedule-barrier} further constrains
the schedule through the index survival probabilities,
excluding polynomial schedules for geometric or Poisson
conditional index laws. This relation motivates the next
construction, which obtains a compatible schedule directly
from these survival probabilities.

\subsection{Survival schedules and full-support SSPs}

Indeed, the inverse criterion suggests choosing the schedule
from the previously defined survival probabilities
$H_j=\sum_{k\geq j}\pi_k$. A strictly convex transformation
of this decreasing sequence yields the increasing totals
required by Proposition~\ref{prop:count-inverse-weights},
whether the index probabilities are monotone or multimodal.

To assess the resulting prior, let $\mathsf X_0=\operatorname{supp}(G_0)$ and write $\mathcal P(\mathsf X_0)$ for the space of probability measures on $\mathsf X_0$. Full weak support allows positive prior probability for every weak neighborhood of every measure in this space. Under our diffuse independent-atom assumptions, \citet[Corollary~3]{BissiriOngaro2014} characterize this property by
\begin{equation}\label{eq:full-support-criterion}
 \Pp\!\left(\sup_j w_j<\varepsilon\right)>0
 \quad\text{for every }\varepsilon>0.
\end{equation}
Thus, full support requires arbitrarily small largest weights.
The next result identifies a class of transformations for which
spreading the index probabilities achieves this.

\begin{theorem}[Survival-transform SSP priors]\label{thm:survival-full-support}
Let $\boldsymbol\pi=\boldsymbol\pi(\eta)$ be a strictly positive probability mass function on $\N$, and let $\phi:[0,1]\to[0,\infty)$ be continuous, strictly increasing and strictly convex, with $\phi(0)=0$. For $H_j=\sum_{k\geq j}\pi_k$, define
\begin{equation}\label{eq:survival-construction}
 \xi_j=\phi(H_j),\qquad
 t_0=0,\quad t_j=\frac{\pi_j}{\phi(H_j)-\phi(H_{j+1})},\qquad
 a_j=t_j-t_{j-1},\quad w_j=\xi_j a_j.
\end{equation}
Then $\xi_1>\xi_2>\cdots\downarrow0$, $a_j,w_j>0$ and $\sum_jw_j=1$. With independent $G_0$ atoms, these weights define a proper SSP whose exact two-stage rule is
\begin{equation}\label{eq:survival-two-stage}
 K\mid\eta\sim\boldsymbol\pi(\eta),\qquad
 \Pp(J=j\mid K=k,\eta)=\frac{a_j}{t_k}\mathbf1\{j\leq k\},\qquad
 \Theta=\theta_J.
\end{equation}
In particular, $G=\sum_k\pi_kG_k^\star$, where $G_k^\star=\sum_{j\leq k}a_j\delta_{\theta_j} /t_k$.

Suppose also that $\phi$ is differentiable on $(0,1)$ and $g=1/\phi'$ extends continuously to $[0,1]$, with $g(1)=0$. If $\eta$ is independent of the atoms and its prior law satisfies
\begin{equation}\label{eq:survival-small-maximum}
 \Pp\{\Delta(\eta)<\varepsilon\}>0\quad\text{for every }\varepsilon>0,
 \qquad \mbox{where}\qquad \Delta(\eta)=\sup_k\pi_k(\eta),
\end{equation}
then the SSP has full weak support on $\mathcal P(\mathsf X_0)$. Its expected index is $\E_\eta\sum_k k\pi_k(\eta)$, with $+\infty$ allowed.
\end{theorem}
For a fixed transformation in this class, $g(1)=0$ is also
necessary for full weak support, since $g(1)>0$ gives
$\sup_jw_j\geq w_1\geq\phi(1)g(1)>0$.
Since $H_j$ decreases, strict convexity of $\phi$ makes the totals
$t_j$ in \eqref{eq:survival-construction} strictly increase, ensuring
$a_j=t_j-t_{j-1}>0$. To establish full weak support, we bound these
increments using the modulus of continuity of $g$ over intervals
of length at most $2\Delta$. Supplement~A.5 provides the proofs and
verifies the required hyperprior conditions for the example hierarchies.

An explicit choice is $\phi(x)=1-(1-x)^\vartheta$, $0<\vartheta<1$. Here $g(x)=(1-x)^{1-\vartheta}/\vartheta$, and
\begin{equation}\label{eq:survival-endpoint-bound}
 w_1=\pi_1^{1-\vartheta},\qquad
 \sup_jw_j\leq\vartheta^{-1}(2\Delta)^{1-\vartheta}.
\end{equation}
For $\vartheta=1/2$, put $\xi_0=1$ and $\xi_j=1-\sqrt{1-H_j}$. The SSP weights and finite allocation probabilities then simplify to
\begin{equation}\label{eq:survival-square-root}
 \begin{aligned}
 w_j&=\xi_j(\xi_{j-1}-\xi_{j+1}),\\
 \widetilde w_j^{(k)}&=\frac{\xi_{j-1}-\xi_{j+1}}{2-\xi_k-\xi_{k+1}},\quad 1\leq j\leq k,
 & R_k&=\xi_k\xi_{k+1}.
 \end{aligned}
\end{equation}
With $\xi_0=\xi_1=1$, \eqref{eq:survival-square-root} also
holds at $j=1$, and finite allocation requires only
$H_1,\ldots,H_{K+1}$. To generate a sample, retain one
realization of $(\eta,\thetab)$ and draw the indices
independently conditional on $\eta$. For the mixture case,
\eqref{eq:survival-two-stage} is followed by a draw from
$f(\cdot\mid\theta_J)$.

The weight formula \eqref{eq:survival-square-root} accommodates light-tailed, heavy-tailed and multimodal index laws.
Table~S.1 in
Supplement~A.5 presents thirteen
illustrative index families whose hyperpriors satisfy
\eqref{eq:survival-small-maximum}, yielding SSP priors with
full weak support, although their weights are not size-biased
invariant in their defining order. Note that the choice of
transformation matters, since $\phi(x)=x^c$ with fixed
$c>1$ forces $w_1>1/c$ and precludes full weak support.

\subsection{Light-tailed index laws}
The negative-binomial and shifted-Poisson index laws both
yield SSP priors under suitable finite allocation rules,
although they differ in their ability to represent an
existing prior. For a DP in its usual stick-breaking order,
the negative-binomial family provides compatible marginal
index laws, whereas the tail of a fixed-parameter
shifted-Poisson law is too light for such calibration.

\begin{example}[Negative-binomial]\label{ex:count-nb}
Let $K-1\sim\NegBin(r,p)$ for $r>1$ and $0<p<1$. Set
\begin{equation}\label{eq:count-nb-pi}
 \pi_k=\frac{\Gamma(k+r-1)}{\Gamma(r)\Gamma(k)}
       p^r(1-p)^{k-1},
 \qquad \xi_j=(1-p)^{j-1}.
\end{equation}
Proposition~\ref{prop:count-inverse-weights} gives
\begin{equation}\label{eq:count-nb-weights}
 w_j=\frac{\Gamma(j+r-2)}{\Gamma(r-1)\Gamma(j)}
       p^{r-1}(1-p)^{j-1}.
\end{equation}
The allocation probabilities $\widetilde w_j^{(k)}$ are
proportional to $\Gamma(j+r-2)/\Gamma(j)$ for $1\leq j\leq k$
and do not depend on $p$. At $r=2$, they are uniform,
$\widetilde w_j^{(k)}=1/k$, while the induced SSP weights
are geometric, $w_j=p(1-p)^{j-1}$. At $r=3$,
\begin{equation*}
 w_j=jp^2(1-p)^{j-1},\qquad
 \widetilde w_j^{(k)}=\frac{2j}{k(k+1)}.
\end{equation*}
For $p<1/2$, the weights initially increase and then decrease.
For $a_0,b_0>0$, a $\Beta(a_0,b_0)$ prior on $p$ gives
full weak support, with finite expected index if and only
if $a_0>1$, as proved in Supplement~A.3.
Supplement~A.2 gives the
corresponding conjugate update.
\end{example}

The geometric case $r=2$ also allows prior elicitation through
the probability that two draws select the same label.
Let $T=\Pp(J_1=J_2\mid p)=p/(2-p)$, where $J_1$ and $J_2$
are labels drawn independently conditional on $p$.
Specifying $T\sim\Beta(a,b)$ induces a prior on
$p=2T/(1+T)$. The negative-binomial mean gives
$\E(K\mid p)=1+2(1-p)/p=1/T$, so averaging over the
beta prior yields
\[
 \E(T)=\frac{a}{a+b},\qquad
 \E(K\mid T)=\frac1T,\qquad\mbox{and}\qquad
 \E(K)=\E(T^{-1})=\frac{a+b-1}{a-1}\quad(a>1).
\]
Every such beta prior gives full weak support because it
assigns positive probability arbitrarily near $T=0$.
For instance, $\Beta(2,8)$ gives a mean coincidence
probability of $0.2$ and $\E(K)=9$. All observations share
the same $T$, whose prior is updated through
\eqref{eq:count-prior-likelihood}.

\begin{example}[Shifted-Poisson]\label{ex:count-poisson}
For $K=1+N$, $N\sim\mathrm{Poisson}(\nu)$, and $a_j=1$,
\begin{equation*}
 \pi_k=e^{-\nu}\frac{\nu^{k-1}}{(k-1)!},\qquad
 w_j=\xi_j=\frac{\Pp(N\geq j)}{\nu},\qquad
 \widetilde w_j^{(k)}=1/k.
\end{equation*}
A prior $\nu\sim\mathrm{Gamma}(a_0,b_0)$, with rate
$b_0$, gives full weak support and $\E(K)=1+a_0/b_0$.
The conjugate update and support argument appear in
Supplements~A.2 and
A.3, respectively. For fixed $\nu$, the shifted-Poisson law is too light-tailed
for marginal calibration to a DP in its usual stick-breaking
order, but gamma mixing can remove this obstruction.
For a DP with total mass $\alpha$, 
$\nu\sim\mathrm{Gamma}(2,1/\alpha)$ gives the compatible
marginal law $K-1\sim\NegBin(2,1/(\alpha+1))$.
The resulting negative-binomial tail remains too light for
a positive-discount PY prior, whose mean residual mass decays
polynomially in its usual stick-breaking order.
These tail restrictions, proved in
Supplement~A.2, concern marginal
calibration while retaining an existing SSP and do not
limit the direct prior constructions.
\end{example}

\subsection{Heavy-tailed index laws and their induced weights}
Let $\rho,b>0$, write $\mathrm B(x,y)=\Gamma(x)\Gamma(y)/\Gamma(x+y)$ for the beta function, and set
\begin{equation}\label{eq:bg-index}
 \pi_k(\rho,b)=\frac{\mathrm B(b+k-1,\rho+1)}{\mathrm B(b,\rho)},\qquad k\geq1.
\end{equation}
It arises by drawing $P\sim\Beta(\rho,b)$ and then $K-1\mid P\sim\Geom_0(P)$, where $Z\sim\Geom_0(p)$ has $\Pp(Z=s)=p(1-p)^s$, $s=0,1,\ldots$. Notice that each observation-specific index uses a fresh beta variable, conditional on the same sample-level parameters $(\rho,b)$. Its survival function and mean are
\begin{equation}
 \Pp(K>k\mid\rho,b)=\frac{\mathrm B(b+k,\rho)}{\mathrm B(b,\rho)},\qquad
 \E(K\mid\rho,b)=1+\frac{b}{\rho-1}\quad(\rho>1).
\end{equation}
At $b=1$ this is the Yule--Simon law. 

\begin{proposition}[An SSP induced by a beta-geometric index law]\label{prop:bg-ssp}
For the index law \eqref{eq:bg-index}, the schedule
\[
 \xi_j=\frac{\mathrm B(b+j-1,\rho+b)}
             {\mathrm B(b,\rho+b)},\qquad j\geq1,
\]
is compatible and, with $v=\rho/(b+\rho)$, induces
\begin{equation}\label{eq:bg-weights}
 w_j=v\mathbf1\{j=1\}+(1-v)\pi_j(\rho,b),\qquad j\geq1.
\end{equation}
These weights decrease strictly and sum to one, with
$\widetilde w_j^{(k)}\propto w_j/\xi_j$ for $1\leq j\leq k$.
Their residual mass $R_k=\sum_{j>k}w_j$ and asymptotic decay
satisfy
\begin{equation}\label{eq:bg-tail}
 R_k=\frac{b}{b+\rho}
      \frac{\mathrm B(b+k,\rho)}{\mathrm B(b,\rho)},\qquad
 w_j\sim\frac{b\Gamma(\rho+1)}
              {(b+\rho)\mathrm B(b,\rho)}j^{-(\rho+1)}.
\end{equation}
Together with i.i.d.\ $G_0$ atoms independent of the weights,
they define a proper SSP.
\end{proposition}
See Supplement~A.1 for the proof.
Keeping the beta-geometric index law fixed, we can obtain
nonmonotone weights by changing only the second-stage
allocation probabilities.

\begin{example}[Nonmonotone weights at a fixed index law]\label{ex:alternating-ssp}
Retain \eqref{eq:bg-index} and choose $a_{2m-1}=1$, $a_{2m}=\gamma>1$. Then
\begin{equation}\label{eq:alternating-weights}
 s_k=\lceil k/2\rceil+\gamma\lfloor k/2\rfloor,\qquad
 \xi_j=\sum_{k\geq j}\frac{\pi_k(\rho,b)}{s_k},\qquad
 w_j=a_j\xi_j,\qquad \widetilde w_j^{(k)}=a_j/s_k.
\end{equation}
Here $\lceil k/2\rceil$ and $\lfloor k/2\rfloor$ count
odd and even labels in $1,\ldots,k$ by rounding $k/2$
up and down, respectively. These weights define a proper SSP by
Theorem~\ref{thm:count-prior}. Writing
$\overline a=(1+\gamma)/2$ and
$C_{\rho,b}=\Gamma(\rho+1)/\mathrm B(b,\rho)$, we have
\begin{equation}\label{eq:alternating-tail}
 \xi_j\sim\frac{C_{\rho,b}}{\overline a(\rho+1)}
              j^{-(\rho+1)},\qquad
 \frac{w_{2m}}{w_{2m-1}}\to\gamma,\qquad
 \frac{w_{2m+1}}{w_{2m}}\to\gamma^{-1}.
\end{equation}
These limits, proved in
Proposition~A.3 of
Supplement~A.3, imply that the weights
eventually alternate between increases and decreases,
although the schedule decreases. This behavior is generated
directly by the allocation rule, without reordering the
preceding family.
\end{example}
The criterion \eqref{eq:full-support-criterion} clarifies
the role of the scale parameter $b$. With fixed parameters,
both families have deterministic weights and lack full weak
support, despite their random atoms. For fixed $\rho>1$ and,
in the alternating family, $\gamma>1$, drawing $b\sim H_b$
once per random measure gives full weak support exactly when
$H_b$ has unbounded support above, allowing the largest weight
to be arbitrarily small with positive probability.
In both families, $\E(K)=1+\E_{H_b}(b)/(\rho-1)$, so a gamma
prior provides full weak support with finite expected index.
These results follow from Theorem~A.2 and
Proposition~A.3 in
Supplement~A.3.

Combining Zipf or discrete-Pareto index laws with $a_j=2j-1$
gives further explicit SSP constructions with polynomially
decaying weights. Their formulas and tail behavior are
derived in Supplement~A.2, with
the implications for full weak support and expected index
examined in Supplement~A.3.

\subsection{Size-biased invariance and marginal index specification}

Size-biased invariance connects the defining weight order
with discovery order \citep{PitmanISBP1996,PermanPitmanYor1992}.
A size-biased permutation samples the labels of the positive
weights $\wb$ successively without replacement, with
probabilities proportional to the remaining weights. Thus, 
an ordered weight law is invariant under size-biased
permutation (ISBP) when this operation leaves the law of
the entire sequence unchanged. This property is appealing as it expresses partition probabilities
through joint moments of finitely many weights, providing
a route to explicit EPPFs and predictive rules. 
When the stick-breaking fractions are independent, ISBP
characterizes the two-parameter PY family among
proper infinite weight sequences \citep{PitmanISBP1996}. Allowing dependent stick-breaking fractions extends the class
of SSPs with ISBP weights beyond the PY family. Indeed, homogeneous normalized CRMs admit ISBP representations when their normalized jumps are taken in size-biased order \citep{PermanPitmanYor1992}. The stick-breaking representation of \citet{FavaroEtAl2016} and Theorem~\ref{thm:main} give their exact two-stage representation. 

We now develop an explicit finite-index specification of an
SSP outside the homogeneous normalized-CRM class, with ISBP
weights and full weak support.
Notice first that an infinite ISBP law satisfies $
 \Pp\!(\sum_{j\leq h}w_j<\varepsilon)>0$,
  for $h\in\N$ and $\varepsilon>0$,
whereas $\Pp(K=k\mid\wb)=\pi_k$ requires
$\sum_{j\leq h}w_j\geq\sum_{k\leq h}\pi_k$ a.s.
This incompatibility, proved in
Supplement~A.4, motivates specifying
$\pi_k=\Pp(K=k)$ and the conditional weight law $\Pi_k$.

For $\wb\sim\Pi$, deterministic $\xib$ and
$\overline w_j=\E_\Pi(w_j)$, these laws take the form
\begin{equation}
 \pi_k=(\xi_k-\xi_{k+1})
       \sum_{j\leq k}\frac{\overline w_j}{\xi_j},
 \qquad
 \Pi_k(\dd\wb)=\frac{p_k(\wb,\xib)}{\pi_k}\,\Pi(\dd\wb).
\end{equation}
Averaging the conditional weight laws $\Pi_k$ over
$K\sim\boldsymbol\pi$ recovers the prior $\Pi$, as expressed
by $\sum_k\pi_k\Pi_k=\Pi$ in
Corollary~\ref{cor:count-marginal}.
For a given sample, all observations share one draw of $\wb$,
with indices sampled independently given $\wb$ according
to $p_k(\wb,\xib)$, so
\eqref{eq:independent-indices-product} applies.

\subsection{An ISBP prior with paired weights}
\label{sec:paired-prior}

Fix $\alpha>0$ and let
$w_m^{\mathrm{DP}}=V_m\prod_{\ell<m}(1-V_\ell)$, with
$V_m\stackrel{\mathrm{iid}}\sim\Beta(1,\alpha)$, be
DP stick-breaking weights. Splitting each weight equally
between two independent atom locations gives
\begin{equation}\label{eq:paired-G}
 G^{(2)}_\alpha
 =\sum_{m\geq1}\frac{w_m^{\mathrm{DP}}}{2}
       \bigl(\delta_{Y_{m,1}}+\delta_{Y_{m,2}}\bigr),
 \qquad Y_{m,r}\stackrel{\mathrm{iid}}\sim G_0,
\end{equation}
where $G_0$ is diffuse and the atoms are independent of the
weights. This equal splitting is a special case of the
fragmentation operator in
\citet[Section~5.4, equation~(5.12)]{Pitman2006}.
It models pairs of equally prevalent species without
requiring similar parameters.

The same prior can be generated directly in size-biased
order. Before generating $w_j$, let $U$ be the combined
weight of pairs with neither member yet selected, and
let $\mathcal H$ record the unselected partner's weight
from each pair represented once among
$w_1,\ldots,w_{j-1}$. Thus $\mathcal H$ has at most
$j-1$ entries. Initially $(U,\mathcal H)=(1,\varnothing)$, and the unselected mass is
$r=U+\sum_{b\in\mathcal H}b=1-\sum_{\ell<j}w_\ell$,
giving selection probabilities
\begin{equation}\label{eq:paired-source}
 \Pp(\text{new pair}\mid U,\mathcal H)=\frac Ur,
 \qquad
 \Pp(\text{partner }b\mid U,\mathcal H)=\frac br,
 \quad b\in\mathcal H.
\end{equation}
For a new pair, draw an independent $V\sim\Beta(1,\alpha)$
and assign $w_j=UV/2$ to its first member, retaining the
same weight in $\mathcal H$ for its partner and update
$U$ to $U(1-V)$. Selecting an existing entry
$b\in\mathcal H$ gives $w_j=b$ and removes that entry,
leaving $U$ unchanged. These independent beta draws follow from the DP's size-biased stick-breaking representation
\citep[Theorem~1.1]{PermanPitmanYor1992}. Write $\wb\sim\Pi^{(2)}_\alpha$ for the resulting sequence and attach independent $G_0$ atoms.

A compatible marginal index law can also be specified
directly. For $0<q<1$, take $\xi_j=q^{j-1}$ and draw
$X\sim\operatorname{Unif}(0,1)$. Conditional on $X=x$, generate independently
\begin{equation*}
 B\sim\operatorname{Bern}\!\left(1-(2-x)^{-2}\right),\quad
\!\!\! N_1\sim\operatorname{Pois}\!\left(2\alpha\log(2-x)\right),\quad
\!\!\! N_2\sim\operatorname{Pois}\!\left(-\alpha\log\{x(2-x)\}\right).
\end{equation*}
With an independent $S_0\sim\Geom_0(1-q)$, put
\begin{equation}\label{eq:paired-index-law}
 D=1+B+N_1+2N_2,\qquad K=D+S_0,\qquad
 \pi_k=\Pp(K=k)
      =(1-q)\sum_{j=1}^k\overline w_jq^{k-j},
\end{equation}
where $\overline w_j=\Pp(D=j)$. Thus $K$ requires finitely many draws. The conditional weight law is
\begin{equation}\label{eq:paired-kernel}
 \Pi^{(2)}_{\alpha,k}(\dd\wb)
 =\frac{(1-q)\sum_{j=1}^kw_jq^{k-j}}{\pi_k}
       \Pi^{(2)}_\alpha(\dd\wb).
\end{equation}
After drawing $K=k$ from $\boldsymbol\pi$, generate the first
$k$ weights from their marginal under $\Pi^{(2)}_{\alpha,k}$,
attach independent $G_0$ atoms and select $J\leq k$ with
$\widetilde w_j^{(k)}\propto w_j/q^{j-1}$.

\begin{theorem}[A full-support invariant prior with paired weights]
\label{thm:paired-prior}
Fix $\alpha>0$ and $q\in(0,1)$, and set $\xi_j=q^{j-1}$.
The weight construction based on $(U,\mathcal H)$ generates
a proper SSP with the same law as $G^{(2)}_\alpha$ and
weights that are ISBP in their defining order.
This SSP has full weak support on $\mathcal P(\mathsf X_0)$
and is not the normalization of a homogeneous drift-free
CRM with deterministic L\'evy intensity and diffuse
independent locations.
Moreover, $\overline w_j=\E_{\Pi^{(2)}_\alpha}(w_j)$, so
\eqref{eq:paired-kernel} defines a conditional weight law
compatible with \eqref{eq:paired-index-law}.
The index $K$ has finite moments of every positive order, with
\begin{equation}\label{eq:paired-index-mean}
 \E(D)=2\alpha+\frac32,\qquad
 \E(K)=2\alpha+\frac32+\frac q{1-q}.
\end{equation}
\end{theorem}

Since $\sup_jw_j=\tfrac12\sup_m w_m^{\mathrm{DP}}$ can be
arbitrarily small with positive probability,
\eqref{eq:full-support-criterion} gives full weak support
at fixed $\alpha$. Each distinct weight occurs exactly
twice almost surely, which is incompatible with the
Poisson jump structure of a homogeneous CRM
\citep{Kingman1967}. Also $2w_1\sim\Beta(1,\alpha)$ implies
$w_1<1/2$ a.s., excluding PY laws with fixed
parameters \citep{PitmanYor1997}.
Supplement~A.4 proves
Theorem~\ref{thm:paired-prior} and gives a sampler for the
required weights. Here $\alpha$ changes the prior, while
$q$ controls only its representation.

The paired construction also gives explicit partition and
predictive laws by averaging over possible pairings of the
observed species. Let $\mathcal P$ be a partition of
$\{1,\ldots,n\}$ with block sizes
$\boldsymbol n=(n_1,\ldots,n_c)$. A partial matching $M$ is a possibly empty collection of disjoint pairs of block indices, each representing two species that share one original DP weight. Write
$\mathcal U(M)$ for the unmatched indices, whose partners
remain unobserved, and set
\[
 W_\alpha(M;\boldsymbol n)
 =(2\alpha)^{c-|M|}
   \prod_{i\in\mathcal U(M)}\Gamma(n_i)
   \prod_{\{i,j\}\in M}\Gamma(n_i+n_j).
\]

\begin{proposition}[Partition law and prediction with latent pairings]
\label{prop:paired-eppf}
The SSP of Theorem~\ref{thm:paired-prior} has EPPF
\begin{equation}\label{eq:paired-eppf}
 \Pp(\mathcal P_n=\mathcal P)
 =\mathfrak p^{(2)}_\alpha(\boldsymbol n)
 =\frac{\Gamma(\alpha)}{2^n\Gamma(\alpha+n)}
       \sum_M W_\alpha(M;\boldsymbol n).
\end{equation}
Given the latent sample
$\thetab_{\zb}=(\theta_{z_1},\ldots,\theta_{z_n})$,
with distinct values $\theta_1^*,\ldots,\theta_c^*$
and counts $\boldsymbol n$, the matching $\mathcal M$
has posterior law
\begin{equation}\label{eq:paired-matching}
 \Pp(\mathcal M=M\mid\thetab_{\zb})
 =\frac{W_\alpha(M;\boldsymbol n)}
        {\sum_{M'}W_\alpha(M';\boldsymbol n)}.
\end{equation}
Write $n_{g_M(i)}=n_i+n_j$ if $\{i,j\}\in M$
and $n_{g_M(i)}=n_i$ if $i\in\mathcal U(M)$.
Conditional prediction satisfies
\begin{equation}\label{eq:paired-predictive}
 \begin{aligned}
 \Pp(\theta_{z_{n+1}}=\theta_i^*
       \mid\thetab_{\zb},\mathcal M=M)
   &=\frac{n_{g_M(i)}}{2(\alpha+n)},\\
 \Pp(\theta_{z_{n+1}}\notin\{\theta_1^*,\ldots,\theta_c^*\}
       \mid\thetab_{\zb},\mathcal M=M)
   &=\frac{\alpha+\tfrac12
             \sum_{i\in\mathcal U(M)}n_i}{\alpha+n}.
 \end{aligned}
\end{equation}
Averaging over \eqref{eq:paired-matching} gives prediction
conditional on $\thetab_{\zb}$, with newly observed locations
distributed according to $G_0$.
\end{proposition}

Discovery can therefore arise from an entirely unobserved
pair or from the unseen partner of an observed species.
Supplement~C.4 proves
Proposition~\ref{prop:paired-eppf} and provides a recursion
for the matching sum.

The parameter $\alpha$ can also be elicited through the
prior coincidence probability
\begin{equation}\label{eq:paired-tie}
 \Pp(\theta_{z_1}=\theta_{z_2})=\frac1{2(\alpha+1)}.
\end{equation}
A specified value $t\in(0,1/2)$ corresponds to
$\alpha=1/(2t)-1$. Averaging over possible latent samples
gives the prior-predictive discovery probability
\begin{equation}\label{eq:paired-discovery}
 d_n=\Pp(\theta_{z_{n+1}}\notin
             \{\theta_{z_1},\ldots,\theta_{z_n}\})
 =\frac{\alpha}{2^n}\sum_{b=0}^n
       \binom nb\frac1{\alpha+b}.
\end{equation}
Summing these probabilities yields the prior expected
number of observed species,
\begin{equation}\label{eq:paired-occupancy}
 \E(c_n)
 =\sum_{m=0}^{n-1}d_m
 =\frac{2\alpha}{2^n}\sum_{b=0}^n\binom nb
       \{\psi(\alpha+b)-\psi(\alpha)\},
\end{equation}
where $\psi=\Gamma'/\Gamma$ is the digamma function.
This equals twice the expected occupancy of a DP with
concentration $\alpha$, averaged over a
$\operatorname{Binomial}(n,1/2)$ sample size.
The same discovery probabilities give
$\E(M_{1,n})=nd_{n-1}$, where $M_{1,n}$ counts species
observed exactly once.

\section{Sampling, reweighting and posterior computation}
With the weights and schedule fixed, \eqref{eq:jointJK} permits
sampling either the finite index $K$ or the atom label $J$ first.
We use this law to derive simulation rules, describe mass
redistribution between successive finite measures, and compare
those measures with $G$. Applying it independently to each
observation then gives the posterior augmentation and identifies
the finite sets of labels required for its updates.

\subsection{Selected label and additional visibility}
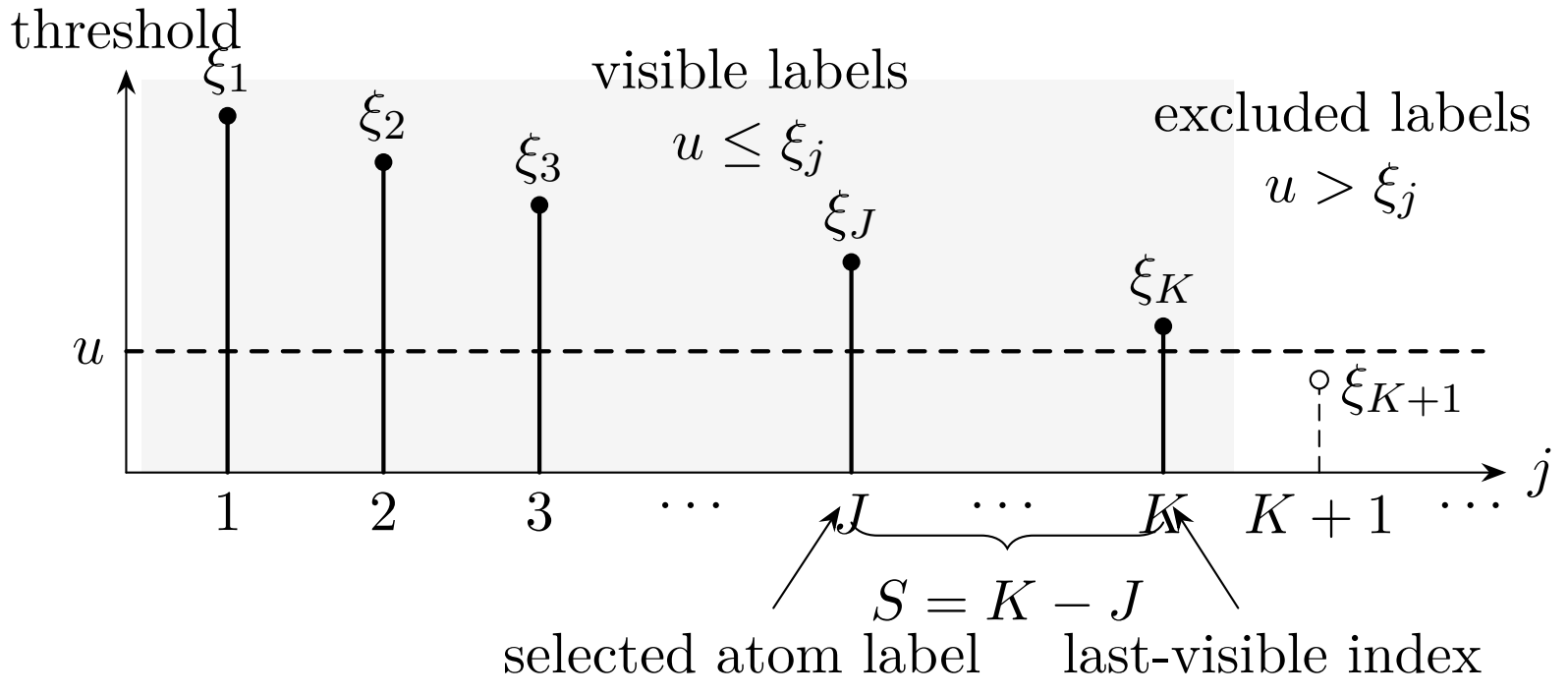
\begin{figure}[!htbp]
\centering
\begin{tikzpicture}[x=1.05cm,y=.48cm,>=Stealth,
  every node/.style={font=\linespread{1}\small},line cap=round]
  \fill[black!4] (.45,0) rectangle (7.45,5.5);
  \draw[->] (.35,0)--(9.2,0) node[right] {$j$};
  \draw[->] (.35,0)--(.35,5.65) node[above] {threshold};
  \foreach \j/\h in {1/5.0,2/4.35,3/3.75,5/2.95,7/2.05}{
    \draw[thick] (\j,0)--(\j,\h);
    \fill (\j,\h) circle (1.7pt);
  }
  \draw[densely dashed] (8,0)--(8,1.30);
  \draw[fill=white] (8,1.30) circle (1.7pt);
  \node[above] at (1,5.0) {$\xi_1$};
  \node[above] at (2,4.35) {$\xi_2$};
  \node[above] at (3,3.75) {$\xi_3$};
  \node[above] at (5,2.95) {$\xi_J$};
  \node[above] at (7,2.05) {$\xi_K$};
  \node[right,yshift=-2pt] at (8,1.30) {$\xi_{K+1}$};
  \draw[dashed,thick] (.35,1.7)--(9.05,1.7);
  \node[left] at (.35,1.7) {$u$};
  \node[align=center] at (4.35,5.05) {visible labels\\$u\leq\xi_j$};
  \node[align=center] at (8.15,4.45) {excluded labels\\$u>\xi_j$};
  \foreach \j/\lab in {1/1,2/2,3/3,4/{\cdots},5/J,6/{\cdots},7/K,8/{K+1},9/{\cdots}}
    \node[below] at (\j,0) {$\lab$};
  \draw[decorate,decoration={brace,mirror,amplitude=5pt}]
    (5,-.70)--(7,-.70) node[midway,below=7pt] {$S=K-J$};
  \draw[->] (4.50,-1.90)--(4.93,-.47);
  \node[anchor=north,align=center] at (4.3,-1.95) {selected atom label};
  \draw[->] (7.48,-1.90)--(7.06,-.47);
  \node[anchor=north,align=center] at (7.7,-1.95) {last-visible index};
\end{tikzpicture}
\caption{A visibility threshold $u$ determines the finite endpoint $K=L(u)$. The selected label $J$ lies in $\{1,\ldots,K\}$, and $S=K-J$ counts the additional label positions beyond $J$ that remain visible.}\label{fig:visibility}
\end{figure}

The decomposition $K=J+S$ separates the selected label $J$ from the additional visible positions $S$. At a slice level $0<u\leq\xi_1$, define
\begin{equation}
 L(u)=\max\{j\geq1:u\leq\xi_j\},\qquad
 L(u)=k\quad\text{if and only if}\quad
 u\in(\xi_{k+1},\xi_k].
\end{equation}
At level $u$, the available labels are $1,\ldots,L(u)$, and the selected label $J$ may be any one of them, as shown in Figure~\ref{fig:visibility}.

For a positive-weight label $j$, conditioning
\eqref{eq:jointJK} on $J=j$ gives
\begin{align}
 \Pp(S=s\mid J=j,\xib)
   &=\frac{\xi_{j+s}-\xi_{j+s+1}}{\xi_j},
   &&s=0,1,\ldots,\label{eq:S-pmf}\\
 \Pp(S\geq t\mid J=j,\xib)
   &=\frac{\xi_{j+t}}{\xi_j},
   &&t=0,1,\ldots.\label{eq:S-tail}
\end{align}

Equivalently, draw $J\sim\wb$, sample $S$ from \eqref{eq:S-pmf}
and set $K=J+S$, so $S=0$ means the selected atom is the last
visible one. For $\xi_j=cq^{j-1}$, $S\sim\Geom_0(1-q)$ is
independent of $J$ given the weights and schedule, although this
independence need not hold marginally when $q$ is random.
Supplement~B.1 gives the general stop--continue decomposition.

Recall that $R_k=\sum_{j>k}w_j$ and $R_0=1$. Summing \eqref{eq:jointJK} over $K>k$ gives
\begin{equation}\label{eq:K-survival}
 \Pp(K>k\mid\wb,\xib)=R_k+\xi_{k+1}s_k.
\end{equation}
Thus, when the weights and schedule values can be generated
successively, a $K$-first sampler locates $K$ by comparing
the survival probabilities in \eqref{eq:K-survival} with
a single uniform draw, accounting for the probability mass
of all larger indices throughout the search. Once $K$ is
determined, normalizing $a_1,\ldots,a_K$ gives the
probabilities for selecting $J$.

\subsection{Recursive reweighting and sequential generation}
Passing from $G_{k-1}^\star$ to $G_k^\star$ adds $\theta_k$ with share $\lambda_k=a_k/s_k$ and preserves the relative proportions of earlier atoms. Assume $w_1>0$, so $\lambda_1=1$, and write $Ph=\int h\dd P$ and $\TV(P,Q)=\sup_{A\in\mathcal B}|P(A)-Q(A)|$.

\begin{proposition}[One-atom recursion]\label{prop:oneatom}
For $k\geq2$,
\begin{equation*}
G_k^\star=(1-\lambda_k)G_{k-1}^\star+\lambda_k\delta_{\theta_k}.
\end{equation*}
Consequently, for every bounded measurable $h$,
\begin{equation*}
G_k^\star h=(1-\lambda_k)G_{k-1}^\star h+\lambda_kh(\theta_k).
\end{equation*}
Furthermore,
\begin{equation*}
\widetilde w_j^{(k)}
=\lambda_j\prod_{\ell=j+1}^k(1-\lambda_\ell),
\qquad 1\leq j\leq k.
\end{equation*}
If the atoms are distinct, as they are almost surely under diffuse $G_0$, then
\begin{equation*}
\TV(G_k^\star,G_{k-1}^\star)=\lambda_k;
\end{equation*}
without distinctness the left-hand side is at most $\lambda_k$.
\end{proposition}
The proof is given in Supplement~B.4.

The recursion updates event probabilities and finite-mixture evaluations, while the reverse product records the rescaling of each retained atom. Under proper stick-breaking with $0<v_j<1$, the tail-mass schedule
\begin{equation}
 \xi_j=R_{j-1}
\end{equation}
gives $a_j=v_j$ and $p_k=w_k\sum_{j\leq k}v_j$. When the stick-breaking fractions are constant, the finite allocation probabilities are uniform, recovering Example~\ref{ex:count-nb} at $r=2$. Further details on sequential generation appear in Supplement~B.4, while Supplement~C.1 uses the reverse product to obtain a parametrization for posterior inference.

\subsection{Exact averaging and a fixed representation level}
Although \eqref{eq:barycenter} recovers $G$ by averaging, increasing one fixed level also reweights every retained atom. To separate this redistribution from omitted mass, for $s_k>0$ put
\[
 G_{\mathrm{ren}}^{(k)}=\sum_{j\leq k}\frac{w_j}{1-R_k}\delta_{\theta_j},\qquad
 M_k=\xi_1/\xi_k,\qquad D_k=\frac{\sqrt{M_k}-1}{\sqrt{M_k}+1}.
\]
Lemma~B.1 in Supplement~B.2 gives
\begin{equation*}
 \TV(G_{\mathrm{ren}}^{(k)},G_k^\star)\leq D_k,\qquad
 \TV(G,G_k^\star)\leq R_k+(1-R_k)D_k.
\end{equation*}
The two terms separate omitted mass from redistribution among retained atoms. Expected bounds are given in Supplement~B.2, and index tails and moments in Supplement~B.3.

\begin{theorem}[Total-variation asymptotics of the finite measures]\label{thm:fixed-limit}
Assume the atoms are distinct and write $s_\infty=\sum_{j\geq1}w_j/\xi_j\in(0,\infty]$.
If $s_\infty<\infty$, define
\begin{equation*}
G_\xi=\sum_{j\geq1}\frac{w_j/\xi_j}{s_\infty}\delta_{\theta_j}.
\end{equation*}
Then
\begin{equation*}
\TV(G_k^\star,G_\xi)=1-\frac{s_k}{s_\infty}\longrightarrow0
\end{equation*}
for $s_k>0$. Moreover $G_\xi=G$ if and only if $1/\xi_j$ is constant on the positive-weight support. If $s_\infty=\infty$, then
\begin{equation*}
\TV(G_k^\star,G)\longrightarrow1.
\end{equation*}
\end{theorem}
The proof is given in Supplement~B.4.
For the DP and geometric weights under the tail-mass schedule, $a_j=v_j$ and $s_\infty=\infty$ almost surely. Thus these finite measures approach total-variation distance one from $G$, although their average still recovers it exactly.

\subsection{Posterior augmentation and full conditional updates}
For observations $x_1,\ldots,x_n$ from \eqref{eq:ssp-mixture},
each auxiliary index $k_i\geq z_i$ restricts the update of $z_i$
to labels $1,\ldots,k_i$, including currently unoccupied ones.
Thus labels up to $k^\star=\max_i k_i$ suffice for all allocation
updates, while the largest occupied label $z^\star=\max_i z_i$
satisfies $c_n\leq z^\star\leq k^\star$. With $\wb\sim\Pi$ and conditional schedule law
$Q_\xi(\dd\xib\mid\wb)$, combining \eqref{eq:jointJK} with
\eqref{eq:ssp-mixture} gives the augmented likelihood
\begin{equation}\label{eq:joint-posterior}
 \mathcal L(\boldsymbol x,\zb,\kb\mid\wb,\thetab,\xib)
 =\prod_{i=1}^n
 \frac{w_{z_i}}{\xi_{z_i}}(\xi_{k_i}-\xi_{k_i+1})
 f(x_i\mid\theta_{z_i})\mathbf1\{z_i\leq k_i\}.
\end{equation}
The joint posterior is proportional to this likelihood times
$\Pi(\dd\wb)Q_\xi(\dd\xib\mid\wb)\prod_jG_0(\dd\theta_j)$,
yielding the full conditional updates
\begin{align}
 \Pp(z_i=j\mid k_i,\cdots)
 &\propto\frac{w_j}{\xi_j}f(x_i\mid\theta_j),
 &&1\leq j\leq k_i,\\
 \Pp(k_i=k\mid z_i,\cdots)
 &=\frac{\xi_k-\xi_{k+1}}{\xi_{z_i}}
   \mathbf1\{k\geq z_i\}.\label{eq:k-update}
\end{align}
Each component parameter $\theta_j$ has conditional law
proportional to
$G_0(\dd\theta_j)\prod_{i:z_i=j}f(x_i\mid\theta_j)$.
Supplement~C.1 derives these conditionals, while Supplement~C.3 derives the component
parameter updates and posterior predictive density for
Gaussian mixtures.

A fixed exogenous schedule preserves the allocation-based weight
update, whereas the tail-mass schedule changes it without altering
the  SSP model. Supplement~C.1 covers
independent-beta updates and the geometric case with a
shared stick-breaking fraction.

\subsection{Slice augmentation and the number of available labels}
For a fixed deterministic schedule, the finite index can be viewed as recording the interval containing the continuous slice variable of \citet{Kalli2011}. Integrating over that interval recovers the same discrete joint law. Supplement~C.2 gives the precise correspondence.

Since  indices \eqref{eq:k-update} are conditionally
independent, the distribution of their maximum $k^\star$
can be obtained explicitly. In particular, a deterministic
geometric schedule $\xi_j=cq^{j-1}$ gives
\[
 \Pp(k^\star-z^\star\geq t\mid\zb)\leq nq^t,
 \qquad t=0,1,\ldots.
\]
For fixed $q$, the extension beyond the largest occupied
label therefore grows at most logarithmically with $n$,
in probability, while the total excess over occupancy also
includes the $z^\star-c_n$ unoccupied labels below $z^\star$.
Supplement~C.2 derives the exact distribution
and a high-probability bound for $k^\star$, while
Supplement~C.5 compares these
quantities with the DP slice overhead studied by
\citet{FranzoliniGaffi2026}.

\section{Numerical considerations}
We examine how the index law and allocation rule shape
prior-predictive behavior, first comparing constructions
with a common index law and then priors calibrated to the
same clustering characteristic. The invariant family
extends these comparisons to prediction conditional on
particular partitions. Numerical illustrations of
fixed-level approximation and posterior representation
size appear in Supplements~B.4
and C.5, respectively.
The displayed calculations use exact formulas or finite
sums with controlled remainders.

\par
\begin{figure}[tbp]\centering
\includegraphics[width=.88\linewidth]{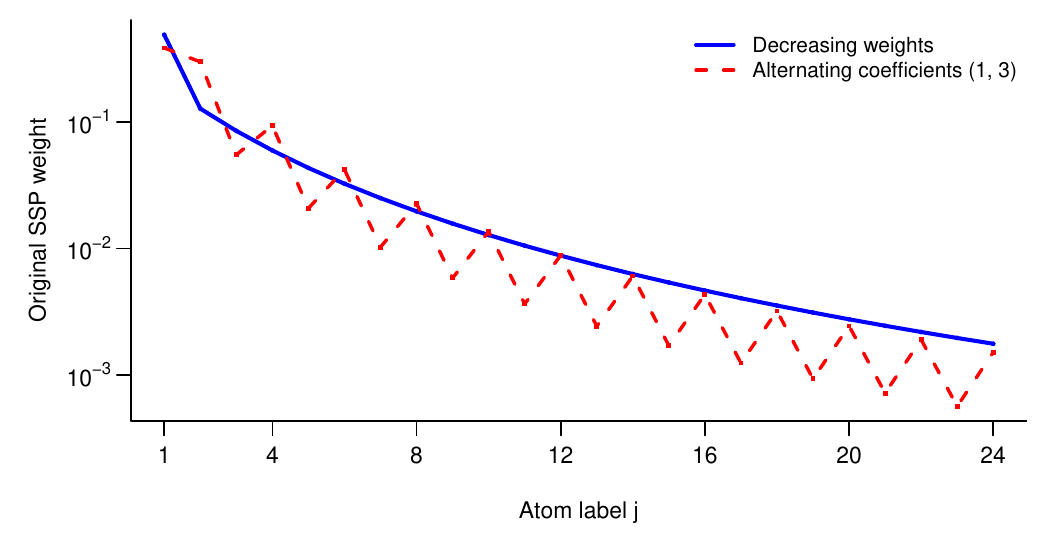}
\caption{Under the same beta-geometric index law,
with $(\rho,b)=(2,5)$ and $\E(K)=6$, the construction in
\eqref{eq:bg-weights} gives decreasing weights, while
alternating coefficients $(1,3,1,3,\ldots)$ give
nonmonotone weights. Weights are shown on a logarithmic
scale.}
\label{fig:weight-geometry}
\end{figure}

\medskip\noindent\textit{Different weights under the same index law.}\enspace
For the beta-geometric law \eqref{eq:bg-index} with
$(\rho,b)=(2,5)$, the constructions in
Proposition~\ref{prop:bg-ssp} and
Example~\ref{ex:alternating-ssp} both have $\E(K)=6$,
but their different second-stage allocation rules yield
different SSP weights. Figure~\ref{fig:weight-geometry}
compares the decreasing construction with the alternating
construction at $\gamma=3$, whose weights repeatedly
increase between adjacent labels despite a decreasing
schedule. Their tie probabilities $\sum_jw_j^2$ are
approximately $0.2717$ and $0.2528$, respectively.
Since this sum is invariant under relabeling, the two
weight sequences cannot be permutations of one another.

Figure~S.1 in
Supplement~A.3 illustrates how the
largest weight in each family decreases as $b$ grows.
A prior on $b$ with unbounded support above and finite
mean then gives full weak support and finite $\E(K)$.

\par\medskip\noindent\textit{Matching pairwise coincidence.}\enspace
Priors with the same pairwise coincidence probability can imply
different occupancy and discovery behavior. Table~\ref{tab:predictive-calibration}
compares a DP with three finite-index constructions, all calibrated
to $T_2=0.2$, through their expected occupancy, singleton counts and
discovery probability at $n=50$. The matching concentration is
$\alpha=4$ for the DP, while the other parameters are obtained by
solving the coincidence equations for the constructions in
Examples~\ref{ex:count-nb} and \ref{ex:count-poisson}.
In mixture models, these summaries concern the partition induced
by the latent allocations.

\begin{table}[tbp]\centering\small
\caption{Prior-predictive summaries for priors with coincidence
probability $T_2=0.2$. Here $M_{1,50}$ counts singleton blocks,
and $d_{50}$ is the probability that latent draw~$51$ takes a
previously unobserved value, averaged over possible samples of
size~$50$. Parameters are fixed in each row.}
\label{tab:predictive-calibration}
\vspace{4pt}
\begin{tabular}{@{}llrrr@{}}
\toprule
Prior & Matching parameter & $\E(c_{50})$ & $\E(M_{1,50})$ & $d_{50}$\\
\midrule
DP & $\alpha=4$ & 10.894 & 3.774 & 0.07407 \\
Geometric & $p=1/3$ & 8.887 & 2.466 & 0.04836 \\
Negative-binomial, $r=3$ & $p=0.524192$ & 7.594 & 1.639 & 0.03211 \\
Shifted-Poisson & $\nu=3.524901$ & 7.012 & 1.170 & 0.02286 \\
\bottomrule
\end{tabular}
\end{table}

The differences across rows show that matching coincidence
alone leaves other aspects of the partition law undetermined.
Assessing occupancy, singleton counts and discovery together
can therefore guide the choice of index and allocation
parameters, or the specification of a joint prior on them.
Supplement~A.6 gives the formulas,
numerical bounds and reproducible calculations.

\par\medskip\noindent\textit{Paired weights, prediction and index size.}\enspace
Table~\ref{tab:paired-prediction} compares prediction under the
paired prior as $\alpha$ varies, with $q=1/2$. A DP with
concentration $3$ provides a comparison with the same coincidence
probability of $1/4$ as the paired prior at $\alpha=1$.
The discovery probabilities obtained from EPPF ratios agree
with those calculated by averaging \eqref{eq:paired-predictive}
over the possible pairings. Although the partitions $(5,1)$
and $(3,3)$ have the same sample size and occupancy, their
block sizes affect the posterior pairing probabilities and
hence discovery under the paired prior. The DP gives the same
discovery probability for both partitions.

\begin{table}[tbp]\centering\small
\caption{Discovery probabilities conditional on the displayed
block sizes and expected representation index, with $q=1/2$.
The paired prior at $\alpha=1$ and the DP with concentration
$3$ both have prior coincidence probability $1/4$.}
\label{tab:paired-prediction}
\begin{tabular}{@{}lrrrr@{}}
\toprule
Prior & \shortstack{Tie\\probability}
& \shortstack{New species\\after $(5,1)$}
& \shortstack{New species\\after $(3,3)$} & $\E(K)$\\
\midrule
Paired, $\alpha=0.5$ & 0.333333 & 0.153846 & 0.091811 & 3.500000\\
Paired, $\alpha=1$ & 0.250000 & 0.265306 & 0.169643 & 4.500000\\
Paired, $\alpha=2$ & 0.166667 & 0.416667 & 0.294118 & 6.500000\\
Paired, $\alpha=4$ & 0.100000 & 0.584615 & 0.463158 & 10.500000\\
DP, concentration $3$ & 0.250000 & 0.333333 & 0.333333 & 5.000000\\
\bottomrule
\end{tabular}
\end{table}

The distinction also appears after averaging over possible
samples. At $n=50$, the paired prior with $\alpha=1$ gives
$\E(c_{50})=7.6121$ and $d_{50}=0.0392$, compared with
$9.1141$ and $0.0566$, respectively, for the DP.
Both priors have full weak support, yet encode different
beliefs about occupancy and discovery despite their matched
coincidence probability. The finite-sum calculations and
discovery curves appear in
Supplement~C.4.
Changing $q$ preserves the prior and these predictions while
varying the geometric contribution $q/(1-q)$ to $\E(K)$.

\Needspace{5\baselineskip}
\section{Discussion}
The finite-index identity represents existing SSPs and constructs new ones whose finite-sample implications can guide prior specification. Coincidence, occupancy and discovery calculations permit assessment without evaluating infinite weights, while a shared hyperprior expresses uncertainty about the inputs and is updated jointly with the latent allocations. This connects Bayesian elicitation, assessment and updating within a coherent random-measure model.

Further work on elicitation should clarify which prior-predictive
constraints are attainable and how to translate them into prior
specifications. Full weak support complements this goal, although
it does not determine the probabilities of sample configurations.
The survival construction achieves full weak support through
hyperpriors, whereas the paired prior retains it at every fixed
positive concentration despite imposing equal weights within each
pair. This illustrates how full weak support can coexist with a
strong belief about relative prevalence. Beyond this support
property, the paired law combines a directly generated size-biased
order, partition prediction and a finite-index representation
outside the class of homogeneous normalized CRMs. Establishing
posterior consistency requires additional arguments tailored to
the likelihood.

For a fixed SSP, changing the schedule affects representation size but preserves the marginal posterior. Its effect on mixing efficiency requires separate assessment. The framework thus combines prior construction with exact finite sampling and posterior augmentation, without identifying a finite conditional measure with a truncation approximation.

\ifnum\anon=1
\section*{Funding}

Rams\'es H. Mena acknowledges support from
UNAM--PAPIIT through project IT100524.
\fi

\section*{Disclosure Statement}

The authors report no competing interests to declare.

\section*{Declaration of Generative AI Use}

During the preparation of this manuscript, the authors used Overleaf's integrated AI Assistant to generate LaTeX table structures and troubleshoot compilation errors. After using this tool, the authors reviewed and edited the content as needed and take full responsibility for the content of the publication.
\bibliography{references_sspfin}

@article{Arbel2019,
  author = {Arbel, Julyan and De Blasi, Pierpaolo and Pr\"unster, Igor},
  year = {2019},
  title = {Stochastic approximations to the Pitman--Yor process},
  journal = {Bayesian Analysis},
  volume = {14},
  number = {4},
  pages = {1201--1219},
  doi = {10.1214/18-BA1127},
  url = {https://arxiv.org/abs/1806.10867},
}

@article{BlackwellMacQueen1973,
  year = {1973},
  author = {Blackwell, David and MacQueen, James B.},
  title = {Ferguson distributions via P\'olya urn schemes},
  journal = {The Annals of Statistics},
  volume = {1},
  number = {2},
  pages = {353--355},
  doi = {10.1214/aos/1176342372},
}

@article{CanaleCorradinNipoti2022,
  author  = {Canale, Antonio and Corradin, Riccardo and Nipoti, Bernardo},
  title   = {Importance conditional sampling for {Pitman--Yor} mixtures},
  journal = {Statistics and Computing},
  year    = {2022},
  volume  = {32},
  pages   = {40},
  doi     = {10.1007/s11222-022-10096-0}
}

@article{DeBlasi2015,
  author = {De Blasi, Pierpaolo and Favaro, Stefano and Lijoi, Antonio and Mena, Rams\'es H. and Pr\"unster, Igor and Ruggiero, Matteo},
  year = {2015},
  title = {Are Gibbs-Type Priors the Most Natural Generalization of the Dirichlet Process?},
  journal = {IEEE Transactions on Pattern Analysis and Machine Intelligence},
  volume = {37},
  number = {2},
  pages = {212--229},
  doi = {10.1109/TPAMI.2013.217},
  url = {https://doi.org/10.1109/TPAMI.2013.217},
}

@article{DeBlasiMartinezMenaPrunster2020,
  author = {De Blasi, Pierpaolo and Mart\'inez, Asael F. and Mena, Rams\'es H. and Pr\"unster, Igor},
  title = {On the inferential implications of decreasing weight structures in mixture models},
  year = {2020},
  journal = {Computational Statistics \& Data Analysis},
  volume = {147},
  pages = {106940},
  doi = {10.1016/j.csda.2020.106940},
  url = {https://www.sciencedirect.com/science/article/abs/pii/S0167947320300311}
}

@article{Ferguson1973,
  year = {1973},
  author = {Ferguson, Thomas S.},
  title = {A Bayesian analysis of some nonparametric problems},
  journal = {The Annals of Statistics},
  volume = {1},
  number = {2},
  pages = {209--230},
  doi = {10.1214/aos/1176342360},
}

@article{FuentesGarcia2010,
  year = {2010},
  author = {Fuentes-Garc\'ia, R. and Mena, Rams\'es H. and Walker, Stephen G.},
  title = {A new Bayesian nonparametric mixture model},
  journal = {Communications in Statistics---Simulation and Computation},
  volume = {39},
  number = {4},
  pages = {669--682},
  doi = {10.1080/03610910903580963},
}

@article{GilLeyvaMena2023,
  author = {Gil-Leyva, Mar\'ia F. and Mena, Rams\'es H.},
  year = {2023},
  title = {Stick-Breaking Processes With Exchangeable Length Variables},
  journal = {Journal of the American Statistical Association},
  volume = {118},
  number = {541},
  pages = {537--550},
  doi = {10.1080/01621459.2021.1941054},
  url = {https://arxiv.org/abs/2008.04475},
}

@article{GilLeyvaMenaNicoleris2020,
  author = {M. F. Gil-Leyva and Rams\'es H. Mena and Theodoros Nicoleris},
  title = {Beta-binomial stick-breaking non-parametric prior},
  year = {2020},
  journal = {Electronic Journal of Statistics},
  volume = {14},
  number = {1},
  pages = {1479--1507},
  doi = {10.1214/20-EJS1694},
  url = {https://arxiv.org/abs/1908.06602}
}

@article{GnedinPitman2006,
  author = {Gnedin, Alexander and Pitman, Jim},
  year = {2006},
  title = {Exchangeable Gibbs partitions and Stirling triangles},
  journal = {Journal of Mathematical Sciences},
  volume = {138},
  number = {3},
  pages = {5674--5685},
  doi = {10.1007/s10958-006-0335-z},
  url = {https://arxiv.org/abs/math/0412494},
}

@article{HatjispyrosMerkatasWalker2023,
  author = {Spyridon J. Hatjispyros and Christos Merkatas and Stephen G. Walker},
  title = {Mixture models with decreasing weights},
  year = {2023},
  journal = {Computational Statistics \& Data Analysis},
  volume = {179},
  pages = {107651},
  doi = {10.1016/j.csda.2022.107651},
  url = {https://research.aalto.fi/en/publications/mixture-models-with-decreasing-weights/}
}

@article{Kalli2011,
  author = {Kalli, Maria and Griffin, Jim E. and Walker, Stephen G.},
  title = {Slice sampling mixture models},
  journal = {Statistics and Computing},
  year = {2011},
  volume = {21},
  number = {1},
  pages = {93--105},
  doi = {10.1007/s11222-009-9150-y}
}

@article{Lee2013,
  year = {2013},
  author = {Lee, Jaeyong and Quintana, Fernando A. and M\"uller, Peter and Trippa, Lorenzo},
  title = {Defining predictive probability functions for species sampling models},
  journal = {Statistical Science},
  volume = {28},
  number = {2},
  pages = {209--222},
  doi = {10.1214/12-STS407},
}

@incollection{LijoiPrunster2010,
  year = {2010},
  author = {Lijoi, Antonio and Pr\"unster, Igor},
  title = {Models beyond the Dirichlet process},
  booktitle = {Bayesian Nonparametrics},
  editor = {Hjort, Nils Lid and Holmes, Chris and M\"uller, Peter and Walker, Stephen G.},
  pages = {80--136},
  publisher = {Cambridge University Press},
  doi = {10.1017/CBO9780511802478.004},
}

@article{McCullaghYang2008,
  author = {Peter McCullagh and Jie Yang},
  title = {How many clusters?},
  year = {2008},
  journal = {Bayesian Analysis},
  volume = {3},
  number = {1},
  pages = {101--120},
  doi = {10.1214/08-BA304},
  url = {https://doi.org/10.1214/08-BA304}
}

@article{Nobile2004,
  author = {Agostino Nobile},
  title = {On the posterior distribution of the number of components in a finite mixture},
  year = {2004},
  journal = {The Annals of Statistics},
  volume = {32},
  number = {5},
  pages = {2044--2073},
  doi = {10.1214/009053604000000788},
  url = {https://arxiv.org/pdf/math/0503673}
}

@article{PapaspiliopoulosRoberts2008,
  author = {Papaspiliopoulos, Omiros and Roberts, Gareth O.},
  year = {2008},
  title = {Retrospective Markov chain Monte Carlo methods for Dirichlet process hierarchical models},
  journal = {Biometrika},
  volume = {95},
  number = {1},
  pages = {169--186},
  doi = {10.1093/biomet/asm086},
  url = {https://doi.org/10.1093/biomet/asm086},
}

@incollection{Pitman1996,
  year = {1996},
  author = {Pitman, Jim},
  title = {Some developments of the Blackwell--MacQueen urn scheme},
  booktitle = {Statistics, Probability and Game Theory},
  series = {IMS Lecture Notes--Monograph Series},
  volume = {30},
  publisher = {Institute of Mathematical Statistics},
  doi = {10.1214/lnms/1215453576},
}

@book{Pitman2006,
  year = {2006},
  author = {Pitman, Jim},
  title = {Combinatorial Stochastic Processes},
  series = {Lecture Notes in Mathematics},
  volume = {1875},
  publisher = {Springer},
  doi = {10.1007/b116015},
}

@article{PitmanYor1997,
  year = {1997},
  author = {Pitman, Jim and Yor, Marc},
  title = {The two-parameter Poisson--Dirichlet distribution derived from a stable subordinator},
  journal = {The Annals of Probability},
  volume = {25},
  number = {2},
  pages = {855--900},
  doi = {10.1214/aop/1024404422},
}

@article{Sethuraman1994,
  year = {1994},
  author = {Sethuraman, Jayaram},
  title = {A constructive definition of Dirichlet priors},
  journal = {Statistica Sinica},
  volume = {4},
  number = {2},
  pages = {639--650},
}

@article{Walker2007,
  author = {Walker, Stephen G.},
  year = {2007},
  title = {Sampling the Dirichlet mixture model with slices},
  journal = {Communications in Statistics---Simulation and Computation},
  volume = {36},
  number = {1},
  pages = {45--54},
  doi = {10.1080/03610910601096262},
  url = {https://doi.org/10.1080/03610910601096262},
}

@article{BissiriOngaro2014,
 author={Bissiri, Pier Giovanni and Ongaro, Andrea}, year={2014},
 title={On the topological support of species sampling priors},
 journal={Electronic Journal of Statistics}, volume={8}, number={1}, pages={861--882}, doi={10.1214/14-EJS912}}

@article{BuddanaKozubowski2014,
 author={Buddana, Amrutha and Kozubowski, Tomasz J.}, year={2014},
 title={Discrete Pareto distributions}, journal={Economic Quality Control},
 volume={29}, number={2}, pages={143--156}, doi={10.1515/eqc-2014-0014}}

@article{PitmanISBP1996,
 author={Pitman, Jim},title={Random discrete distributions invariant under size-biased permutation},journal={Advances in Applied Probability},year={1996},volume={28},number={2},pages={525--539},doi={10.2307/1428070}}

@article{Kingman1967,
 author={Kingman, J. F. C.},title={Completely random measures},journal={Pacific Journal of Mathematics},year={1967},volume={21},number={1},pages={59--78},doi={10.2140/pjm.1967.21.59}}

@article{JamesLijoiPrunster2009,
 author={James, Lancelot F. and Lijoi, Antonio and Pr{\"u}nster, Igor},title={Posterior analysis for normalized random measures with independent increments},journal={Scandinavian Journal of Statistics},year={2009},volume={36},number={1},pages={76--97},doi={10.1111/j.1467-9469.2008.00609.x}}

@article{PermanPitmanYor1992,
 author={Perman, Mihael and Pitman, Jim and Yor, Marc},
 title={Size-biased sampling of {Poisson} point processes and excursions},
 journal={Probability Theory and Related Fields},year={1992},volume={92},pages={21--39}}

@article{GilLeyvaLijoiMenaPrunster2026,
  author = {Gil-Leyva, Mar\'ia F. and Lijoi, Antonio and Mena, Rams\'es H. and Pr\"unster, Igor},
  title = {Markov stick-breaking processes},
  journal = {The Annals of Statistics},
  year = {2026},
  volume = {54},
  number = {3},
  pages = {1345--1371},
  doi = {10.1214/25-AOS2607}
}

@article{FavaroEtAl2016,
  author = {Favaro, Stefano and Lijoi, Antonio and Nava, Consuelo and Nipoti, Bernardo and Pr\"unster, Igor and Teh, Yee Whye},
  title = {On the stick-breaking representation for homogeneous {NRMIs}},
  journal = {Bayesian Analysis},
  year = {2016},
  volume = {11},
  number = {3},
  pages = {697--724},
  doi = {10.1214/15-BA964}
}

@inproceedings{FranzoliniGaffi2026,
  author = {Franzolini, Beatrice and Gaffi, Francesco},
  title = {Complexity Bounds for {Dirichlet} Process Slice Samplers},
  booktitle = {Proceedings of the 43rd International Conference on Machine Learning},
  series = {Proceedings of Machine Learning Research},
  volume = {306},
  pages = {31553--31574},
  year = {2026},
  publisher = {PMLR},
  url = {https://proceedings.mlr.press/v306/franzolini26a.html}
}

@misc{MenaSSPSupp2026,
 author = {Mena, Rams{\'e}s H. and Merkatas, Christos and Nicoleris, Theodoros and Rodr{\'i}guez, Carlos E.},
 year = {2026},
 title = {Supplement to ``Species sampling priors through exact augmentation''},
 note = {Accompanying supplementary material}
}

@article{Binette2019,
  author = {Binette, Olivier},
  title = {A Note on Reverse {Pinsker} Inequalities},
  journal = {IEEE Transactions on Information Theory},
  year = {2019},
  volume = {65},
  number = {7},
  pages = {4094--4096},
  doi = {10.1109/TIT.2019.2896192}
}

@article{NavarreteQuintanaMueller2008,
  author = {Navarrete, Carlos and Quintana, Fernando A. and M{\"u}ller, Peter},
  title = {Some issues in nonparametric {Bayesian} modeling using species sampling models},
  journal = {Statistical Modelling},
  year = {2008},
  volume = {8},
  number = {1},
  pages = {3--21},
  doi = {10.1177/1471082X0700800102}
}

@book{LastPenrose2017,
  author = {Last, G{\"u}nter and Penrose, Mathew},
  title = {Lectures on the {Poisson} Process},
  year = {2017},
  series = {Institute of Mathematical Statistics Textbooks},
  volume = {7},
  publisher = {Cambridge University Press},
  address = {Cambridge},
  doi = {10.1017/9781316104477}
}
\end{document}


\spacingset{1}

\ifnum\anon=1

\title{\bfseries
Supplement to ``Species sampling priors
through exact augmentation''}

\author{
Rams\'es H. Mena$^{1}$\,\orcidlink{0000-0001-9608-8059},
Christos Merkatas$^{2}$\,\orcidlink{0000-0001-7880-7443},\\
Theodoros Nicoleris$^{3}$,
Carlos E. Rodr\'iguez$^{1}$\,\orcidlink{0000-0002-4796-0784}\thanks{Corresponding author:
\href{mailto:carloserwin@sigma.iimas.unam.mx}
{carloserwin@sigma.iimas.unam.mx}}
}

\date{}
\maketitle

\vspace{-1.5em}

\begin{center}
\small
$^{1}$ Department of Probability and Statistics,
IIMAS, Universidad Nacional Aut\'onoma de M\'exico\\[0.3em]
$^{2}$ Department of Statistics and
Actuarial--Financial Mathematics,
University of the Aegean\\[0.3em]
$^{3}$ Department of Economics,
National and Kapodistrian University of Athens
\end{center}

\else

\bigskip
\bigskip
\bigskip
\begin{center}
{\LARGE\bfseries
Supplement to ``Species sampling priors
through exact augmentation''
}
\end{center}
\medskip

\fi

\bigskip

\begin{abstract}
This supplement supplies the proofs and additional calculations for the main paper. Supplement A develops compatibility, inverse maps, further index distributions and support results, including the paired-weight prior and its invariant defining order. It also treats survival-based hierarchies and finite-index prior-predictive assessment. Supplement B studies the sample law, recursive sampling and total-variation comparisons. Supplement C derives posterior updates, the matched slice comparison, Gaussian-mixture calculations and the paired prior's partition and predictive laws, together with the numerical evaluation procedures.
\end{abstract}

\vfill
\newpage

\spacingset{1.8}

Equation numbers beginning with S refer to this supplement.
Other equation numbers refer to the main paper.
Results numbered by a lettered section are local to this
supplement, whereas results with numeric section numbers
are in the main paper.

\appendix

\renewcommand{\appendixname}{Supplement}
\renewcommand{\theequation}{S.\arabic{equation}}
\renewcommand{\thefigure}{S.\arabic{figure}}
\renewcommand{\thetable}{S.\arabic{table}}

\section{Prior construction and support}\label{app:prior}
This appendix supplies the compatibility arguments used to choose between specifying a prior and representing an existing one. It then develops the additional index distributions and proves the support results for the resulting SSPs.

\subsection{Compatibility and inverse maps}\label{sec:inverse-proofs}
We first prove Theorem~3.1, then establish
compatibility criteria for arbitrary finite allocation rules
and for the common-coefficient representation, before
deriving the beta-geometric weight family.

\par\medskip\Needspace{3\baselineskip}\noindent\textit{Proof of Theorem~3.1.}\label{sec:construction-proof}\enspace
\begin{proof}
The sums defining $\xi_j$ are finite as
\[
0<\xi_j\leq\frac{1}{s_j}\sum_{k\geq j}\pi_k
\leq\frac{1}{a_1}\sum_{k\geq j}\pi_k\longrightarrow0.
\]
Their successive differences are $\xi_k-\xi_{k+1}=\pi_k/s_k>0$. Tonelli's theorem gives
\[
\sum_{j\geq1}w_j
=\sum_{j\geq1}a_j\sum_{k\geq j}\frac{\pi_k}{s_k}
=\sum_{k\geq1}\frac{\pi_k}{s_k}\sum_{j=1}^k a_j
=1.
\]
The same interchange, with $\delta_{\theta_j}(A)$ included, proves (9) for all measurable $A$. Since $w_j/\xi_j=a_j$, the successive-difference identity proves (8). For the converse,
\[
\sum_{k\geq j}\frac{p_k(\wb,\xib)}{s_k}
=\sum_{k\geq j}(\xi_k-\xi_{k+1})=\xi_j.
\]
The weights are positive and sum to one for every admissible input pair. Drawing the atoms i.i.d. from $G_0$, independently of $(\boldsymbol\pi,\boldsymbol a)$, makes them independent of these induced weights and completes the SSP specification.
\end{proof}

\begin{lemma}[Feasibility of a prescribed index law]\label{lem:count-order}
There are probabilities $b_j^{(k)}\geq0$, $1\leq j\leq k$, with $\sum_{j=1}^k b_j^{(k)}=1$ and
\begin{equation}\label{eq:count-triangular}
w_j=\sum_{k\geq j}\pi_k b_j^{(k)}
\end{equation}
if and only if $Q_k\leq W_k$ for every $k$.
\end{lemma}
\begin{proof}
The joint probabilities $\Pp(J=j,K=k)=\pi_kb_j^{(k)}$ have $J\leq K$, so $\{K\leq k\}\subseteq\{J\leq k\}$ and $Q_k\leq W_k$. Conversely, take a common $U\sim\mathrm{Unif}(0,1)$ and set
\[
J=\min\{j:W_j\geq U\}\qquad\mbox{and}\quad
K=\min\{k:Q_k\geq U\}.
\]
Both variables are finite almost surely. The cumulative inequality gives $J\leq K$, with marginals $\wb$ and $\boldsymbol\pi$. Taking $b_j^{(k)}=\Pp(J=j\mid K=k)$ establishes \eqref{eq:count-triangular}.
\end{proof}

\par\medskip\Needspace{3\baselineskip}\noindent\textit{Proof of Proposition~3.2.}\enspace
\begin{proof}
Summing the probabilities in (3) through $k$ and collecting the contribution of each selected label gives
\begin{equation*}
 \sum_{\ell=1}^k p_\ell
 =\sum_{j=1}^k\frac{w_j}{\xi_j}
       \sum_{\ell=j}^k(\xi_\ell-\xi_{\ell+1})
 =W_k-\xi_{k+1}s_k.
\end{equation*}
If $p_\ell=\pi_\ell$ for every $\ell$, then
$W_k-Q_k=\xi_{k+1}s_k>0$. At level $k-1$,
$W_{k-1}-Q_{k-1}=\xi_ks_{k-1}$, and adding $w_k$ together with
$\xi_ks_k=\xi_ks_{k-1}+w_k$ gives $W_k-Q_{k-1}=\xi_ks_k$. Taking the ratio of these identities eliminates $s_k$ and yields (15). Since its right-hand side is determined entirely by the prescribed cumulative sums, fixing $\xi_1$ determines the remaining schedule recursively and proves uniqueness up to a common positive scale.

Conversely, suppose (14) holds, put $\xi_1=1$, and define the remaining schedule values by (15). Since $\pi_k>0$, each ratio is strictly between zero and one. Define $s_k=\sum_{j=1}^k w_j/\xi_j$. We prove successively that
\[
 \xi_ks_k=W_k-Q_{k-1},\qquad
 \xi_{k+1}s_k=W_k-Q_k.
\]
For $k=1$, the first equality is $\xi_1s_1=w_1=W_1-Q_0$. Multiplying it by the schedule ratio gives the second equality at $k=1$. If the second equality holds at a general $k$, then
\[
 \xi_{k+1}s_{k+1}
 =\xi_{k+1}s_k+w_{k+1}
 =W_k-Q_k+w_{k+1}
 =W_{k+1}-Q_k.
\]
Together, these steps establish both identities at $k+1$,
completing the induction. Taking their difference then gives
\[
 (\xi_k-\xi_{k+1})s_k=Q_k-Q_{k-1}=\pi_k.
\]
Finally, $s_k\geq w_1>0$ and both cumulative sums tend to one, so
\[
 0<\xi_{k+1}=\frac{W_k-Q_k}{s_k}
 \leq\frac{W_k-Q_k}{w_1}\longrightarrow0.
\]
The constructed schedule therefore satisfies all the required properties.
\end{proof}

\par\medskip\Needspace{3\baselineskip}\noindent\textit{Proof of Corollary~3.3.}\label{sec:marginal-proof}\enspace
\begin{proof}
For a deterministic schedule, taking expectation in (3) gives $\Pp(K=k)=(\xi_k-\xi_{k+1})\overline s_k$. The positive mean weights sum to one by Tonelli's theorem. Proposition~3.2 therefore applies to $(\overline w_j)$ and gives both the cumulative condition and the unique schedule up to scale. Bayes' theorem yields
\[
 \Pi(\dd\wb\mid K=k)
 =\frac{p_k(\wb,\xib)}{\pi_k}\Pi(\dd\wb)
 =\frac{(\xi_k-\xi_{k+1})s_k(\wb)}
        {(\xi_k-\xi_{k+1})\overline s_k}\Pi(\dd\wb),
\]
which is (16). Averaging these conditional laws with probabilities $\pi_k$ returns $\Pi$.
\end{proof}

\par\medskip\Needspace{3\baselineskip}\noindent\textit{Proof of Proposition~3.4.}\label{sec:inverse-weights-proof}\enspace
\begin{proof}
The identity $p_k=(\xi_k-\xi_{k+1})s_k=\pi_k$ forces $s_k=t_k$. Since $s_j-s_{j-1}=w_j/\xi_j$, the differences of $t_j$ must be nonnegative and determine the weights uniquely. Conversely, suppose those differences are nonnegative and set $w_j=\xi_j(t_j-t_{j-1})$. Writing $\xi_j=\sum_{k\geq j}(\xi_k-\xi_{k+1})$ and applying Tonelli's theorem gives
\[
 \sum_{j\geq1}w_j
 =\sum_{k\geq1}(\xi_k-\xi_{k+1})
       \sum_{j=1}^k(t_j-t_{j-1})
 =\sum_{k\geq1}(\xi_k-\xi_{k+1})t_k
 =\sum_{k\geq1}\pi_k=1.
\]
The partial sums $\sum_{j\leq k}w_j/\xi_j$ equal $t_k$, so (3) gives the required probabilities. The weights are strictly positive exactly when all the differences are positive. Finally, $t_k\geq t_1$ gives
$H_j=\sum_{k\geq j}(\xi_k-\xi_{k+1})t_k\geq t_1\xi_j$,
which proves (17).
\end{proof}

Equivalently,
\begin{equation*}
 \begin{aligned}
 w_1&=\frac{\pi_1}{1-\xi_2/\xi_1},\\
 w_j&=\frac{\pi_j}{1-\xi_{j+1}/\xi_j}
       -\frac{\pi_{j-1}}{\xi_{j-1}/\xi_j-1},\qquad j\geq2.
 \end{aligned}
\end{equation*}

\par\medskip\Needspace{3\baselineskip}\noindent\textit{Proof of Proposition~3.8.}\enspace
For the schedule in Proposition~3.8, the reweighted totals and expanded allocation probabilities are
\begin{equation}\label{eq:bg-normalizer}
 t_k=\frac{\mathrm B(b,\rho+b)}{\mathrm B(b,\rho)}
 \frac{\mathrm B(b+k-1,\rho+1)}{\mathrm B(b+k-1,\rho+b+1)},\qquad t_0=0.
\end{equation}
\begin{equation}\label{eq:bg-expanded-weights}
 \begin{aligned}
 w_1&=1-\left(\frac{b}{b+\rho}\right)^2,\qquad
 w_j=\frac{b}{b+\rho}\frac{\mathrm B(b+j-1,\rho+1)}{\mathrm B(b,\rho)},\quad j\geq2,\\
 \widetilde w_1^{(k)}&=\frac{t_1}{t_k},\qquad
 \widetilde w_j^{(k)}=\frac{b\,t_j}{(2b+\rho+j-1)t_k},\quad 2\leq j\leq k.
 \end{aligned}
\end{equation}
These formulas supplement the compact mixture expression for the weights in the main paper.

\begin{proof}
The beta integral gives
\[
 \xi_k-\xi_{k+1}
 =\frac{\mathrm B(b+k-1,\rho+b+1)}{\mathrm B(b,\rho+b)}>0.
\]
The expression in \eqref{eq:bg-normalizer} therefore equals $\pi_k/(\xi_k-\xi_{k+1})$, which simplifies through gamma-function ratios to
\[
 \frac{t_{k+1}}{t_k}=\frac{2b+\rho+k}{b+\rho+k}>1,
 \qquad
 t_j-t_{j-1}=\frac{b\,t_j}{2b+\rho+j-1},\quad j\geq2.
\]
Applying Proposition~3.4 with $t_0=0$ yields (28).
The normalization and monotonicity of these weights become
clearer by setting $v=\rho/(b+\rho)=\pi_1$ and writing
\begin{equation}\label{eq:bg-mixture-weights}
 w_j=v\mathbf1\{j=1\}+(1-v)\pi_j,\qquad j\geq1.
\end{equation}
This representation shows that the weights sum to one and
decrease strictly, since
$\pi_{j+1}/\pi_j=(b+j-1)/(b+\rho+j)<1$.
Summing \eqref{eq:bg-mixture-weights} over $j>k$ gives the
residual mass in (29), while the stated weight asymptotic
follows from gamma-function ratios. Finally, the finite
allocation probabilities are $(t_j-t_{j-1})/t_k$ and sum to
one as $\sum_{j=1}^k(t_j-t_{j-1})=t_k-t_0=t_k$.
\end{proof}

\subsection{Additional index specifications}\label{sec:additional-examples}\label{sec:light-details}\label{sec:powerlaw-details}
These calculations complement the main examples with explicit conditional probabilities and tail restrictions. Each distribution is specified on the positive integers, and its induced weights are obtained from the same construction identities.

\par\medskip\Needspace{3\baselineskip}\noindent\textit{Negative-binomial and Poisson calculations.}\enspace
For (24), the normalizing constants required by the fixed schedule are
\[
 t_k=p^{r-1}\frac{\Gamma(k+r-1)}{\Gamma(r)\Gamma(k)},\qquad t_0=0.
\]
Since $t_{k+1}/t_k=(k+r-1)/k>1$, the sequence $(t_k)$ is
strictly increasing, with successive differences given by
\[
 t_j-t_{j-1}=p^{r-1}
 \frac{\Gamma(j+r-2)}{\Gamma(r-1)\Gamma(j)},
 \qquad j\geq1,
\]
where $t_0=0$ gives the same expression at $j=1$. Multiplying by $\xi_j=(1-p)^{j-1}$ proves (25). The finite allocation probabilities are
\begin{equation*}
 \widetilde w_j^{(k)}=
 \frac{\Gamma(j+r-2)/\{\Gamma(r-1)\Gamma(j)\}}
 {\Gamma(k+r-1)/\{\Gamma(r)\Gamma(k)\}},\quad j\leq k.
\end{equation*}
Their independence from $p$ yields, provided $p\sim\Beta(a_0,b_0)$,
\[
 p\mid\kb,\cdots\sim\Beta\!\left(a_0+nr,\ b_0+\sum_i(k_i-1)\right).
\]

For a DP with total mass $\alpha$, the usual stick-breaking
order has mean weights $\overline w_j=p(1-p)^{j-1}$, where
$p=(\alpha+1)^{-1}$. Choosing $\xi_j=(1-p)^{j-1}$ then gives
the marginal index probabilities $kp^2(1-p)^{k-1}$.
This marginal calibration extends to
$K-1\sim\NegBin(r,p)$ for every $r>1$, since a
$\NegBin(r,p)$ variable can be represented as the sum of
independent $\NegBin(1,p)$ and $\NegBin(r-1,p)$ variables
and therefore strictly stochastically dominates the first
summand. Corollary~3.3 therefore
provides a compatible deterministic schedule.

The conditional index probabilities nevertheless depend on
the realized DP weights. Indeed, prescribing a fixed
$\pi_1>0$ conditionally on every realization would require
$w_1\geq\pi_1$ almost surely, which is incompatible with
$w_1\sim\Beta(1,\alpha)$. Thus, marginal calibration is
possible without imposing a common conditional index law
across weight realizations.

For the shifted-Poisson construction with $a_j=1$, let
$N_\nu\sim\operatorname{Poisson}(\nu)$. The inverse relation gives
\[
 w_j=\sum_{k\geq j}\frac{e^{-\nu}\nu^{k-1}}{k!}
     =\frac{\Pp(N_\nu\geq j)}{\nu},
 \qquad
 \sum_{j\geq1}w_j=\frac{\E(N_\nu)}{\nu}=1.
\]
It follows that
$\E(K\mid\nu)=1+\nu$ and
$\E(J\mid\nu)=\E\{N_\nu(N_\nu+1)\}/(2\nu)=1+\nu/2$.

With a prior $\nu\sim\mathrm{Gamma}(a_0,b_0)$, uniform finite allocations give
\[
 \nu\mid\kb,\cdots\sim\mathrm{Gamma}\!\left(a_0+\sum_i(k_i-1),\ b_0+n\right).
\]

Every first-$K$-atom construction, even with weight-dependent schedules or arbitrary finite allocation probabilities, satisfies
\begin{equation}\label{eq:count-tail-constraint}
 \Pp(K>k)\geq\Pp(J>k)=\E_\Pi\!\left(\sum_{j>k}w_j\right).
\end{equation}
For a DP in its usual stick-breaking order, the mean residual mass
beyond label $k$ is $q_0^k$, where $q_0=\alpha/(\alpha+1)$.
A shifted-Poisson index $K=1+N$, with
$N\sim\operatorname{Poisson}(\nu)$, instead satisfies the
Chernoff bound
\begin{equation*}
 \Pp(K>k)=\Pp(N\geq k)
 \leq e^{-\nu}(e\nu/k)^k,\qquad k>\nu.
\end{equation*}
Since this bound decays faster than $q_0^k$, it eventually
violates the necessary condition
$\Pp(K>k)\geq\E(R_k)$. Thus, no fixed $\nu>0$ yields a
compatible marginal index law for the DP in this order.

Gamma mixing can remove this obstruction. Integrating the
 prior $\nu\sim\mathrm{Gamma}(a_0,b_0)$ gives
\[
 \Pp(K=k)=\frac{\Gamma(k+a_0-1)}{\Gamma(a_0)\Gamma(k)}
 \left(\frac{b_0}{b_0+1}\right)^{a_0}
 \left(\frac1{b_0+1}\right)^{k-1},\qquad k\geq1,
\]
so $K-1\sim\NegBin(a_0,b_0/(b_0+1))$ marginally.
In particular, $a_0=2$ and $b_0=1/\alpha$ recover the
$\NegBin(2,1/(\alpha+1))$ calibration above, with
$\xi_j=\{\alpha/(\alpha+1)\}^{j-1}$. This agreement concerns
the marginal law of one index. Preserving the DP also requires
its weight prior and corresponding conditional allocation law,
whereas the direct equal-allocation construction retains its
own induced weight law.

For a Pitman--Yor process with $0<\sigma<1$, the mean residual
mass in its usual stick-breaking order has polynomial decay,
\begin{equation*}
 \E(R_k)\sim\frac{C_{\sigma,\alpha}}{1/\sigma-1}
 k^{-(1-\sigma)/\sigma}.
\end{equation*}
Since the tails of fixed-parameter Poisson and negative-binomial
laws decay faster than any power, neither can satisfy the necessary
constraint \eqref{eq:count-tail-constraint} at every cutoff.
This restriction applies when preserving an existing SSP
in its defining order, whereas Theorem~3.1
takes the index distribution and finite allocation probabilities
as starting points and derives the corresponding weight law.

\par\medskip\Needspace{3\baselineskip}\noindent\textit{Zipf and discrete-Pareto specifications.}\enspace
Power-tailed index laws provide further examples and illustrate how
specifying point probabilities or survival probabilities can lead to
different SSPs, even with matching tail exponents. To make this
distinction explicit, we use $a_j=2j-1$ in both constructions, so
$s_k=k^2$ and the finite allocation probabilities are
$(2j-1)/k^2$, favoring later eligible labels.

For $\tau>1$, consider the Zipf law
$\pi_k=k^{-\tau}/\zeta(\tau)$, where
$\zeta(s,q)=\sum_{\ell=0}^{\infty}(q+\ell)^{-s}$ is the Hurwitz
zeta function for $s>1$, $q>0$, and $\zeta(s)=\zeta(s,1)$.
Theorem~3.1 then yields
\begin{equation}\label{eq:zipf-ssp}
 \xi_j=\frac{\zeta(\tau+2,j)}{\zeta(\tau)},\qquad
 w_j=(2j-1)\frac{\zeta(\tau+2,j)}{\zeta(\tau)},\qquad\mbox{and}\qquad
 \widetilde w_j^{(k)}=\frac{2j-1}{k^2},\quad j\leq k.
\end{equation}
Indeed, $\xi_k-\xi_{k+1}=k^{-(\tau+2)}/\zeta(\tau)$,
which, multiplied by $s_k=k^2$, recovers $\pi_k$.
Integral comparison of the defining series gives
\begin{equation*}
 \begin{aligned}
 \Pp(K\geq k)&=\frac{\zeta(\tau,k)}{\zeta(\tau)}
       \sim\frac{k^{1-\tau}}{(\tau-1)\zeta(\tau)},\\
 w_j&\sim\frac{2j^{-\tau}}{(\tau+1)\zeta(\tau)},\qquad\mbox{with}\qquad
 \E(K)=\frac{\zeta(\tau-1)}{\zeta(\tau)}\quad(\tau>2).
 \end{aligned}
\end{equation*}
The mean index is infinite for $1<\tau\leq2$, while the
conditional mean of the selected label is
$\E(J\mid K=k)=(k+1)(4k-1)/(6k)$, reflecting the preference
for later labels.

Alternatively, a discrete-Pareto law can be specified through
its survival probabilities,
\begin{equation}\label{eq:pareto-index}
 \Pp(K\geq k)=k^{-\rho},\qquad\mbox{and}\qquad
 \pi_k=k^{-\rho}-(k+1)^{-\rho},\qquad \rho>0.
\end{equation}
This gives $K=\lfloor U^{-1/\rho}\rfloor$ for
$U\sim\mathrm{Unif}(0,1)$. Using the same coefficients,
we obtain
\begin{equation}\label{eq:pareto-ssp}
 \xi_j=\sum_{k\geq j}\frac{k^{-\rho}-(k+1)^{-\rho}}{k^2},
 \qquad
 w_j=(2j-1)\xi_j,\qquad\mbox{and}\qquad
 \widetilde w_j^{(k)}=\frac{2j-1}{k^2}.
\end{equation}
Here $\xi_j\sim\rho j^{-(\rho+2)}/(\rho+2)$ and
$w_j\sim2\rho j^{-(\rho+1)}/(\rho+2)$, while
$\E(K)=\zeta(\rho)$ for $\rho>1$ and is infinite otherwise.
Taking $\tau=\rho+1$ matches the Zipf tail exponent,
but the two laws retain different probabilities at finite
labels. Among the various discrete-Pareto conventions
\citep{BuddanaKozubowski2014}, we adopt
\eqref{eq:pareto-index}.

Both constructions define proper SSP priors when combined
with independent i.i.d.\ atoms from $G_0$. Assigning a prior
to the exponent makes the weights random, while leaving the
finite allocation probabilities unchanged. Consequently,
the augmented likelihood depends on the exponent only
through $\prod_i\pi_{k_i}$, allowing inference directly from
the  auxiliary indices. 
\subsection{Weak support of the induced priors}\label{sec:prior-support}\label{sec:power-support}
The preceding constructions yield proper SSPs, although full weak
support requires additional flexibility in their weight laws.
Let $\mathsf X_0=\operatorname{supp}(G_0)$, with atoms drawn i.i.d.\
from diffuse $G_0$, independently of the weights. The criterion
of \citet[Corollary~3]{BissiriOngaro2014}, recalled in (18),
reduces full weak support to the possibility that the largest
weight is arbitrarily small with positive probability. We examine this property for the light- and heavy-tailed families in the main paper and the preceding examples. Theorem~3.5 of the main paper establishes the general survival-transformation result.

\Needspace{3\baselineskip}\textit{Restrictions under a fixed conditional index law.}
Suppose that the same probabilities $\pi_k$ specify the conditional
law of $K$ for every realization of the weights. Since $J\leq K$,
the inclusion $\{K\leq h\}\subseteq\{J\leq h\}$ implies, for
$Q_h=\sum_{k\leq h}\pi_k>0$,
\[
 \sup_jw_j\geq\frac1h\sum_{j=1}^h w_j
 \geq\frac{Q_h}{h}>0
 \quad\text{almost surely}.
\]
This bound precludes full weak support, regardless of the
randomness assigned to the finite allocation probabilities.
In particular, the strictly positive index laws considered
here imply $w_1\geq\pi_1>0$.

Full weak support can instead be achieved by allowing the index
probabilities to vary with a parameter $\eta\sim\Lambda$.
For the induced weights $w_j(\eta)$, the required condition becomes
\[
 \Lambda\{\eta:\sup_jw_j(\eta)<\varepsilon\}>0
 \quad\text{for every }\varepsilon>0.
\]
The parameter $\eta$ is drawn once per random measure and shared
across observations, so the two-stage representation and augmented
likelihood (10) remain valid conditionally on $\eta$. The following
examples establish this property under suitable hyperpriors.
By contrast, the marginal calibration of Corollary~3.3 preserves
the original weight prior and hence its support.

\Needspace{3\baselineskip}\textit{Light-tailed index families.}
For the negative-binomial construction in Example~3.6, fix
$r>1$ and put $s=r-1$. The weights in (25) are the shifted
probabilities of a $\NegBin(s,p)$ distribution. When
$0<s\leq1$, they decrease and their maximum is $p^s$.
For $s>1$, convolution with a $\NegBin(1,p)$ distribution
bounds the maximum by $p$. Thus,
\[
 p^{r-1}\leq\sup_jw_j\leq p^{\min(r-1,1)}.
\]
Full weak support holds precisely when the prior on $p$
assigns positive probability arbitrarily close to zero,
as does any $\Beta(a_0,b_0)$ prior with positive shapes.
Under this prior, $\E(K)=1+rb_0/(a_0-1)$ when $a_0>1$,
and is infinite otherwise. The case $r=2$ recovers the
geometric-weight support condition of
\citet[Section~3.3]{BissiriOngaro2014}.

For the equal-allocation Poisson construction in Example~3.7,
the largest weight is $(1-e^{-\nu})/\nu$, which approaches zero
as $\nu\to\infty$ but remains bounded away from zero when
$\nu$ is bounded above. Consequently, full weak support holds
exactly when the hyperprior on $\nu$ has unbounded support above.
The $\mathrm{Gamma}(a_0,b_0)$ prior satisfies this condition
and gives $\E(K)=1+a_0/b_0$. We next consider how analogous
arguments apply to the beta-geometric scale families.

\begin{theorem}[Full weak support with finite mean index]\label{thm:bg-support}
Fix $\rho>1$, give $b>0$ a prior $H_b$, and construct $G$ with weights (28). This SSP prior has full weak support on $\mathcal P(\mathsf X_0)$ if and only if $H_b((M,\infty))>0$ for every finite $M>0$. Moreover,
\begin{equation}\label{eq:bg-mean-index}
 \E(K)=1+\frac{\E_{H_b}(b)}{\rho-1},
\end{equation}
with $+\infty$ permitted. In particular, every gamma prior $b\sim\mathrm{Gamma}(c,d)$ with $c>0$ and rate $d>0$ gives full weak support and $\E(K)=1+c/\{d(\rho-1)\}<\infty$.
\end{theorem}

\begin{proof}
Since the weights are decreasing, their maximum is
$w_1=1-\{b/(b+\rho)\}^2$. For $0<\varepsilon<1$, this gives
\begin{equation*}
 \Pp\!\left(\sup_jw_j<\varepsilon\right)
 =H_b\!\left(\left(\frac{\rho\sqrt{1-\varepsilon}}
                        {1-\sqrt{1-\varepsilon}},\infty\right)\right).
\end{equation*}
As the threshold tends to infinity when $\varepsilon\downarrow0$,
the full-support criterion (18) holds if and only if $H_b$ has
unbounded support above. Moreover, for $\rho>1$, the beta-geometric
representation gives $\E(K\mid b)=1+b/(\rho-1)$, whose expectation over $H_b$ yields \eqref{eq:bg-mean-index}. In particular, a
$\mathrm{Gamma}(c,d)$ prior has unbounded support and mean $c/d$.
\end{proof}

\begin{proposition}[Alternating weights and their support]
\label{prop:alternating-support}
For fixed $\rho>1$ and $\gamma>1$, the weights in (30)
satisfy the asymptotic relations in (31). If $b\sim H_b$,
the induced SSP has full weak support if and only if
$H_b((M,\infty))>0$ for every finite $M>0$.
Its expected representation index is given by
\eqref{eq:bg-mean-index}.
\end{proposition}

\begin{proof}
Theorem~3.1 ensures positivity, normalization and the prescribed
index law for every $b>0$. Since
$s_k/k\to\overline a=(1+\gamma)/2$ and
$\pi_k\sim C_{\rho,b}k^{-(\rho+1)}$, we have
\[
 \frac{\pi_k}{s_k}
 \sim\frac{C_{\rho,b}}{\overline a}k^{-(\rho+2)},
 \qquad
 \xi_j\sim\frac{C_{\rho,b}}{\overline a(\rho+1)}
 j^{-(\rho+1)},
\]
where the second equivalence follows by comparing the series
with its corresponding integral. Consequently,
$\xi_{j+1}/\xi_j\to1$, and the alternating coefficients
yield the ratios in (31), establishing infinitely many
strict increases and decreases.

To establish full weak support, notice that
$1\leq a_j\leq\gamma$ and $s_k\geq k$ imply
\begin{equation}\label{eq:alternating-support-bound}
 \sup_j w_j\leq\gamma\sum_{k\geq1}\frac{\pi_k(\rho,b)}{k}
 =\gamma\E(K^{-1}\mid b).
\end{equation}
For any fixed integer $L$, the beta-geometric survival
probabilities satisfy
\[
 \Pp(K>L\mid b)
 =\prod_{h=0}^{L-1}\frac{b+h}{b+\rho+h}
 \longrightarrow1\qquad(b\to\infty).
\]
Moreover,
\[
 \E(K^{-1}\mid b)
 \leq\Pp(K\leq L\mid b)+\frac{1}{L+1}.
\]
Letting $b\to\infty$ and then $L\to\infty$ shows that
the upper bound in \eqref{eq:alternating-support-bound}
approaches zero. Thus any $H_b$ with unbounded support
above yields full weak support.

Conversely, $w_1\geq\pi_1(\rho,b)=\rho/(b+\rho)$.
If $H_b$ is supported on $(0,M]$, the largest weight
is thus bounded below by $\rho/(M+\rho)$, ruling
out full  support. Finally, changing the allocation
coefficients leaves the index law unchanged, so its
expected value remains \eqref{eq:bg-mean-index}.
\end{proof}

\par\medskip\Needspace{3\baselineskip}\noindent\textit{Power-law examples.}\enspace
For the Zipf weights in \eqref{eq:zipf-ssp}, the finite allocation
probabilities give
\[
 \sup_jw_j
 \leq\sum_{k\geq1}\pi_k
       \max_{1\leq j\leq k}\frac{2j-1}{k^2}
 \leq\frac{2\zeta(\tau+1)}{\zeta(\tau)}
 \longrightarrow0
 \qquad(\tau\downarrow1).
\]
The limit follows from the divergence of $\zeta(\tau)$ as
$\tau\downarrow1$. Thus, a prior on $\tau$ yields full weak
support precisely when it assigns positive probability arbitrarily
close to one. Indeed, if $\tau\geq1+\delta$ almost surely for
some $\delta>0$, then
$w_1\geq1/\zeta(1+\delta)>0$, ruling out full weak support.
Such a prior necessarily assigns positive probability to the
range $1<\tau\leq2$, where the conditional mean index is
infinite, and consequently $\E(K)=\infty$.

For the discrete-Pareto law in \eqref{eq:pareto-index}, the same
allocation bound gives
$\sup_jw_j\leq2\E(K^{-1}\mid\rho)$. Since, for every fixed $L$,
\[
 \Pp(K>L\mid\rho)=(L+1)^{-\rho}\longrightarrow1
 \qquad(\rho\downarrow0),
\]
we have $\E(K^{-1}\mid\rho)\to0$. Hence full weak support
holds whenever the exponent prior assigns positive probability
arbitrarily close to zero. This condition is also necessary,
because $w_1\geq1-2^{-\rho}$ remains bounded away from zero
when $\rho$ is bounded below. As in the Zipf construction,
full weak support entails positive prior probability of an
infinite conditional mean index, here for $0<\rho\leq1$.

The asymptotic relations following \eqref{eq:pareto-ssp}
follow from
$k^{-\rho}-(k+1)^{-\rho}\sim\rho k^{-(\rho+1)}$
and integral comparison. In both constructions, the conditional
mean of the selected label is
\[
 \E(J\mid K=k)
 =\frac{1}{k^2}\sum_{j=1}^k j(2j-1)
 =\frac{(k+1)(4k-1)}{6k}.
\]

\par\medskip\Needspace{3\baselineskip}\noindent\textit{Largest weights and expected index.}\enspace
\label{sec:support-cost-numerics}

\begin{figure}[tbp]
\centering
\includegraphics[width=.88\linewidth]{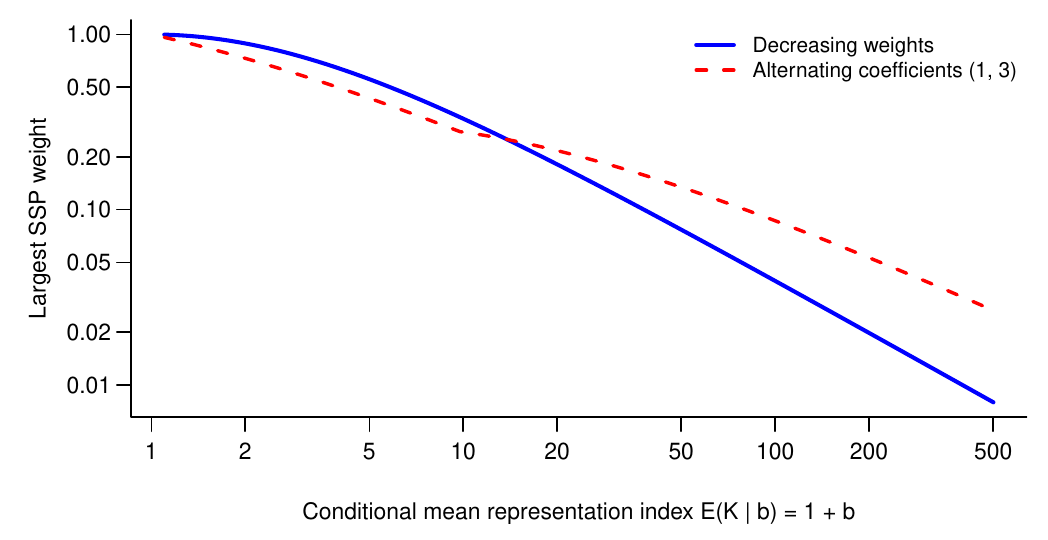}
\caption{Largest SSP weight against conditional mean index
for the decreasing and alternating families, with $\rho=2$
and varying $b$. Each point corresponds to a fixed parameter value.}
\label{fig:support-cost}
\end{figure}

Figure~\ref{fig:support-cost} illustrates how the largest
weight decreases as the conditional mean index grows, with
$\rho=2$ and varying $b$. In the decreasing family, the maximum
is $1-\{b/(b+2)\}^2$, while in the alternating family it occurs
among the first two weights, since each parity subsequence
decreases. Both maxima approach zero as $b\to\infty$.
An unbounded prior on $b$ therefore gives full weak support
without necessarily imposing an infinite expected index.
For instance, $b\sim\mathrm{Gamma}(2,1/5)$  yields $\E(K)=11$ in both families.

\subsection{Paired weights and their invariant order}\label{app:paired-proofs}
We establish Theorem~3.10, including the marginal index law and its compatible conditional weight kernel. The atoms remain independent $G_0$ draws throughout. The equal splitting of each original DP weight is a special case of the fragmentation operator in \citet[Section~5.4, equation~(5.12)]{Pitman2006}, with fragmentation law concentrated on $(1/2,1/2,0,\ldots)$. The arguments below identify the resulting ordered prior and its finite sampling representation.

\par\medskip\Needspace{3\baselineskip}\noindent\textit{Index law specified marginally.}\enspace
For any infinite sequence of positive weights and any $h\geq1$, there are $h$ distinct weights whose sum is smaller than $\varepsilon$. A size-biased permutation places these weights in its first $h$ positions with positive probability. Conditioning on the masses and then averaging therefore shows that an ISBP law satisfies
\[
 \Pp\!\left(\sum_{j=1}^h w_j<\varepsilon\right)>0
 \quad\text{for every }h\geq1,\ \varepsilon>0.
\]
If $\Pp(K=k\mid\wb)=\pi_k$ at every realization, the restriction $J\leq K$ requires $Q_h=\sum_{k\leq h}\pi_k\leq\sum_{j\leq h}w_j$ almost surely. Since $Q_h>0$ for some finite $h$, this contradicts the preceding property. The argument permits arbitrary allocation probabilities supported on the first $K$ labels. Specifying the index law marginally avoids the obstruction because its conditional probabilities may vary with $\wb$.

\par\medskip\Needspace{3\baselineskip}\noindent\textit{Defining order and properness.}\enspace
The original DP weights $w_m^{\mathrm{DP}}$ in
(33) are generated by independent
$\Beta(1,\alpha)$ stick-breaking fractions \citep{Sethuraman1994}.
When size-biased sampling first enters a pair, it selects that
pair proportionally to its total weight, with either atom equally
likely to appear first. The DP's size-biased stick-breaking representation therefore gives
an independent $\Beta(1,\alpha)$ stick-breaking fraction of the mass
in pairs not yet encountered \citep[Theorem~1.1]{PermanPitmanYor1992}.
After normalization, the remaining intact pairs retain their
original DP weight law \citep{Pitman2006}.

Before generating $w_j$, let $U$ be the combined weight of
pairs with neither member selected, and let $\mathcal H$ list the
unselected partner's weight from each pair represented once among
$w_1,\ldots,w_{j-1}$. Thus $\mathcal H$ has at most $j-1$ entries.
Initially $(U,\mathcal H)=(1,\varnothing)$, and the next choice
follows (34), with
$r=U+\sum_{b\in\mathcal H}b=1-\sum_{\ell<j}w_\ell$.
Selecting an entry $b\in\mathcal H$ has probability $b/r$,
removes that entry and leaves $U$ unchanged. A new pair is selected
with probability $U/r$. An independent $V\sim\Beta(1,\alpha)$
then gives $w_j=UV/2$, adds the equal partner to $\mathcal H$
and replaces $U$ by $U(1-V)$. Each selected weight receives
an independent $G_0$ location, so equal weights impose no
similarity on the atom locations.

This rule selects the atoms successively without replacement with probabilities proportional to their remaining masses. For a fixed atom of mass $p>0$, its conditional selection probability at every preceding step is at least $p$. Its probability of remaining unselected after $k$ steps is therefore at most $(1-p)^k$. Every atom is eventually selected almost surely, and countability ensures that the resulting weights sum to one.

Given the ranked masses, their generated order is size-biased. Applying an independent size-biased permutation preserves this conditional ordering law, proving ISBP \citep{PitmanISBP1996}. Since the ordering depends only on the masses, the locations remain independent of the weights, and the unlabelled measure has law (33).

\par\medskip\Needspace{3\baselineskip}\noindent\textit{Full weak support and separation from normalized CRMs.}\enspace
The largest paired weight is $\frac12\sup_m w_m^{\mathrm{DP}}$. Fix $\varepsilon>0$ and choose $M$ such that $1/(M+1)<2\varepsilon$. If the first $M$ DP stick-breaking fractions are $1/(M+2-m)$, for $m=1,\ldots,M$, then their weights and the remaining mass all equal $1/(M+1)$. The beta densities are positive on $(0,1)$, so continuity gives positive probability that these stick-breaking fractions lie sufficiently close to the indicated values for every generated weight and the remaining mass to be smaller than $2\varepsilon$. Every later DP weight is bounded by that remainder. Hence
\[
 \Pp\!\left(\sup_jw_j<\varepsilon\right)>0
 \quad\text{for every }\varepsilon>0.
\]
Corollary~3 of \citet{BissiriOngaro2014}, applied with diffuse independent atoms, proves full weak support on $\mathcal P(\mathsf X_0)$ for each fixed $\alpha>0$.

To exclude normalized homogeneous CRMs, first note that the original DP weights are pairwise distinct almost surely, by continuity of their beta stick-breaking fractions. Thus every distinct positive weight in the paired measure occurs exactly twice. Normalization preserves equality and multiplicity of jump sizes. A homogeneous drift-free CRM with deterministic L\'evy intensity $\rho(\dd t)G_0(\dd\theta)$ and diffuse $G_0$ has distinct atom locations, while its jump sizes form a Poisson random measure with intensity $\rho$ \citep{Kingman1967}.

If the atomless part of $\rho$ is nonzero, it assigns finite positive mass to some set of sizes bounded away from zero and infinity. With positive probability this set contributes a jump, whose size differs from all other jumps almost surely. Its normalized weight occurs only once. If $\rho$ is purely atomic, choose a size $t_0>0$ with $c=\rho(\{t_0\})>0$. L\'evy integrability implies $c<\infty$, and its Poisson count equals one with probability $ce^{-c}>0$. Both cases contradict the required multiplicity two. Therefore no such CRM with finite strictly positive total mass has the paired normalized law. This includes superpositions of independent homogeneous CRMs of this form, whose L\'evy intensities add.

The first generated weight satisfies $2w_1\sim\Beta(1,\alpha)$ and is always smaller than $1/2$. The first size-biased weight of an infinite Pitman--Yor law has a $\Beta(1-\sigma,\theta+\sigma)$ distribution, with positive probability above $1/2$ \citep{PitmanYor1997}. Finite-support Pitman--Yor cases are excluded by their support cardinality. The paired law thus also differs from every Pitman--Yor law with fixed parameters.

\par\medskip\Needspace{3\baselineskip}\noindent\textit{The marginal index distribution.}\enspace
Let $\overline w_j=\E_{\Pi^{(2)}_\alpha}(w_j)$. Consider an independent size-biased ordering of the paired masses, obtained by ordering independent exponential clocks with rates proportional to those masses. Independently select a tagged atom proportionally to its mass. Its rank $D$ has probabilities $\Pp(D=j)=\overline w_j$ after averaging over the masses and their ordering.

To derive this rank distribution, represent the original DP weights
by normalized gamma-process jumps, with intensity
$\alpha t^{-1}e^{-t}\dd t$ and total mass $T$ \citep{Ferguson1973}.
A jump of size $t$ creates a pair with total weight $t/T$.
Both atoms may be assigned exponential clock rate $t$, since
common rescaling preserves the clock order. Selecting the pair
containing the tagged atom contributes $t/T$, and either of its
atoms is tagged with probability $1/2$.

For $0\leq z\leq1$, let $F(z)=\E(z^{D-1})$. Conditional on
the tagged jump $t$ and the tagged atom's clock time $s$,
$D-1$ counts the other clocks that ring before $s$.
The partner contributes the factor $z+(1-z)e^{-st}$,
while each other original jump $u$ contributes
$\{z+(1-z)e^{-su}\}^2$, since its two atoms have
independent clocks.

Using $1/T=\int_0^\infty e^{-rT}\dd r$, the Poisson
Mecke identity and Laplace functional
\citep[Theorems~3.9 and~4.1]{LastPenrose2017} give the
factor
\[
 \exp\left\{
 \alpha\int_0^\infty\frac{e^{-u}}u
 \left[e^{-ru}\{z+(1-z)e^{-su}\}^2-1\right]\dd u
 \right\}
\]
for the remaining jumps. After expanding the square and
evaluating the resulting logarithmic integrals, we integrate
over the tagged jump and set $A=1+r$ to obtain
\begin{equation}\label{eq:paired-pgf-double}
 \begin{aligned}
 F(z)=\alpha\int_1^\infty\int_0^\infty
 &\left\{\frac{z}{(A+s)^2}
          +\frac{1-z}{(A+2s)^2}\right\}\\
 &\quad\times A^{-\alpha z^2}
 (A+s)^{-2\alpha z(1-z)}
 (A+2s)^{-\alpha(1-z)^2}\dd s\dd A.
 \end{aligned}
\end{equation}
Tonelli's theorem justifies the exchange of the nonnegative
sampling integrals.

With $s=Av$, the integral over $A$ becomes
$\int_1^\infty A^{-\alpha-1}\dd A=1/\alpha$,
which cancels the leading factor $\alpha$ and leaves
\[
 F(z)=\int_0^\infty
 \left\{\frac z{(1+v)^2}
          +\frac{1-z}{(1+2v)^2}\right\}
 (1+v)^{-2\alpha z(1-z)}
 (1+2v)^{-\alpha(1-z)^2}\dd v.
\]
Substituting $x=(1+v)^{-1}$ and writing
$b(x)=1-(2-x)^{-2}$,
$\lambda_1(x)=2\alpha\log(2-x)$ and
$\lambda_2(x)=-\alpha\log\{x(2-x)\}$ gives
\begin{equation}\label{eq:paired-pgf}
 \E(z^{D-1})=\int_0^1\{1-b(x)+b(x)z\}
 \exp\{\lambda_1(x)(z-1)+\lambda_2(x)(z^2-1)\}\dd x.
\end{equation}
For each $x$, the integrand is the conditional generating
function of $B+N_1+2N_2$ in the finite experiment.
Averaging over $X\sim\operatorname{Unif}(0,1)$ therefore
proves the stated law of $D$ by uniqueness of probability
generating functions. Adding an independent
$S_0\sim\Geom_0(1-q)$ then gives the marginal probabilities
of $K=D+S_0$ in (35).

The same experiment also gives the mean directly.
Since $\int_0^1b(x)\dd x=1/2$ and
$\lambda_1(x)+2\lambda_2(x)=-2\alpha\log x$,
\[
 \E(D)=1+\frac12-2\alpha\int_0^1\log x\dd x
 =2\alpha+\frac32.
\]
Together with $\E(S_0)=q/(1-q)$, this proves
(37).

Finally, for $z>1$, the conditional generating function
of the finite experiment is bounded by
$z x^{-\alpha(z^2-1)}$ on $0<x<1$. Choosing
$1<z<\min\{\sqrt{1+1/\alpha},q^{-1}\}$ makes this bound
integrable and ensures $\E(z^{S_0})<\infty$.
Independence therefore gives $\E(z^K)<\infty$, and hence
finite moments of every positive order for both $D$ and $K$.
This establishes the distributional assertions of
Theorem~3.10.

\par\medskip\Needspace{3\baselineskip}\noindent\textit{Conditional weights and finite simulation.}\enspace
For $\xi_j=q^{j-1}$, the conditional index probabilities
$p_k(\wb,\xib)=(1-q)\sum_{j\leq k}w_jq^{k-j}$ satisfy
$\E_{\Pi^{(2)}_\alpha}\{p_k(\wb,\xib)\}=\pi_k$.
Thus the kernel in (36) integrates
to one, and averaging over $K\sim\boldsymbol\pi$ recovers
the weight prior,
$\sum_k\pi_k\Pi^{(2)}_{\alpha,k}=\Pi^{(2)}_\alpha$.
Combining these conditional laws with the allocation
probabilities $\widetilde w_j^{(k)}=a_j/s_k$ gives
\[
 \sum_{k\geq j}\pi_k\widetilde w_j^{(k)}
       \Pi^{(2)}_{\alpha,k}(\dd\wb)
 =w_j\Pi^{(2)}_\alpha(\dd\wb).
\]
Consequently, drawing $K$, then its conditional weights,
and finally $J$ with probabilities $\widetilde w_j^{(K)}$
recovers the original joint law of the weights and selected
label. For a sample, all observations share one draw of
$\wb$, with their indices conditionally independent and
distributed according to $p_k(\wb,\xib)$. The auxiliary
$D$ used to generate the marginal index can be discarded
once $K$ is known and need not equal the final label $J$.

For simulation, this joint law can also be generated
sequentially, determining $K$ as successive weights are
drawn from the paired prior recursion. After step $k$,
the state gives the residual mass
$R_k=U+\sum_{b\in\mathcal H}b$, while the scaled total
$A_k=q^ks_k$ can be updated through
\[
 A_0=0,\qquad A_k=q(A_{k-1}+w_k).
\]
Together, these quantities give the conditional survival
probability
\begin{equation}\label{eq:paired-lazy-survival}
 \Pp(K>k\mid\wb,\xib)=R_k+A_k.
\end{equation}
Draw $U_0\sim\operatorname{Unif}(0,1)$ independently of
the weight recursion and generate weights until
\[
 K=\min\{k\geq1:R_k+A_k\leq U_0\}.
\]
Then select $J\leq K$ with probabilities
$\widetilde w_j^{(K)}=a_j/s_K$.
Properness and $q<1$ imply $R_k+A_k\to0$, so the
procedure terminates almost surely after generating
exactly $K$ weights, with $\E(K)$ given by
(37). Maintaining and sampling
from $\mathcal H$ entails additional computational work.

For a prescribed $K=k$, an alternative is rejection
sampling. Propose $w_1,\ldots,w_k$ from the prior recursion
and accept with probability
\[
 \Pp(\text{accept}\mid w_1,\ldots,w_k,K=k)
 =\frac{p_k(\wb,\xib)}{1-q}
 =\sum_{j=1}^k w_jq^{k-j}\leq1.
\]
The acceptance rate is $\pi_k/(1-q)>0$, and the accepted
weights have the required marginal under
$\Pi^{(2)}_{\alpha,k}$. Although this sampler terminates
almost surely for every fixed $k$, drawing
$K\sim\boldsymbol\pi$ first gives an infinite expected
number of proposals. Indeed, writing $N_{\rm prop}$ for
that count,
\[
 \E(N_{\rm prop})
 =\sum_{k\geq1}\pi_k\frac{1-q}{\pi_k}
 =\infty.
\]
The sequential construction avoids this rejection step,
illustrating why finite $\E(K)$ alone does not guarantee
finite expected computational cost for a chosen
conditional-weight sampler.

\subsection{Survival transformations and full-support hierarchies}
\label{sec:survival-proofs}
We prove Theorem~3.5 and apply its
square-root specialization to thirteen illustrative index families.
Full weak support is verified through
(18), following
\citet{BissiriOngaro2014}.

\par\medskip\Needspace{3\baselineskip}\noindent\textit{Proof of Theorem~3.5.}\enspace
Since $H_j-H_{j+1}=\pi_j>0$, the quantities
\[
 t_j^{-1}
 =\frac{\phi(H_j)-\phi(H_{j+1})}{H_j-H_{j+1}}
\]
are positive secant slopes. Strict convexity gives
$t_{j+1}^{-1}<t_j^{-1}$, hence
$0=t_0<t_1<t_2<\cdots$ and $a_j=t_j-t_{j-1}>0$.
Continuity at zero gives $\xi_j=\phi(H_j)\downarrow0$,
so Proposition~3.4 yields
$w_j>0$ and $\sum_jw_j=1$. Moreover,
$\pi_k/t_k=\xi_k-\xi_{k+1}$ implies
\begin{equation}\label{eq:survival-joint-check}
 \Pp(J=j\mid\eta)
 =\sum_{k\geq j}\pi_k\frac{a_j}{t_k}
 =a_j\sum_{k\geq j}(\xi_k-\xi_{k+1})
 =a_j\xi_j=w_j.
\end{equation}
Multiplying by $\delta_{\theta_j}(A)$ and summing over $j$
therefore gives $\sum_k\pi_kG_k^\star(A)=G(A)$ for every
$A\in\mathcal B$.

To verify full weak support, put $\Delta=\sup_k\pi_k$ and
define the modulus of continuity
\[
 \omega_g(r)=\sup_{\substack{x,y\in[0,1]\\|x-y|\leq r}}
 |g(x)-g(y)|,\qquad 0\leq r\leq1.
\]
Convexity gives $g(H_j)\leq t_j\leq g(H_{j+1})$.
Together with $H_1=1$ and $g(1)=0$, these inequalities yield
\[
 \begin{aligned}
 a_1&=t_1\leq g(H_2)-g(1)\leq\omega_g(\Delta),\\
 0<a_j&\leq g(H_{j+1})-g(H_{j-1})
       \leq\omega_g(\pi_{j-1}+\pi_j),\qquad j\geq2.
 \end{aligned}
\]
Since $w_j=\xi_ja_j$ and $\xi_j\leq\phi(1)$,
\begin{equation}\label{eq:survival-uniform-bound}
 \sup_jw_j\leq\phi(1)\omega_g\!\left(\min\{2\Delta,1\}\right).
\end{equation}
Uniform continuity of $g$ makes the right-hand side tend
to zero with $\Delta$. For every $\varepsilon>0$, a sufficiently
small $\delta>0$ therefore gives
\[
 \Pp\!\left(\sup_jw_j<\varepsilon\right)
 \geq\Pp\{\Delta(\eta)<\delta\}>0,
\]
where the last inequality is
(21). This proves full weak
support by (18).

For $\phi(x)=1-(1-x)^\vartheta$, $0<\vartheta<1$,
\[
 g(x)=\vartheta^{-1}(1-x)^{1-\vartheta},\qquad
 \omega_g(r)=\vartheta^{-1}r^{1-\vartheta},\qquad
 w_1=\pi_1^{1-\vartheta},
\]
which gives (22).
Finally, Tonelli's theorem yields
$\E(K)=\E_\eta\{\sum_k k\pi_k(\eta)\}$, with $+\infty$
allowed.\hfill$\square$

\par\medskip\Needspace{3\baselineskip}\noindent\textit{Explicit square-root weights and remaining mass.}\enspace
For $\phi(x)=1-\sqrt{1-x}$, set
$\xi_j=1-\sqrt{1-H_j}$, $\xi_0=1$ and $t_0=0$.
Rationalizing the denominator in $t_j$ gives
\begin{equation}\label{eq:survival-square-root-totals}
 t_j=\sqrt{1-H_j}+\sqrt{1-H_{j+1}}
     =2-\xi_j-\xi_{j+1},\qquad j\geq1.
\end{equation}
Thus $a_j=\xi_{j-1}-\xi_{j+1}$, including $j=1$, and
$s_k=\sum_{j\leq k}a_j=t_k$. Summing the weights gives
\[
 \sum_{j=1}^kw_j
 =\sum_{j=1}^k\xi_j(\xi_{j-1}-\xi_{j+1})
 =1-\xi_k\xi_{k+1},
 \qquad R_k=\xi_k\xi_{k+1},
\]
proving the remaining-mass identity in
(23).

For example, $K\sim\operatorname{Geom}(p)$ has
$H_j=r^{j-1}$, where $r=1-p$, and induces
\begin{equation}\label{eq:survival-geometric-weights}
 \begin{aligned}
 w_1&=\sqrt p,\\
 w_j&=\bigl(1-\sqrt{1-r^{j-1}}\bigr)
       \bigl(\sqrt{1-r^j}-\sqrt{1-r^{j-2}}\bigr),
       \qquad j\geq2.
 \end{aligned}
\end{equation}
The finite probabilities are
$\widetilde w_j^{(k)}=a_j/t_k$, with
$t_k=\sqrt{1-r^{k-1}}+\sqrt{1-r^k}$.
These weights differ from the geometric weights in
Example~3.6, which arise from uniform
allocation and a shifted negative-binomial index law.

The transformation also determines whether full weak support
is possible. For $\phi(x)=x^c$ with fixed $c>1$,
\[
 w_1=\frac{1-H_2}{1-H_2^c}>\frac1c,
\]
so randomizing the index law cannot remove this restriction.
More generally, a fixed differentiable convex transformation
with endpoint limit $g(1)>0$ satisfies
\[
 \sup_jw_j\geq w_1=\phi(1)t_1\geq\phi(1)g(1)>0.
\]
Thus $g(1)=0$ is necessary within this class, while
Theorem~3.5 supplies sufficient
regularity and hyperprior conditions.

\par\medskip\Needspace{3\baselineskip}\noindent\textit{Illustrative index families.}\enspace
We use $\operatorname{Geom}(p)$ on $\N$ and
$\NegBin(r,p)$ on $\{0,1,\ldots\}$, with the latter
counting failures. The beta-geometric law is specified by
$P\sim\Beta(a,b)$ and $K\mid P\sim\operatorname{Geom}(P)$.
The scale families take the form $K=1+\lfloor bY\rfloor$,
where
\[
 \begin{array}{c|c}
 \text{Index family}&\Pp(Y>y),\quad y\geq0\\ \hline
 \operatorname{dLomax}(3,b)&(1+y)^{-3}\\
 \operatorname{dWeibull}(2,b)&e^{-y^2}\\
 \operatorname{dLN}(0,1;b)&\overline\Phi(\log y),\quad y>0\\
 \operatorname{dGamma}(2,1;b)&e^{-y}(1+y)
 \end{array}
\]
and $\overline\Phi$ is the standard normal survival function.
The Pareto families are defined directly by their discrete
survivals, using $\rho$ for the tail exponent.

Write $Z_{r,p}\sim\NegBin(r,p)$,
$N_\nu\sim\operatorname{Poisson}(\nu)$ and
$C\sim\operatorname{Bernoulli}(1/2)$, and put
\[
 \mathcal B_r(m;p)=\Pp(Z_{r,p}\geq m),\qquad
 \mathcal P_\nu(m)=\Pp(N_\nu\geq m),\qquad m\geq0.
\]
Both survivals equal one at $m=0$.
In the mixture specifications, first draw $C$ and then the
indicated $Z$ or $N$ conditional on $C$.

For each row, $a_1,b_1>0$ are fixed hyperparameters and gamma
laws use the shape--rate convention. The shifted scales have
$b-b_0\sim\mathrm{Gamma}(a_1,b_1)$, with fixed $b_0\geq0$
and $m_b=b_0+a_1/b_1$.
The index parameter $\eta$ is drawn once, independently of
the atoms, and shared across observations. The mixture
indicator $C$ and beta-geometric variable $P$ are instead
drawn afresh for each index, conditional on $\eta$.

\begin{table}[p]
\centering
\begingroup
\renewcommand{\baselinestretch}{1}
\fontsize{10}{12}\selectfont
\setlength{\tabcolsep}{3pt}
\renewcommand{\arraystretch}{1.9}
\caption[Illustrative full-support index hierarchies]{Thirteen
illustrative index hierarchies for
Theorem~3.5. With
$\xi_j=1-\sqrt{1-H_j}$ and $\xi_0=1$, the SSP weights are
$w_j=\xi_j(\xi_{j-1}-\xi_{j+1})$ and
$\widetilde w_j^{(k)}=(\xi_{j-1}-\xi_{j+1})/(2-\xi_k-\xi_{k+1})$,
$j\leq k$. Each hyperprior gives full weak support. The last
column contains exact mean indices or upper bounds.}
\label{tab:survival-all}
\vspace{5pt}
\begin{tabular}{@{}lllll@{}}
\toprule
No. & $\mathcal L(K\mid\eta)$ & Survival $H_j$ & Hyperprior & $\E(K)$ or bound\\
\midrule
1 & $\operatorname{Geom}(p)$ & $(1-p)^{j-1}$ & $p\sim\Beta(a_1,b_1)$ & $1+\dfrac{b_1}{a_1-1}$\\
2 & $\mathcal L(1+Z_{4,p})$ & $\mathcal B_4(j-1;p)$ & $p\sim\Beta(a_1,b_1)$ & $1+\dfrac{4b_1}{a_1-1}$\\
3 & $\mathcal L(1+N_\nu)$ & $\mathcal P_\nu(j-1)$ & $\nu\sim\mathrm{Gamma}(a_1,b_1)$ & $1+\dfrac{a_1}{b_1}$\\
4 & $\mathcal L(1+Z_{2+28C,p})$ & $\dfrac{\mathcal B_2(j-1;p)+\mathcal B_{30}(j-1;p)}{2}$ & $p\sim\Beta(a_1,b_1)$ & $1+\dfrac{16b_1}{a_1-1}$\\
5 & $\mathcal L(1+N_{(1+4C)\nu})$ & $\dfrac{\mathcal P_\nu(j-1)+\mathcal P_{5\nu}(j-1)}{2}$ & $\nu\sim\mathrm{Gamma}(a_1,b_1)$ & $1+\dfrac{3a_1}{b_1}$\\
6 & $\operatorname{BetaGeom}(2,b)$ & $\dfrac{b(b+1)}{(b+j-1)(b+j)}$ & $b\sim\mathrm{Gamma}(a_1,b_1)$ & $1+\dfrac{a_1}{b_1}$\\
7 & $\operatorname{dLomax}(3,b)$ & $(1+y_j)^{-3}$ & $b-b_0\sim\mathrm{Gamma}(a_1,b_1)$ & $\leq1+m_b/2$\\
8 & $\operatorname{dWeibull}(2,b)$ & $e^{-y_j^2}$ & $b-b_0\sim\mathrm{Gamma}(a_1,b_1)$ & $\leq1+m_b\Gamma(3/2)$\\
9 & $\operatorname{dLN}(0,1;b)$ & $\overline\Phi(\log y_j)$ & $b-b_0\sim\mathrm{Gamma}(a_1,b_1)$ & $\leq1+m_b e^{1/2}$\\
10 & $\operatorname{dGamma}(2,1;b)$ & $e^{-y_j}(1+y_j)$ & $b-b_0\sim\mathrm{Gamma}(a_1,b_1)$ & $\leq1+2m_b$\\
11 & $\operatorname{Zipf}(\tau)$ & $\zeta(\tau,j)/\zeta(\tau)$ & $\tau-1\sim\Beta(a_1,b_1)$ & $\infty$\\
12 & $\operatorname{dPareto}(\rho)$ & $j^{-\rho}$ & $\rho\sim\mathrm{Gamma}(a_1,b_1)$ & $\infty$\\
13 & $\operatorname{logPareto}(\rho)$ & $(1+\log j)^{-\rho}$ & $\rho\sim\mathrm{Gamma}(a_1,b_1)$ & $\infty$\\
\bottomrule
\end{tabular}
\par\medskip
\raggedright
\noindent
In rows 4 and 5, $C\sim\operatorname{Bernoulli}(1/2)$
is drawn separately for each index.
In rows 1, 2 and 4, the displayed means require $a_1>1$
and are infinite otherwise.
In rows 7--10, $y_j=(j-1)/b$ and $m_b=b_0+a_1/b_1$,
with fixed $0\leq b_0<\infty$.
In row 9, $H_1=1$ and the formula applies for $j\geq2$.
Throughout, $a_1,b_1>0$, with gamma laws parametrized by
shape and rate.
\par
\endgroup
\end{table}

\par\medskip\Needspace{3\baselineskip}\noindent\textit{Verification of full weak support.}\enspace
The required small-maximum limits follow from
\[
 \Delta(p)\leq p\longrightarrow0
 \quad(p\downarrow0,\ \text{rows }1,2,4),
 \qquad
 \Delta(\nu)\longrightarrow0
 \quad(\nu\to\infty,\ \text{rows }3,5).
\]
For the scale families,
\[
 \begin{aligned}
 \Delta(b)&=\frac{2}{b+2} &&\text{in row }6,\\
 \Delta(b)&=1-(1+1/b)^{-3} &&\text{in row }7,\\
 \Delta(b)&\leq\frac{\|f_Y\|_\infty}{b}
     &&\text{in rows }8\text{--}10,
 \end{aligned}
\]
so $\Delta(b)\to0$ as $b\to\infty$. The last bound uses
the bounded density of $Y$ and the interval width $1/b$.
For the remaining families,
\[
 \begin{aligned}
 \Delta(\tau)&=\zeta(\tau)^{-1}\longrightarrow0
       &&(\tau\downarrow1,\ \text{row }11),\\
 \Delta(\rho)&=1-2^{-\rho}\longrightarrow0
       &&(\rho\downarrow0,\ \text{row }12),\\
 \Delta(\rho)&=1-(1+\log2)^{-\rho}\longrightarrow0
       &&(\rho\downarrow0,\ \text{row }13).
 \end{aligned}
\]
Each displayed hyperprior assigns positive probability
arbitrarily close to the relevant boundary, ensuring
(21) and hence full weak
support by Theorem~3.5.

\par\medskip\Needspace{3\baselineskip}\noindent\textit{Mean-index calculations.}\enspace
For $p\sim\Beta(a_1,b_1)$,
\[
 \E\!\left(\frac{1-p}{p}\right)
 =\begin{cases}
 b_1/(a_1-1),&a_1>1,\\
 \infty,&0<a_1\leq1.
 \end{cases}
\]
Combining this with $\E(K\mid p)=1+r(1-p)/p$, where
$r=1,4,16$, gives rows 1, 2 and 4.
The conditional means in rows 3, 5 and 6 are
$1+\nu$, $1+3\nu$ and $1+b$, respectively, so their
gamma hyperpriors give the stated entries.
For rows 7--10, one common bound suffices,
\[
 K=1+\lfloor bY\rfloor\leq1+bY,
 \qquad \E(K)\leq1+m_b\E(Y),
\]
with $\E(Y)=1/2,\Gamma(3/2),e^{1/2},2$, respectively.
In row 10, summing the survival probabilities also gives
\[
 \E(K\mid b)
 =\frac1{1-e^{-1/b}}
  +\frac{e^{-1/b}}{b(1-e^{-1/b})^2}.
\]
Finally, the hyperpriors in rows 11 and 12 satisfy
$\Pp(1<\tau<2)=1$ and $\Pp(\rho\leq1)>0$,
respectively, placing positive probability on an infinite
conditional mean. In row 13,
\[
 \E(K\mid\rho)=\sum_{j\geq1}(1+\log j)^{-\rho}
 =\infty,\qquad \rho>0.
\]
Thus all three unconditional means are infinite.

\par\medskip\Needspace{3\baselineskip}\noindent\textit{Defining-order size-biased invariance.}\enspace
The square-root weights are not ISBP in their defining
order, even after mixing over the index parameters.
Their first two stick-breaking fractions satisfy
\[
 v_1=w_1=\sqrt{\pi_1},\qquad
 v_2=\frac{w_2}{1-w_1}=\sqrt{\pi_1+\pi_2}>v_1
 \quad\text{almost surely}.
\]
For any infinite positive weight sequence, choose a mass
$x$ and another mass $y<x(1-x)$. A size-biased permutation
selects them first and second with probability
$xy/(1-x)>0$, giving fractions
$v'_1=x$ and $v'_2=y/(1-x)<v'_1$.
Since the event $v_2<v_1$ has zero probability under the
square-root construction, its ordered law cannot be ISBP
\citep{PitmanISBP1996}.
Jointly relabeling weights and atoms preserves $G$, but
can change the finite measures $G_k^\star$.

\par\medskip\Needspace{3\baselineskip}\noindent\textit{Finite evaluation.}\enspace
Given $K=k$, only $H_1,\ldots,H_{k+1}$ are needed to
evaluate $\widetilde w_j^{(k)}=a_j/t_k$.
Numerically stable expressions are
\[
 \xi_j=\frac{H_j}{1+\sqrt{1-H_j}},\qquad
 a_1=\sqrt{\pi_1},\qquad
 a_j=\frac{H_{j-1}-H_{j+1}}
 {\sqrt{1-H_{j+1}}+\sqrt{1-H_{j-1}}},\quad j\geq2.
\]
Within each SSP realization, draw an atom from $G_0$
when its label is first selected and reuse it thereafter.
For a prior $\eta\sim\Lambda$, the augmented update is
\[
 \Lambda(\dd\eta\mid\boldsymbol x,\zb,\kb,\thetab)
 \propto
 \left\{\prod_{i=1}^n\pi_{k_i}(\eta)
       \frac{a_{z_i}(\eta)}{t_{k_i}(\eta)}\right\}
 \Lambda(\dd\eta).
\]

\subsection{Prior-predictive calibration through finite allocations}
\label{sec:predictive-calibration}
Condition on the shared inputs $(\boldsymbol\pi,\boldsymbol a)$
of Theorem~3.1, and draw
$K_i\mid(\boldsymbol\pi,\boldsymbol a)
\stackrel{\mathrm{iid}}\sim\boldsymbol\pi$, $i\geq1$.
For realized indices $k_1,\ldots,k_n$, write
\[
 b_{ij}=\Pp(J_i=j\mid K_i=k_i,\boldsymbol a)
       =\frac{a_j}{s_{k_i}}\mathbf1\{j\leq k_i\}\quad\mbox{and}\quad
  k^\star=\max_{1\leq i\leq n}k_i.
\]

\begin{proposition}[Finite-index prior-predictive identities]
\label{prop:finite-predictive}
For every integer $r\geq2$,
\begin{equation}\label{eq:finite-r-tie}
 \Pp(J_1=\cdots=J_r\mid\boldsymbol\pi,\boldsymbol a)
 =\E\left\{
       \frac{\sum_{j\leq\min_{1\leq i\leq r}K_i}a_j^r}
            {\prod_{i=1}^r s_{K_i}}
       \,\middle|\,\boldsymbol\pi,\boldsymbol a
       \right\}
 =\sum_{j\geq1}w_j^r.
\end{equation}
Let $c_n=|\{J_1,\ldots,J_n\}|$ and let $M_{1,n}$ count
labels selected exactly once. Conditional on
$K_{1:n}=k_{1:n}$,
\begin{align}
 \E(c_n\mid K_{1:n}=k_{1:n},\boldsymbol a)
 &=\sum_{j=1}^{k^\star}
       \left\{1-\prod_{i=1}^n(1-b_{ij})\right\},
       \label{eq:finite-occupancy}\\
 \E(M_{1,n}\mid K_{1:n}=k_{1:n},\boldsymbol a)
 &=\sum_{j=1}^{k^\star}\sum_{i=1}^n
       b_{ij}\prod_{\substack{1\leq h\leq n\\h\ne i}}
       (1-b_{hj}).\label{eq:finite-singletons}
\end{align}
Applying the same definition of $b_{ij}$ to the next draw gives
\begin{equation}\label{eq:finite-discovery}
 \begin{aligned}
 &\Pp\!\left(J_{n+1}\notin\{J_1,\ldots,J_n\}
       \,\middle|\,K_{1:n+1}=k_{1:n+1},\boldsymbol a\right)\\
 &\qquad=\sum_{j=1}^{k_{n+1}}b_{n+1,j}
                 \prod_{i=1}^n(1-b_{ij}).
 \end{aligned}
\end{equation}
Averaging over the indices yields
(13) and
\[
 \E(M_{1,n}\mid\boldsymbol\pi,\boldsymbol a)
 =n\sum_{j\geq1}w_j(1-w_j)^{n-1}.
\]
\end{proposition}

\begin{proof}
Conditional independence gives the $r$-fold coincidence
probability $\sum_j\prod_{i=1}^rb_{ij}$.
For $n_j=\sum_{i=1}^n\mathbf1\{J_i=j\}$, it also gives
\[
 \begin{aligned}
 \Pp(n_j=0\mid K_{1:n}=k_{1:n},\boldsymbol a)
   &=\prod_{i=1}^n(1-b_{ij}),\\
 \Pp(n_j=1\mid K_{1:n}=k_{1:n},\boldsymbol a)
   &=\sum_{i=1}^nb_{ij}
       \prod_{\substack{1\leq h\leq n\\h\ne i}}(1-b_{hj}).
 \end{aligned}
\]
Summing the indicators in
$c_n=\sum_j\mathbf1\{n_j>0\}$ and
$M_{1,n}=\sum_j\mathbf1\{n_j=1\}$ proves the two count
identities. Multiplying the probability of $n_j=0$ by
$b_{n+1,j}$ and summing gives \eqref{eq:finite-discovery}.
Finally,
\[
 \E\!\left\{\frac{a_j}{s_{K_i}}\mathbf1\{j\leq K_i\}
       \,\middle|\,\boldsymbol\pi,\boldsymbol a\right\}
 =\sum_{k\geq j}\pi_k\frac{a_j}{s_k}=w_j.
\]
Independence of the indices therefore yields all marginal
formulas, with Tonelli's theorem justifying the nonnegative
interchanges.
\end{proof}

Drawing one index per observation, and hence two for pairwise
coincidence, allows each evaluation to use only finitely many
allocation coefficients. Averaging independent repetitions gives
unbiased Monte Carlo estimates with finite variance, since the
coincidence and discovery terms lie in $[0,1]$ and the occupancy
and singleton terms in $[0,n]$. These estimates retain finite
variance even when $\E(K)=\infty$, although their expected
computational work may then be infinite. Under a hyperprior,
each repetition begins with a single draw of
$(\boldsymbol\pi,\boldsymbol a)$, shared by all observations
in that simulated sample.

Writing $T_2=\sum_jw_j^2$ and
$d_n=\Pp(J_{n+1}\notin\{J_1,\ldots,J_n\})$,
the common sampling law gives, after averaging over any
shared hyperprior,
\[
 \begin{aligned}
 \E(c_{n+1})-\E(c_n)&=d_n,\\
 d_n-d_{n+1}
   &=\E\!\left\{\sum_jw_j^2(1-w_j)^n\right\}\geq0,\\
 \E(c_2)&=2-\E(T_2).
 \end{aligned}
\]
These identities restrict the combinations of summary values
that a prior can reproduce, and even attainable values may
correspond to several priors. In mixture models, the summaries
concern latent allocations, so assessing their implications
for observed data also requires considering $G_0$ and the
sampling density $f(\cdot\mid\theta)$.

\par\medskip\Needspace{3\baselineskip}\noindent\textit{A shared hyperprior on coincidence.}\enspace
With $K-1\mid p\sim\NegBin(2,p)$ and uniform finite
allocation, Theorem~3.1 gives
$w_j=p(1-p)^{j-1}$. Consequently,
\[
 T=\Pp(J_1=J_2\mid p)
  =\sum_jw_j^2=\frac{p}{2-p},\qquad
 p=\frac{2T}{1+T},\quad\mbox{and}\quad \E(K\mid T)=\frac{1}{T}.
\]
Choosing $T\sim\Beta(a,b)$ therefore specifies uncertainty
directly about the coincidence probability of the realized
weights, with
\[
 \E(T)=\frac{a}{a+b}\quad\mbox{and}\quad 
 \Var(T)=\frac{ab}{(a+b)^2(a+b+1)}.
\]
Since $\sup_jw_j=p\to0$ as $T\downarrow0$, every
$a,b>0$ gives full weak support by
(18). Its mean index is
\[
 \E(K)=\E(T^{-1})
 =\begin{cases}
 (a+b-1)/(a-1),&a>1,\\
 \infty,&0<a\leq1.
 \end{cases}
\]

For the augmented update, put $L=\sum_{i=1}^n(k_i-1)$.
Since $\pi_k(p)=kp^2(1-p)^{k-1}$ and
$\widetilde w_j^{(k)}=1/k$, the parameter-dependent
likelihood factor is
\[
 \prod_{i=1}^n\frac{\pi_{k_i}(p)}{k_i}
 =p^{2n}(1-p)^L
 \propto T^{2n}(1-T)^L(1+T)^{-(2n+L)}.
\]
Multiplying by the beta prior density gives the posterior density
\begin{equation}\label{eq:coincidence-posterior}
 f_T(t\mid\bm x,\zb,\kb,\thetab)
 \propto t^{a+2n-1}(1-t)^{b+L-1}(1+t)^{-(2n+L)},
 \qquad 0<t<1.
\end{equation}
The component likelihood is constant in $t$ under this
conditioning. This updates the elicited prior directly in
$T$, without requiring its induced prior on $p$ to be Beta.

\par\medskip\Needspace{3\baselineskip}\noindent\textit{Comparisons at a common coincidence probability.}\enspace
Table~1 matches
$\Pp(J_1=J_2)=0.2$ across priors at fixed parameter values.
For the DP,
\[
 \Pp(J_1=J_2)=\E(T_2)=\frac1{\alpha+1},
\]
so the matching concentration is $\alpha=4$.
Its prior-predictive summaries are
\[
 \E(c_n)=\sum_{i=0}^{n-1}\frac{\alpha}{\alpha+i},\qquad
 \E(M_{1,n})=\frac{n\alpha}{\alpha+n-1},\quad\mbox{and}\quad 
 d_n=\frac{\alpha}{\alpha+n}.
\]
The DP predictive rule gives $d_n$, whose sum gives
$\E(c_n)$. By exchangeability, the singleton expectation
is $n$ times
\[
 \Pp(\{1\}\text{ is a block of }\mathcal P_n)
 =\prod_{i=1}^{n-1}\frac{\alpha+i-1}{\alpha+i}
 =\frac{\alpha}{\alpha+n-1}.
\]

For geometric weights, the matching value is $p=1/3$.
At $r=3$ in Example~3.6,
\[
 T_2(p)=p^4\sum_{j\geq1}j^2(1-p)^{2(j-1)}
       =\frac{p\{1+(1-p)^2\}}{(2-p)^3},
\]
so $T_2(p)=0.2$ gives $p\approx0.52$.
For the shifted-Poisson construction, solving
\[
 T_2(\nu)=\sum_{j\geq1}
       \left\{\frac{\Pp(N_\nu\geq j)}{\nu}\right\}^2=0.2
\]
gives $\nu\approx3.53$.
The remaining summaries are evaluated from their weight
sums, retaining the original normalization.

For these three directly specified weight sequences, we use $H=500$ terms. The exact residual
masses are
\[
 R_H=
 \begin{cases}
 (1-p)^H,&\text{geometric},\\
 (1-p)^H(1+Hp),&r=3,\\
 \Pp(N_\nu\geq H)-\dfrac{H}{\nu}\Pp(N_\nu\geq H+1),
     &\text{shifted Poisson}.
 \end{cases}
\]
Write $C_n=\E(c_n)$ and $S_n=\E(M_{1,n})$, and use
an additional subscript $H$ for each sum restricted to
$j\leq H$. The omitted contributions satisfy
\[
 \begin{aligned}
 0\leq T_2-T_{2,H}&\leq R_H^2,
 &0\leq d_n-d_{n,H}&\leq R_H,\\
 0\leq C_n-C_{n,H}&\leq nR_H,
 &0\leq S_n-S_{n,H}&\leq nR_H.
 \end{aligned}
\]
For $H>\nu$, the Chernoff inequality gives
\[
 \log R_H\leq-\nu+H\log(e\nu/H),
\]
providing a bound on the omitted Poisson mass even when the
tail probability is too small to represent numerically.
For $n=50$ and the parameter values used in the comparison,
all bounds on omitted contributions are below $10^{-80}$.
The calibrated coincidence probabilities differ from the
target $0.2$ by less than $10^{-11}$.

Under a shared hyperprior, prior-predictive assessment averages
the conditional summaries over the parameter distribution.
For example, $T\sim\Beta(2,8)$ and fixed $p=1/3$ both yield
$\Pp(J_1=J_2)=0.2$, although the beta hierarchy allows the
coincidence probability to vary across realizations of the
random measure. Agreement at the  mean level therefore need not
extend to occupancy, singleton counts or discovery.

\section{Sampling identities and finite-measure comparisons}\label{app:sampling}
We begin with the joint sampling law and its decomposition in
terms of additional visibility, then examine how the finite
conditional measures relate to the original SSP and characterize
the index distribution through its moments and tail probabilities.

\subsection{The sample law and additional visibility}
\label{sec:notation}\label{sec:supporting}\label{sec:layers}

The auxiliary index $K$ determines the labels available to an
individual draw, rather than a number of mixture components shared
by the sample. This is important to distinguish this index from the classical mixture order in finite mixture modeling. Although observations from an SSP mixture share
the same random weights and atoms, each requires its own index
to preserve the original conditional sampling law. The distinction
becomes clear by considering what happens when two draws use
a common index.

Fix the weights, atoms and schedule $(\wb,\thetab,\xib)$, and
write $q_k=G_k^\star(A)$ for a measurable set $A$. Suppose that
a single index $K$ is drawn with probabilities $(p_k)$, followed
by two conditionally independent observations
\[
 Y_1,Y_2\mid K,\wb,\thetab,\xib
 \stackrel{\mathrm{iid}}{\sim}G_K^\star.
\]
Although each observation has marginal probability
$G(A)=\sum_kp_kq_k$ of belonging to $A$, their joint
probability is
\begin{equation}\label{eq:shared-K-joint}
\begin{aligned}
 \Pp(Y_1\in A,Y_2\in A\mid\wb,\thetab,\xib)
 &=\sum_{k\geq1}p_kq_k^2,\\
 G(A)^2&=\left(\sum_{k\geq1}p_kq_k\right)^2.
\end{aligned}
\end{equation}
The difference reflects the dependence introduced by sharing
$K$. Indeed, conditional independence given $K$ and the law
of total covariance yield
\begin{equation}\label{eq:shared-K-covariance}
\begin{aligned}
 \Cov\{\mathbf1_A(Y_1),\mathbf1_A(Y_2)\mid\wb,\thetab,\xib\}
 &=\sum_{k\geq1}p_kq_k^2-G(A)^2\\
 &=\Var\{G_K^\star(A)\mid\wb,\thetab,\xib\}.
\end{aligned}
\end{equation}
This covariance is strictly positive whenever $q_k\ne q_\ell$
for some $k,\ell$ with $p_kp_\ell>0$.

A simple example illustrates the effect. Take two distinct atoms
with $w_1=w_2=1/2$, $w_j=0$ for $j>2$, and
$\xi_j=2^{1-j}$. For $A=\{\theta_1\}$, we have
$p_1=1/4$, $q_1=1$, and $q_k=1/3$ for $k\geq2$,
whose combined index probability is $3/4$. Consequently,
\[
 \Pp(Y_1\in A,Y_2\in A\mid\wb,\thetab,\xib)
 =\frac14+\frac34\left(\frac13\right)^2
 =\frac13,
\]
whereas $G(A)^2=1/4$. Sharing the index therefore introduces
an additional covariance of $1/12$.

To recover the original SSP sampling law, draw independent
indices $K_1,\ldots,K_n$ conditional on the shared weights
and schedule, and then sample each $Y_i$ independently from
$G_{K_i}^\star$. For measurable sets $A_1,\ldots,A_n$,
\begin{equation*}
\begin{aligned}
 &\Pp(Y_i\in A_i\text{ for all }i\mid\wb,\thetab,\xib)\\
 &\quad=\sum_{k_1,\ldots,k_n}
       \prod_{i=1}^n p_{k_i}G_{k_i}^\star(A_i)\\
 &\quad=\prod_{i=1}^n
       \left\{\sum_{k\geq1}p_kG_k^\star(A_i)\right\}
       =\prod_{i=1}^nG(A_i).
\end{aligned}
\end{equation*}
Thus, independent observation-specific indices preserve
conditional independence given $G$ and, after averaging over
the random measure, the usual SSP dependence. In particular,
the covariance of two event indicators is $\Var\{G(A)\}$,
whereas a shared index would add
$\E[\Var\{G_K^\star(A)\mid\wb,\thetab,\xib\}]$.

\par\medskip\Needspace{3\baselineskip}\noindent\textit{Stop--continue layers.}\enspace
The decomposition $K=J+S$ also represents $G$ through the
additional visibility $S=K-J$. For integers $t\geq0$ and
$A\in\mathcal B$, define
\begin{equation*}
 \begin{aligned}
 \mathcal C_t(A)
 &=\Pp(\theta_J\in A,S\geq t\mid\wb,\thetab,\xib)
   =\sum_{j\geq1}w_j\frac{\xi_{j+t}}{\xi_j}\delta_{\theta_j}(A),\\
 \mathcal M_t(A)
 &=\Pp(\theta_J\in A,S=t\mid\wb,\thetab,\xib)
   =\mathcal C_t(A)-\mathcal C_{t+1}(A).
 \end{aligned}
\end{equation*}
These are subprobability measures, not necessarily supported
on finitely many atoms. They describe the joint distribution
of the selected atom and the event that at least $t$
additional positions remain visible, or exactly $t$ in
the case of $\mathcal M_t$. In particular,
\begin{equation*}
 \mathcal C_t(\mathsf X)
 =\Pp(S\geq t\mid\wb,\xib),
 \qquad
 \mathcal M_t(\mathsf X)
 =\Pp(S=t\mid\wb,\xib).
\end{equation*}

Since $\mathcal C_t=\mathcal M_t+\mathcal C_{t+1}$
and $\mathcal C_0=G$, successive applications give,
for every integer $r\geq1$,
\begin{equation*}
 G=\sum_{t=0}^{r-1}\mathcal M_t+\mathcal C_r.
\end{equation*}
Moreover, $\mathcal C_r(A)\leq\mathcal C_r(\mathsf X)$
for every measurable $A$. Since
$\xi_{j+r}/\xi_j\to0$ for each fixed $j$, dominated
convergence under the probability weights $(w_j)$ gives
$\mathcal C_r(\mathsf X)\to0$. Consequently,
\begin{equation*}
 G(A)=\sum_{t\geq0}\mathcal M_t(A),
 \qquad A\in\mathcal B.
\end{equation*}
For a constant-ratio schedule $\xi_j=cq^{j-1}$, this
decomposition simplifies to
$\mathcal C_t=q^tG$ and
$\mathcal M_t=(1-q)q^tG$.

\subsection{Total-variation bounds}\label{sec:tv}\label{sec:tilting}\label{sec:tvbounds}

The exact two-stage construction averages over $K$.  Here we compare a particular conditional component $G_k^\star$, under the natural coupling, with the original realized measure $G$.

For $s_k>0$, let
\begin{equation*}
R_k=\sum_{j>k}w_j,
\qquad
G_{\mathrm{ren}}^{(k)}=\sum_{j=1}^k\frac{w_j}{1-R_k}\delta_{\theta_j}.
\end{equation*}

\par\medskip\Needspace{3\baselineskip}\noindent\textit{A tilting inequality.}\enspace
The following sharp bound specializes the likelihood-ratio-range inequality of \citet[Corollary~5]{Binette2019}. A short proof is included.

\begin{lemma}[Sharp bounded-tilt inequality]\label{lem:tilt}
Let $P$ be a probability mass function and $0<a\leq f\leq b<\infty$. Define
\begin{equation}\label{eq:tilt}
Q(i)=f(i)P(i)/\E_Pf.
\end{equation}
Then
\begin{equation}\label{eq:tilt-bound}
\TV(P,Q)\leq\frac{\sqrt b-\sqrt a}{\sqrt b+\sqrt a}.
\end{equation}
The bound is sharp uniformly over $P$ and $f$ with these bounds.
\end{lemma}
\begin{proof}
For $a=b$ the assertion is immediate. Otherwise put $\mu=\E_Pf$. Since $\E_P(f-\mu)=0$,
\[
\TV(P,Q)=\frac{\E_P(f-\mu)_+}{\mu}.
\]
Convexity gives $(f-\mu)_+\leq(f-a)(b-\mu)/(b-a)$ on $[a,b]$, and hence
\[
\TV(P,Q)\leq\frac{(\mu-a)(b-\mu)}{(b-a)\mu}.
\]
The right-hand side is maximized at $\mu=\sqrt{ab}$, yielding \eqref{eq:tilt-bound}. Equality is attained when $f$ takes only the values $a,b$, with $P(f=b)=\sqrt a/(\sqrt a+\sqrt b)$.
\end{proof}
The coarser inequality $(b-a)/(b+a)$ is also valid but is not generally sharp.

\par\medskip\Needspace{3\baselineskip}\noindent\textit{Pathwise and expected bounds.}\enspace
Put
\begin{equation*}
M_k=\frac{\xi_1}{\xi_k},
\qquad
D_k=\frac{\sqrt{M_k}-1}{\sqrt{M_k}+1}.
\end{equation*}

\begin{proposition}\label{prop:tv}
For every $k\geq1$ with $s_k>0$, the finite reweighted
measure satisfies
\begin{equation*}
 \TV(G_{\mathrm{ren}}^{(k)},G_k^\star)\leq D_k,
 \qquad
 \TV(G,G_k^\star)\leq R_k+(1-R_k)D_k.
\end{equation*}
These bounds hold regardless of whether atom locations are
distinct. If $K$ follows the conditional law in (3), averaging
over its distribution gives
\begin{equation}\label{eq:conditional-expected-tv}
 \E_{K\mid\wb,\xib}\{\TV(G,G_K^\star)\}
 \leq
 \E_{K\mid\wb,\xib}\{R_K+(1-R_K)D_K\}.
\end{equation}

When the schedule is independent of the weights, conditioning
on its realized value yields
\begin{align}
 \Pp(K=k\mid\xib)
 &=(\xi_k-\xi_{k+1})
   \sum_{h=1}^k\frac{\E(w_h)}{\xi_h},
 \label{eq:marginalK}\\
 \E(R_K\mid\xib)
 &=\sum_{1\leq h<j}
   \left(1-\frac{\xi_j}{\xi_h}\right)\E(w_hw_j).
 \label{eq:ERK}
\end{align}
For a random independent schedule, further averaging over
$\xib$ gives the unconditional quantities. In particular,
replacing $(1-R_K)D_K$ by $D_K$ provides the simpler bound
\begin{equation*}
 \E\{\TV(G,G_K^\star)\}
 \leq
 \E_{\xib}\!\left[
 \sum_{h<j}\left(1-\frac{\xi_j}{\xi_h}\right)\E(w_hw_j)
 +\sum_{k\geq1}D_k\Pp(K=k\mid\xib)
 \right].
\end{equation*}
A tighter bound follows by retaining the factor $1-R_K$
in \eqref{eq:conditional-expected-tv}. For schedules depending
on the weights, the corresponding expectations must instead
be evaluated under their joint law.
\end{proposition}

\begin{proof}
Fix the weights, atoms and schedule. For $s_k>0$, define
$P(j)=w_j/(1-R_k)$ and $f(j)=1/\xi_j$ on
$\{1,\ldots,k\}$. The tilted probabilities in
\eqref{eq:tilt} then satisfy
$Q(j)=\widetilde w_j^{(k)}$. Since
$f(j)\in[1/\xi_1,1/\xi_k]$, Lemma~\ref{lem:tilt}
bounds their total-variation distance by $D_k$.
Passing from labels to atom locations cannot increase
this distance, so the first inequality follows even
when some atoms coincide.

For $R_k>0$, write
\[
 G=(1-R_k)G_{\mathrm{ren}}^{(k)}
   +R_kG_{\mathrm{tail}}^{(k)},
\]
where $G_{\mathrm{tail}}^{(k)}$ is the normalized
remaining mass. Convexity of total variation gives
\[
 \begin{aligned}
 \TV(G,G_k^\star)
 &\leq(1-R_k)
   \TV(G_{\mathrm{ren}}^{(k)},G_k^\star)
   +R_k\TV(G_{\mathrm{tail}}^{(k)},G_k^\star)\\
 &\leq(1-R_k)D_k+R_k.
 \end{aligned}
\]
The same bound follows directly from the first inequality
when $R_k=0$. Averaging over the conditional distribution
of $K$ proves \eqref{eq:conditional-expected-tv}.

To obtain the moment identities, Tonelli's theorem and
$\sum_{k=h}^{j-1}(\xi_k-\xi_{k+1})=\xi_h-\xi_j$
give, for fixed $(\wb,\xib)$,
\begin{align}
 \E(R_K\mid\wb,\xib)
 &=\sum_{k\geq1}(\xi_k-\xi_{k+1})
   \sum_{h\leq k}\sum_{j>k}\frac{w_hw_j}{\xi_h}
   \notag\\
 &=\sum_{h<j}
   \left(1-\frac{\xi_j}{\xi_h}\right)w_hw_j.
 \label{eq:pathwise-ERK}
\end{align}
When $\xib$ is independent of $\wb$, taking expectations
over the weights conditional on the schedule yields
\eqref{eq:ERK}. Applying the same argument to $p_k$
gives \eqref{eq:marginalK}. The tower property then
establishes the unconditional total-variation bound.

For an endogenous schedule, the pathwise identity
\eqref{eq:pathwise-ERK} remains valid, although its
expectation generally cannot be expressed solely through
the unconditional moments $\E(w_hw_j)$.
\end{proof}

For $\xi_j=cq^{j-1}$,
\begin{equation*}
D_k=\frac{1-q^{(k-1)/2}}{1+q^{(k-1)/2}}
=\tanh\!\left\{\frac{(k-1)|\log q|}{4}\right\}.
\end{equation*}
For the tail-mass schedule $\xi_k=R_{k-1}=\prod_{\ell<k}(1-v_\ell)$,
\begin{equation*}
D_k=\frac{1-\sqrt{R_{k-1}}}{1+\sqrt{R_{k-1}}}.
\end{equation*}
As $k$ increases, these bounds may become less informative,
reflecting the fact that the finite reweighted measures need
not converge to $G$, as established in Theorem~4.2 of the
main paper.

\begin{corollary}[A universal mean residual bound]\label{cor:mean-residual}
For every proper weight sequence and every admissible schedule, including endogenous schedules,
\begin{equation*}
\E(R_K\mid\wb,\xib)\leq\frac12\left(1-\sum_jw_j^2\right)<\frac12.
\end{equation*}
\end{corollary}
\begin{proof}
In \eqref{eq:pathwise-ERK}, $0\leq1-\xi_j/\xi_h\leq1$ for $h<j$. Also $\sum_{h<j}w_hw_j=(1-\sum_jw_j^2)/2$. At least one weight is positive, so $\sum_jw_j^2>0$.
\end{proof}
The bound controls the omitted mass of the original SSP and
cannot be improved uniformly over schedules, as illustrated
by constant-ratio schedules with $q\downarrow0$.

For the DP under the tail-mass schedule, the reweighted
coefficients $a_j=v_j$ are i.i.d.\ $\Beta(1,\alpha)$, so
$s_\infty=\infty$ almost surely. With the same schedule,
geometric stick-breaking gives $a_j=v>0$ and the empirical
measure $G_k^\star=k^{-1}\sum_{j=1}^k\delta_{\theta_j}$, again
with $s_\infty=\infty$. For distinct atoms, both constructions
therefore satisfy $\TV(G_k^\star,G)\to1$, while averaging over
the index still recovers $G$ exactly through
$\sum_{k\geq1}p_kG_k^\star=G$.

\subsection{Index tails and moment identities}\label{sec:tails}\label{sec:mean-dimension}
For exogenous schedules, averaging over the SSP weight prior
expresses the tail probabilities of the representation index
in terms of the mean weights, yielding explicit formulas for
several familiar stick-breaking priors.

\par\medskip\Needspace{3\baselineskip}\noindent\textit{Tail behavior of the representation index.}\enspace
Define, by continuous extension at $x=y$,
\begin{equation*}
\Delta_k(x,y)=
\begin{cases}
(x^k-y^k)/(x-y),&x\neq y,\\
ky^{k-1},&x=y.
\end{cases}
\end{equation*}

\begin{proposition}\label{prop:regularvar}
Suppose $\xi_j=cq^{j-1}$ with $q\in(0,1)$ and $\E(w_k)\sim Ck^{-p}$ for $C>0$ and $p>1$.  Then
\begin{equation*}
\Pp(K=k)\sim Ck^{-p},
\end{equation*}
so $\E(K)<\infty$ if and only if $p>2$.
\end{proposition}

\begin{proof}
From \eqref{eq:marginalK},
\[
\Pp(K=k)=(1-q)\sum_{m=0}^{k-1}\E(w_{k-m})q^m.
\]
For any fixed $M$, regular variation gives
$\E(w_{k-m})/\E(w_k)\to1$ uniformly over $0\leq m\leq M$.
For $M<m\leq k/2$, these ratios remain uniformly bounded,
so their contribution is controlled by the summable geometric
tail. When $m>k/2$, the factor $q^m$ makes the contribution
exponentially small relative to $\E(w_k)\asymp k^{-p}$.
Letting $k\to\infty$ and then $M\to\infty$ therefore yields
\[
 \frac{\Pp(K=k)}{\E(w_k)}
 \longrightarrow (1-q)\sum_{m\geq0}q^m=1.
\]
The condition for a finite expected index follows by comparing
$\sum_k k\Pp(K=k)$ with $\sum_k k^{1-p}$.
\end{proof}

\par\medskip\Needspace{3\baselineskip}\noindent\textit{Dirichlet process.}\enspace
For a Dirichlet process with total mass $\alpha$, put $a=\alpha/(\alpha+1)$.  Since $\E(w_k)=(1-a)a^{k-1}$,
\begin{equation*}
\Pp(K=k)=\frac{1-q}{\alpha+1}\Delta_k(q,a).
\end{equation*}
At $q=a$ this becomes $k a^{k-1}/(\alpha+1)^2$.  The expected tail term in \eqref{eq:ERK} is
\begin{equation*}
\E(R_K)=\frac{\alpha}{2}\frac{1-q}{(\alpha+1)-\alpha q}.
\end{equation*}
A matched-rate calibration $q\approx a$ keeps the schedule decay near the prior mean decay of the weights.  Note that this calibrates computation while preserving the target SSP.

\par\medskip\Needspace{3\baselineskip}\noindent\textit{Geometric stick-breaking.}\enspace
Let $w_k=V(1-V)^{k-1}$ and $a_V=1-V$.  Conditionally on $V$,
\begin{equation*}
\Pp(K=k\mid V)=V(1-q)\Delta_k(q,a_V).
\end{equation*}
When $q=a_V$, this reduces to
\begin{equation*}
\Pp(K=k\mid V)=kV^2(1-V)^{k-1},
\qquad K-1\sim\NegBin(2,V),
\end{equation*}
where the negative-binomial variable counts failures before two successes.  Consequently
\begin{equation*}
\E(K\mid V)=\frac{2}{V}-1,
\qquad
\E(R_K\mid V)=\frac{1-V}{(2-V)^2}.
\end{equation*}
If $V\sim\Beta(a_0,b_0)$,
\begin{equation*}
\Pp(K=k)=k\frac{B(a_0+2,b_0+k-1)}{B(a_0,b_0)}.
\end{equation*}

\par\medskip\Needspace{3\baselineskip}\noindent\textit{Pitman--Yor process.}\enspace
For discount $\sigma\in(0,1)$ and strength $\alpha>-\sigma$,
\begin{equation}\label{eq:PYmeanw}
\E(w_k)=
\frac{1-\sigma}{\alpha+1+(k-1)\sigma}
\prod_{\ell=1}^{k-1}
\frac{\alpha+\ell\sigma}{\alpha+1+(\ell-1)\sigma}
\sim C_{\sigma,\alpha}k^{-1/\sigma}.
\end{equation}
Thus Proposition~\ref{prop:regularvar} gives
\begin{equation*}
\Pp(K=k)\sim C_{\sigma,\alpha}k^{-1/\sigma},
\qquad
\E(K)<\infty\quad\Longleftrightarrow\quad\sigma<\tfrac12.
\end{equation*}
This tail behavior helps explain why generating PY weights
in their original stick-breaking order can be computationally
demanding. Under a fixed geometric schedule, $K$ remains finite almost surely but has infinite prior mean for $\sigma\geq1/2$.
Since $k^\star=\max_i k_i\geq k_1$, the number of labels
required by a sample also has infinite prior expectation.
The resulting need to generate potentially long sequences
of weights offers an explanation for the computational
difficulties of PY slice samplers reported by
\citet{CanaleCorradinNipoti2022}. Clearly, this connection concerns
representation apriori, while posterior running
times and mixing depend additionally on the likelihood.

\par\medskip\Needspace{3\baselineskip}\noindent\textit{Moment identities and the DP case.}\enspace
For a fixed constant-ratio schedule, Tonelli's theorem and $K=J+S$ give, allowing $+\infty$,
\begin{equation*}
\E(K\mid\wb)=\sum_{j\geq1}jw_j+\frac{q}{1-q}.
\end{equation*}
For a general realized schedule,
\[
\E(S\mid J=j,\xib)=\sum_{t\geq1}\frac{\xi_{j+t}}{\xi_j}.
\]
For PY weights with $0\leq\sigma<1/2$,
\begin{equation}\label{eq:PY-exact-EK}
\E(J)=\frac{\alpha+1-\sigma}{1-2\sigma},\qquad
\E(K)=\frac{\alpha+1-\sigma}{1-2\sigma}+\frac{q}{1-q}.
\end{equation}
To prove this, let $\overline R_k=\E(R_k)$, with $\overline R_0=1$, and put $u_k=(\alpha+1-\sigma+k\sigma)\overline R_k$. Independence of the PY sticks gives
\[
\overline R_{k+1}=\overline R_k\frac{\alpha+(k+1)\sigma}{\alpha+1+k\sigma},
\qquad u_k-u_{k+1}=(1-2\sigma)\overline R_k.
\]
For $\sigma<1/2$, $u_k\to0$ (geometrically for $\sigma=0$, and by the product asymptotics for $\sigma>0$). Summing proves \eqref{eq:PY-exact-EK}, since $\E(J)=\sum_{k\geq0}\overline R_k$. For $\sigma\geq1/2$, $K$ is still finite almost surely but has infinite prior mean under a fixed geometric schedule. 

For a DP with $v_j\stackrel{\mathrm{iid}}\sim\Beta(1,\alpha)$,
consider the tail-mass schedule $\xi_j=\prod_{\ell<j}(1-v_\ell)$ and put $\rho=\alpha/(\alpha+1)$. The resulting marginal index law and its expectation are
\begin{equation*}
 \Pp(K=k)=\frac{k+1}{(\alpha+1)(\alpha+2)}\rho^{k-1},
 \qquad
 \E(K)=\frac{2(\alpha+1)^2}{\alpha+2}.
\end{equation*}
These expressions follow by averaging
$p_k=\sum_{h=1}^k w_kv_h$ over the independent beta
stick-breaking fractions, whose moments give
\[
 \E(w_kv_h)
 =\frac{\rho^{k-1}}{(\alpha+1)(\alpha+2)}
 \begin{cases}
 1,&1\leq h<k,\\
 2,&h=k.
 \end{cases}
\]
Summing over $h$ gives the stated index law, whose first
moment yields $\E(K)$. Since the weights and the tail-mass
schedule are functions of the same fractions, this averaging
accounts for their dependence.

\subsection{Recursive reweighting and fixed-level illustrations}\label{sec:recursive-detail}
\par\medskip\Needspace{3\baselineskip}\noindent\textit{Proof of Proposition~4.1.}\label{sec:oneatom-proof}\enspace
\begin{proof}
For $k\geq2$, separating the last term in the defining sum gives
\[
 \begin{aligned}
 G_k^\star
 &=\sum_{j=1}^{k-1}\frac{a_j}{s_k}\delta_{\theta_j}
   +\frac{a_k}{s_k}\delta_{\theta_k}
 =\frac{s_{k-1}}{s_k}G_{k-1}^\star
   +\lambda_k\delta_{\theta_k}
  =(1-\lambda_k)G_{k-1}^\star
   +\lambda_k\delta_{\theta_k}.
 \end{aligned}
\]
Integrating this identity against $h$ gives the functional
recursion. Moreover, for $1\leq j\leq k$,
\[
 \lambda_j\prod_{\ell=j+1}^k(1-\lambda_\ell)
 =\frac{a_j}{s_j}
  \prod_{\ell=j+1}^k\frac{s_{\ell-1}}{s_\ell}
 =\frac{a_j}{s_k}
 =\widetilde w_j^{(k)},
\]
with empty products equal to $1$. Finally, subtracting
$G_{k-1}^\star$ from the first identity yields
\[
 \TV(G_k^\star,G_{k-1}^\star)
 =\lambda_k\TV(\delta_{\theta_k},G_{k-1}^\star)
 \leq\lambda_k.
\]
When the atoms are distinct, $G_{k-1}^\star(\{\theta_k\})=0$,
so $\TV(\delta_{\theta_k},G_{k-1}^\star)=1$ and equality holds.
\end{proof}

Once $a_k$ is available, the recursion updates $G_k^\star h$
from its previous value and $h(\theta_k)$, covering event
probabilities and finite-mixture density evaluations whenever
the required kernel values are finite. The reverse product
shows how an atom's share changes as later labels become
visible, with the joint prior of $(\lambda_j)$ induced by the
original weights and schedule.

For proper stick-breaking weights $w_j=v_jR_{j-1}$, with
$R_j=\prod_{\ell\leq j}(1-v_\ell)$ and $0<v_j<1$ almost surely
for every $j$, the \emph{tail-mass schedule}
$\xi_j=R_{j-1}$
is positive, strictly decreasing and converges to zero. It
reduces the reweighted coefficients to $a_j=v_j$, giving
$p_k=w_k\sum_{j\leq k}v_j$ and
$G_k^\star=\sum_{j\leq k}v_j\delta_{\theta_j}/\sum_{j\leq k}v_j$.
The same stick-breaking fractions generate the overshoot through successive
stopping probabilities $v_J,v_{J+1},\ldots$. In the geometric
case $v_j\equiv v$, the finite probabilities become $1/k$ and
$p_k=kv^2(1-v)^{k-1}$, recovering
Example~3.6 with $r=2$.

More generally, a sequential $J$-first sampler is available
when successive stick-breaking fractions can be drawn from their conditional
laws given the quantities already generated. The first $r$
stick-breaking fractions determine $w_1,\ldots,w_r$ and $R_r$, so the next
draw $v_{r+1}$ extends the realization through
$w_{r+1}=v_{r+1}R_r$ and $R_{r+1}=(1-v_{r+1})R_r$.
The sampler uses an independent $U_J\sim\mathrm{Unif}(0,1)$
to find $J=\min\{j:R_j\leq U_J\}$, then samples $S$ using
(46), extending the same realization through
$K=J+S$ as needed. Previously generated weights and atoms
remain fixed, while each new atom is drawn independently
from $G_0$. Provided each schedule comparison can be resolved
from finitely many generated quantities, directly or through
rigorous bounds, properness and $\xi_j\downarrow0$ ensure
almost-sure termination.

\par\medskip\Needspace{3\baselineskip}\noindent\textit{Proof of Theorem~4.2.}\label{sec:fixed-limit-proof}\enspace
\begin{proof}
For finite $s_\infty$, $G_k^\star$ is the normalized restriction of $G_\xi$ to its first $k$ atoms. For infinite $s_\infty$, the mass $G_k^\star(\{\theta_1,\ldots,\theta_r\})=s_r/s_k$ tends to zero for every fixed $r$, whereas $G$ assigns these sets masses increasing to one. This proves both limits. The equality condition follows by comparing each positive atomic mass.
\end{proof}

\par
\begin{figure}[tbp]\centering
\includegraphics[width=.88\linewidth]{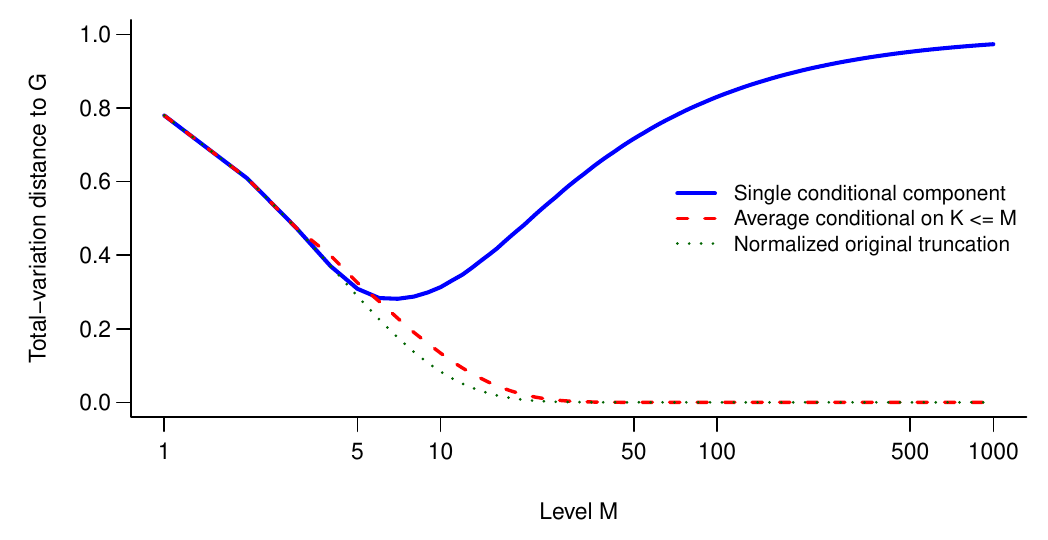}
\caption{Exact total-variation distances under one geometric weight sequence. A single finite conditional component $G_M^\star$ eventually becomes increasingly different from $G$, whereas the conditional average over $K\leq M$ approaches $G$. The normalized original truncation is included as a separate comparison and has distance $(0.78)^M$. }\label{fig:averaging-level}
\end{figure}

\medskip\Needspace{3\baselineskip}\noindent\textit{Averaging and fixing the finite level.}\label{sec:illustration-recursion}\label{sec:averaging-illustration}\enspace
For $w_j=v(1-v)^{j-1}$ with $v=0.22$ and the tail-mass schedule, $G_k^\star$ is uniform on its first $k$ atoms and $p_k=kv^2(1-v)^{k-1}$. To display the  identity alongside the fixed-level limit, define
\begin{equation*}
 \overline G_M=\frac{\sum_{k=1}^M p_kG_k^\star}{Q_M},\qquad
 Q_M=\Pp(K\leq M)=1-(1-v)^M(1+Mv).
\end{equation*}
This is the conditional average given $K\leq M$, used here only to visualize convergence of the averaging identity. Its masses are $\{w_j-v(1-v)^M\}/Q_M$ for $j\leq M$. The complete identity gives $\TV(G,\overline G_M)\leq1-Q_M$ by separating the remaining mixture mass. Figure~\ref{fig:averaging-level} compares that conditional average with a single $G_M^\star$ and the normalized original truncation $G_{\mathrm{ren}}^{(M)}$. Averaging approaches $G$, whereas a fixed component eventually has distance approaching one. Thus the effect of deterministic extension is materially different from the exact randomized representation.

\section{Posterior computation}\label{app:posterior}
We derive the parameter-specific posterior updates from the
augmented likelihood and examine their computational implications,
including the relationship between discrete and continuous
augmentations and the number of labels required for allocation
updates. We also provide the Gaussian reference calculations
and details of the numerical implementation.

\subsection{Full conditional distributions}\label{sec:conditionals}\label{sec:latent-cond}\label{sec:exo-cond}\label{sec:endo-cond}\label{sec:gsb-cond}\label{sec:reverse-posterior}
Combining the augmented likelihood in (49)
with the original priors on the weights, schedule and atoms
give the posterior distribution, with $z_i\leq k_i$ for every
observation. Conditional on the allocations, atoms with no
assigned observations retain their prior distribution $G_0$.

\par\medskip\Needspace{3\baselineskip}\noindent\textit{Component parameters and latent variables.}\enspace
For every instantiated component,
\begin{equation*}
p(\theta_j\mid\cdots)\propto p(\theta_j)\prod_{i:z_i=j}f(x_i\mid\theta_j).
\end{equation*}
The generic allocation and truncation updates are
\begin{align*}
\Pp(z_i=j\mid k_i,\cdots)
&\propto \frac{w_j}{\xi_j}f(x_i\mid\theta_j),
&&j=1,\ldots,k_i,\\
\Pp(k_i=k\mid z_i=j,\cdots)
&=\frac{\xi_k-\xi_{k+1}}{\xi_j}\mathbf1\{k\geq j\}.
\end{align*}
Thus $k_i=z_i+S_i$ with the overshoot law (45).
For an exogenous constant-ratio schedule $\xi_j=cq^{j-1}$,
\begin{equation*}
k_i=z_i+S_i,
\qquad S_i\sim\Geom_0(1-q).
\end{equation*}
Equivalently, if $\xi_j=\exp(-\eta j)$ and $U_i\sim\mathrm{Unif}(0,1)$,
\begin{equation*}
k_i=\left\lfloor z_i-\eta^{-1}\log U_i\right\rfloor.
\end{equation*}
For a generic deterministic schedule, draw $U_i$ and choose the smallest $k\geq z_i$ satisfying
\begin{equation*}
\xi_{k+1}\leq(1-U_i)\xi_{z_i}.
\end{equation*}

\par\medskip\Needspace{3\baselineskip}\noindent\textit{Exogenous schedules and stick-breaking weights.}\enspace
Suppose
\begin{equation*}
w_1=v_1,
\qquad w_j=v_j\prod_{\ell<j}(1-v_\ell),
\qquad v_j\sim\Beta(\alpha_j,\beta_j)
\end{equation*}
independently, and the schedule does not depend on the $v_j$'s.  Define
\begin{equation*}
r_j=\sum_{i=1}^n\mathbf1\{z_i>j\}.
\end{equation*}
Collecting powers of $v_j$ and $1-v_j$ in (49) gives
\begin{equation*}
v_j\mid\zb,\cdots\sim\Beta(\alpha_j+n_j,\beta_j+r_j).
\end{equation*}
For a Dirichlet process, $(\alpha_j,\beta_j)=(1,\alpha)$.
When the schedule is independent of the hyperparameters being
updated, their full conditionals retain the form of the
underlying stick-breaking sampler. If the schedule depends
on an updated hyperparameter, that dependence must also be
accounted for in its full conditional.

\par\medskip\Needspace{3\baselineskip}\noindent\textit{Endogenous tail-mass schedule.}\enspace
For proper stick-breaking weights with $0<v_j<1$ almost surely for every $j$, choose $\xi_j=R_{j-1}=\prod_{\ell<j}(1-v_\ell)$.  Then $w_j/\xi_j=v_j$ and $\xi_k-\xi_{k+1}=w_k$, so
\begin{align}
\Pp(z_i=j\mid k_i,\cdots)&\propto v_jf(x_i\mid\theta_j),\qquad j\leq k_i, \notag\\
\Pp(k_i=k\mid z_i=j,\wb)&=\frac{w_k}{R_{j-1}}\mathbf1\{k\geq j\}.
\label{eq:k-endo-supp}
\end{align}
Let
\begin{equation*}
m_j=\sum_{i=1}^n\mathbf1\{k_i=j\},
\qquad
h_j=\sum_{i=1}^n\mathbf1\{k_i>j\}.
\end{equation*}
Under the same independent-beta prior, the factor $v_{z_i}$ contributes $n_j$ successes, $w_{k_i}=v_{k_i}R_{k_i-1}$ contributes $m_j$ additional successes and $h_j$ failures, and therefore
\begin{equation*}
v_j\mid\zb,\kb,\cdots
\sim\Beta(\alpha_j+n_j+m_j,\beta_j+h_j).
\end{equation*}

To sample \eqref{eq:k-endo-supp}, draw $U_i\sim\mathrm{Unif}(0,1)$ and start with tail mass $R_{z_i-1}$.  Move forward until the cumulative mass $\sum_{\ell=z_i}^k w_\ell$ exceeds $U_iR_{z_i-1}$.  If the search passes the current $k^\star$, generate new sticks and atoms from their prior conditionals.  

\par\medskip\Needspace{3\baselineskip}\noindent\textit{Geometric stick-breaking specialization.}\enspace
If $v_j\equiv v$ and $w_j=v(1-v)^{j-1}$, the tail-mass schedule is $\xi_j=(1-v)^{j-1}$.  Then
\begin{equation*}
\Pp(z_i=j\mid k_i,\cdots)\propto f(x_i\mid\theta_j),
\qquad
k_i=z_i+S_i,
\quad S_i\sim\Geom_0(v).
\end{equation*}
With $v\sim\Beta(a,b)$, each observation contributes $v^2(1-v)^{k_i-1}$ to the joint kernel, giving
\begin{equation}\label{eq:v-gsb-supp}
v\mid\kb,\cdots
\sim\Beta\!\left(a+2n,\ b+\sum_{i=1}^n(k_i-1)\right).
\end{equation}

\par\medskip\Needspace{3\baselineskip}\noindent\textit{Posterior coordinates from reverse stick-breaking.}\enspace
Fix a schedule $\xib$ independent of the weights and put
$H=k^\star$. Suppose the marginal prior of $(w_1,\ldots,w_H)$
has a Lebesgue density $\varpi_H$ on
$w_j>0$, $\sum_{j=1}^H w_j<1$.
Write $s=s_H$, $\lambda_j=a_j/s_j$ and
$\boldsymbol\lambda=(\lambda_2,\ldots,\lambda_H)$.
Since $\lambda_1=1$, the reverse product in
Proposition~4.1 gives
\[
 b_j(\boldsymbol\lambda)
 =\lambda_j\prod_{\ell=j+1}^H(1-\lambda_\ell)
 =\widetilde w_j^{(H)},\qquad
 w_j=\xi_jsb_j(\boldsymbol\lambda),\qquad 1\leq j\leq H.
\]
The transformation from $(s,\boldsymbol\lambda)$ to
$(w_1,\ldots,w_H)$ has absolute Jacobian
\begin{equation}\label{eq:reverse-jacobian}
 \mathcal J_H
 =s^{H-1}\left(\prod_{j=1}^H\xi_j\right)
       \prod_{\ell=2}^H(1-\lambda_\ell)^{\ell-2}.
\end{equation}
Indeed, $b_j$ depends only on $\lambda_j,\ldots,\lambda_H$,
so the Jacobian matrix of $(b_2,\ldots,b_H)$ with respect
to $(\lambda_2,\ldots,\lambda_H)$ is upper triangular, with
\[
 \left|\frac{\partial(b_2,\ldots,b_H)}
                   {\partial(\lambda_2,\ldots,\lambda_H)}\right|
 =\prod_{j=2}^H\prod_{\ell=j+1}^H(1-\lambda_\ell)
 =\prod_{\ell=2}^H(1-\lambda_\ell)^{\ell-2}.
\]
The scale transformation $a_j=sb_j$, with
$b_1=1-\sum_{j=2}^Hb_j$, contributes $s^{H-1}$, while
the fixed rescaling $w_j=\xi_ja_j$ contributes
$\prod_{j=1}^H\xi_j$, proving \eqref{eq:reverse-jacobian}.

Since $\sum_{j=1}^Hw_j=s\sum_{j=1}^H\xi_jb_j(\boldsymbol\lambda)$,
the transformed domain is
\begin{equation}\label{eq:reverse-domain}
 0<\lambda_\ell<1\quad(2\leq\ell\leq H),\qquad
 0<s<\left\{\sum_{j=1}^H\xi_jb_j(\boldsymbol\lambda)\right\}^{-1}.
\end{equation}
The induced prior density is therefore
$\varpi_H((\xi_jsb_j(\boldsymbol\lambda))_{j\leq H})\mathcal J_H$
on \eqref{eq:reverse-domain}, and vanishes wherever the
original prior density does.

Writing $n_j=\sum_{i=1}^n\mathbf1\{z_i=j\}$, the
weight-dependent factor in (49) becomes
\[
 \prod_{j=1}^H(w_j/\xi_j)^{n_j}
 =s^n\prod_{j=1}^H b_j(\boldsymbol\lambda)^{n_j}.
\]
Combining this factor with the induced prior gives
\begin{equation*}
 \begin{aligned}
 &p(s,\boldsymbol\lambda\mid
      \bm x,\zb,\kb,\thetab,\xib)\\
 &\qquad\propto
 \varpi_H\bigl((\xi_jsb_j(\boldsymbol\lambda))_{j\leq H}\bigr)
 \mathcal J_H\,s^n\prod_{j=1}^H b_j(\boldsymbol\lambda)^{n_j},
 \end{aligned}
\end{equation*}
on \eqref{eq:reverse-domain}. Here $H$ is fixed by $\kb$,
and the weights beyond $H$ have been integrated out in
$\varpi_H$. All remaining likelihood factors are constant
in $(s,\boldsymbol\lambda)$ under the stated conditioning.

For the paired-weight prior, $\Pp(w_1=w_2)>0$, so the law
of $(w_1,\ldots,w_H)$ has a singular component whenever
$H\geq2$ and does not admit the assumed Lebesgue density.
The transformation nevertheless induces a prior measure
$\Pi_H^{\rm rev}$, whose posterior satisfies
\[
 \begin{aligned}
 &\Pi_H^{\rm rev}(\dd s,\dd\boldsymbol\lambda\mid
       \bm x,\zb,\kb,\thetab,\xib)\\
 &\qquad\propto s^n\prod_{j=1}^H
       b_j(\boldsymbol\lambda)^{n_j}\,
       \Pi_H^{\rm rev}(\dd s,\dd\boldsymbol\lambda).
 \end{aligned}
\]
This measure formulation also accommodates weights generated
by a shared geometric stick-breaking fraction. Finite-support
laws require the appropriate boundary domain, while endogenous
schedules require a separate calculation of both the
transformation and the schedule-dependent likelihood factors.

\subsection{Slice updates and the number of available labels}\label{sec:matched}\label{sec:kstar-proof}
For a fixed deterministic schedule, the discrete index $K$
identifies the interval containing the continuous slice variable
$U$. We establish this correspondence and show that the two
augmentations yield identical Gibbs transitions when their
update rules are matched.

\par\medskip\Needspace{3\baselineskip}\noindent\textit{Matching the discrete and continuous scans.}\enspace
For a fixed deterministic schedule $\xib$, the slice
augmentation of \citet{Kalli2011} pairs each allocation
$z_i$ with a continuous variable $U_i$ through
\begin{equation*}
 \Pp(z_i=j,U_i\in\dd u\mid x_i,\wb,\thetab,\xib)
 \propto \frac{w_j}{\xi_j}f(x_i\mid\theta_j)
          \mathbf1\{0<u\leq\xi_j\}\dd u.
\end{equation*}
With $L(u)=k$ for $u\in(\xi_{k+1},\xi_k]$, integration
over the interval gives
\begin{equation*}
 \Pp(z_i=j,L(U_i)=k\mid x_i,\wb,\thetab,\xib)
 \propto \frac{w_j}{\xi_j}f(x_i\mid\theta_j)
       (\xi_k-\xi_{k+1})\mathbf1\{j\leq k\}.
\end{equation*}
Thus the finite-index augmentation records the interval
containing each slice variable. For
$\boldsymbol u=(u_1,\ldots,u_n)$, define the corresponding
map between sampler states by
\[
 \mathcal T(\wb,\thetab,\zb,\boldsymbol u)
 =\bigl(\wb,\thetab,\zb,(L(u_i))_{i=1}^n\bigr).
\]

\begin{proposition}[Equivalence of matched finite-index and slice Gibbs samplers]
\label{prop:matched-supp}
Fix a deterministic schedule and compare two Gibbs scans
with the same blocks and update order. The finite-index
scan updates $k_i$ using (51), while the
slice scan draws the $U_i$ independently from
$\mathrm{Unif}(0,\xi_{z_i})$. All other updates use the
same transition rules, with slice values entering only
through $(L(u_i))_{i=1}^n$, and the same rules for generating
and retaining weights and atoms.

Their one-iteration transition kernels satisfy
\[
 P_{\rm slice}\bigl(y,\mathcal T^{-1}(A)\bigr)
 =P_{\rm fin}\bigl(\mathcal T(y),A\bigr)
\]
for every admissible slice state $y$ and measurable set $A$
of finite-index states. Starting from states related by
$\mathcal T$, the two chains can be coupled so that this
relation holds after every iteration.
\end{proposition}

\begin{proof}
For $\kb=(L(u_i))_{i=1}^n$, the finite-index and slice
posterior kernels differ only by the factor
$\prod_{i=1}^n(\xi_{k_i}-\xi_{k_i+1})$, which is constant
in $(\wb,\thetab,\zb)$ because the schedule is fixed.
Their corresponding Gibbs conditionals therefore agree.
Moreover, the independent uniform slice draws give
\[
 \Pp\{L(U_i)=k\mid z_i,\xib\}
 =\frac{\xi_k-\xi_{k+1}}{\xi_{z_i}}
       \mathbf1\{k\geq z_i\},
\]
so the auxiliary updates also agree after applying
$\mathcal T$. Composing these matched updates proves
the kernel identity. The pathwise coupling follows by
setting $k_i=L(U_i)$ at each auxiliary update and coupling
all remaining updates to have identical outputs.
\end{proof}

\par\medskip\Needspace{3\baselineskip}\noindent\textit{Conditional law of the active dimension.}\enspace
For a fixed realized schedule and fixed allocations $z_{1:n}$, the variables $k_i$ are conditionally independent and
\begin{equation*}
\Pp(k_i\leq H\mid z_i=j,\xib)
=\sum_{k=j}^H\frac{\xi_k-\xi_{k+1}}{\xi_j}
=1-\frac{\xi_{H+1}}{\xi_j},
\qquad H\in\N,\ H\geq j.
\end{equation*}
Thus, for $k^\star=\max_i k_i$ and integer $H\geq z^\star=\max_i z_i$,
\begin{equation}\label{eq:kstar-cdf}
\Pp(k^\star\leq H\mid z_{1:n},\xib)
=\prod_{i=1}^n\left(1-\frac{\xi_{H+1}}{\xi_{z_i}}\right).
\end{equation}
When $\xi_j=cq^{j-1}$, $S_i=k_i-z_i$ are i.i.d. geometric-zero with $\Pp(S_i\geq t)=q^t$ for integers $t\geq0$.  The union bound gives
\begin{equation*}
\Pp\!\left(\max_iS_i\geq t\right)\leq nq^t.
\end{equation*}

Under the deterministic
geometric schedule $\xi_j=cq^{j-1}$, the union bound yields
$\Pp(k^\star-z^\star\geq t\mid\zb)\leq nq^t$ for integers
$t\geq0$. Thus, for $0<\delta<1$,
\begin{equation}\label{eq:tdelta}
 t_\delta=\left\lceil\frac{\log(n/\delta)}{\log(1/q)}\right\rceil
 \quad\text{implies}\quad
 \Pp(k^\star\leq z^\star+t_\delta-1\mid\zb)\geq1-\delta.
\end{equation}
For fixed $q$, the extension beyond $z^\star$ therefore grows
at most logarithmically in $n$, in probability. The total excess
over the occupied count also includes unoccupied labels below
$z^\star$, giving
\begin{equation}\label{eq:overhead-decomposition}
 k^\star-c_n=(k^\star-z^\star)+(z^\star-c_n).
\end{equation}

\subsection{Gaussian mixture calculations and numerical evaluation}\label{sec:gaussian}\label{sec:repro}
Consider a Gaussian mixture with observation standard deviation
$\omega>0$ and component means drawn independently from a
normal base distribution,
\begin{equation*}
 x_i\mid z_i,\thetab\sim N(\theta_{z_i},\omega^2),
 \qquad
 \theta_j\stackrel{\mathrm{iid}}\sim N(\mu_0,\tau_0^2).
\end{equation*}
The parameters $\mu_0$ and $\tau_0>0$ are fixed, while the
mixing weights follow either a Dirichlet or Pitman--Yor prior
with fixed parameters $(\alpha,\sigma)$, or a geometric
stick-breaking prior with $v_j\equiv V\sim\Beta(a_0,b_0)$.
Note that here $\sigma$ denotes the Pitman--Yor discount.
\par\medskip\Needspace{3\baselineskip}\noindent\textit{Conjugate component updates and the predictive tail.}\enspace
Given the allocations, normal conjugacy yields
\begin{equation}\label{eq:normal-fullconditional}
 \tau_j^2=\left(\tau_0^{-2}+n_j\omega^{-2}\right)^{-1},
 \qquad
 \mu_j=\tau_j^2\left(\frac{\mu_0}{\tau_0^2}
 +\frac{\sum_{i:z_i=j}x_i}{\omega^2}\right).
\end{equation}
The component parameters are then conditionally independent,
with $\theta_j\mid\zb,\bm x\sim N(\mu_j,\tau_j^2)$,
reducing to the base distribution $G_0$ when $n_j=0$.
The corresponding weight and auxiliary-index updates are
given in Supplement~\ref{sec:conditionals}. In particular,
for geometric stick-breaking under a fixed exogenous schedule,
the shared stick-breaking fraction has conditional distribution
\begin{equation*}
 V\mid\zb,\cdots\sim
 \Beta\!\left(a_0+n,\ b_0+\sum_i(z_i-1)\right),
\end{equation*}
which differs from \eqref{eq:v-gsb-supp} under the
endogenous tail-mass schedule.

For predictive inference, the contribution of atoms not yet
generated can be integrated analytically. At the current
maximum index $H=k^\star$, let
$R_H=1-\sum_{j\leq H}w_j$ denote their remaining mass.
Conditional on the generated weights and atoms, integrating
over the independent $G_0$ locations gives
\begin{equation}\label{eq:predictive-tail}
 \begin{aligned}
 \overline f_H(t)
 &=\E\{f_G(t)\mid (w_j,\theta_j)_{j\leq H},\zb,\kb,\bm x\}\\
 &=\sum_{j=1}^H w_j\,\varphi(t;\theta_j,\omega^2)
 +R_H\,\varphi(t;\mu_0,\omega^2+\tau_0^2),
 \end{aligned}
\end{equation}
where $\varphi(\cdot;\mu,s^2)$ is the normal density with
variance $s^2$. Averaging \eqref{eq:predictive-tail} over
posterior draws provides a Rao--Blackwellized estimate of
the posterior predictive density. Credible bands for the
random density $f_G$, however, must also account for the
conditional variability of the uninstantiated tail, which
is absent from the variation of these conditional means.

\par\medskip\Needspace{3\baselineskip}\noindent\textit{Exact posterior enumeration.}\enspace
For a partition $\mathcal P$ of the observation indices, the
component parameters can be integrated out separately within
each block. For every nonempty $B\in\mathcal P$, let
$n_B=|B|$, $d_i=x_i-\mu_0$,
$D_B=\sum_{i\in B}d_i$ and $Q_B=\sum_{i\in B}d_i^2$.
Normal conjugacy gives the log marginal likelihood
\begin{equation*}
\begin{split}
\ell_B
&=\log\int\prod_{i\in B}\varphi(x_i;\theta,\omega^2)
             \varphi(\theta;\mu_0,\tau_0^2)\,\dd\theta\\
&=-\frac{n_B}{2}\log(2\pi\omega^2)
-\frac12\log\left(1+\frac{n_B\tau_0^2}{\omega^2}\right)\\
&\quad-\frac1{2\omega^2}\left(
 Q_B-\frac{\tau_0^2D_B^2}{\omega^2+n_B\tau_0^2}\right).
\end{split}
\end{equation*}
Under a PY prior, the probability of a particular
partition $\mathcal P=\{B_1,\ldots,B_c\}$ is given by its
exchangeable partition probability function
\begin{equation}\label{eq:PY-eppf-reference}
p_{\alpha,\sigma}(\mathcal{P})=
\frac{\Gamma(\alpha+1)}{\Gamma(\alpha+n)}
\left\{\prod_{r=1}^{c-1}(\alpha+r\sigma)\right\}
\prod_{B\in\mathcal{P}}\frac{\Gamma(|B|-\sigma)}{\Gamma(1-\sigma)}.
\end{equation}
Here an empty product equals one, and $\sigma=0$ with
$\alpha>0$ recovers the DP. Since the EPPF assigns
probability to each individual set partition, every
partition is included exactly once in the enumeration
\citep{Pitman2006,PitmanYor1997}.

Combining these prior probabilities with the integrated
block likelihoods gives
\begin{equation}\label{eq:exact-partition-posterior}
\Pi(\mathcal{P}\mid\bm x)=
\frac{p_{\alpha,\sigma}(\mathcal{P})\exp\{\sum_{B\in\mathcal{P}}\ell_B\}}
 {\sum_{\mathcal{P}'}p_{\alpha,\sigma}(\mathcal{P}')\exp\{\sum_{B\in\mathcal{P}'}\ell_B\}}.
\end{equation}
The denominator sums over all partitions of
$\{1,\ldots,n\}$, of which there are $203$ when $n=6$.
The posterior occupancy probabilities
$\Pi(c_n=c\mid\bm x)$ follow by summing over partitions
with $c$ blocks.

For prediction, let $(\mu_B,\tau_B^2)$ denote the posterior
mean and variance of the component parameter associated
with block $B$, obtained from \eqref{eq:normal-fullconditional}.
The Pitman--Yor predictive rule then gives
\begin{equation*}
\begin{split}
f(t\mid\mathcal{P},\bm x)
&=\frac{\alpha+\sigma c}{\alpha+n}
\varphi(t;\mu_0,\omega^2+\tau_0^2)\\
&\quad+\sum_{B\in\mathcal{P}}\frac{|B|-\sigma}{\alpha+n}
\varphi(t;\mu_B,\omega^2+\tau_B^2).
\end{split}
\end{equation*}
Averaging this expression with respect to
\eqref{eq:exact-partition-posterior} yields the posterior
predictive density. Both the occupancy probabilities and
the predictive density are therefore available by finite
summation, independently of the admissible visibility
schedule used to represent the same SSP prior.

\par\medskip\Needspace{3\baselineskip}\noindent\textit{Numerical evaluation.}\enspace
To maintain numerical stability, survival probabilities, stick
products and categorical weights are evaluated on the log scale.
In the $K$-first sampler, (47) accounts for the remaining
probability mass of the infinite index law, while normalization
is performed only after a finite conditional distribution, such
as $J\mid K=k$, has been determined. If the remaining mass
cannot yet be resolved, the computation must continue rather
than redistribute it among the generated atoms, which would
alter the sampling law. Numerical draws outside their prescribed
support are treated as errors. The original stick-breaking labels
are preserved, with order-of-appearance labels introduced
separately when displaying partitions.

For a fixed geometric schedule $\xi_j=cq^{j-1}$, the update
\[
 k_i=z_i+\left\lfloor\frac{\log U_i}{\log q}\right\rfloor,
 \qquad U_i\sim\mathrm{Unif}(0,1),
\]
and the corresponding continuous slice intervals can both be
evaluated in log coordinates, avoiding underflow when
$\xi_{z_i}$ is very small. Matched implementations retain the
same Gibbs scan order. Finally, levels with $s_k=0$ have zero
index probability and require no normalization by $s_k$.

\subsection{Partition and predictive calculations for paired weights}
\label{app:paired-predictive-proof}
Summing over the possible latent pairings gives finite
expressions for partition probabilities and prediction.
We use these expressions to prove Proposition~3.11, derive the associated sample summaries and describe their numerical evaluation. Since the atom locations are independent draws from diffuse $G_0$, coincidences among sampled locations correspond almost surely to repeated observations of the same species.

\par\medskip\Needspace{3\baselineskip}\noindent\textit{Partition probabilities and the latent matching.}\enspace
Let $\mathcal P=\{B_1,\ldots,B_c\}$ be a specified partition with
block sizes $\boldsymbol n=(n_1,\ldots,n_c)$. Write
$\thetab_{\zb}=(\theta_{z_1},\ldots,\theta_{z_n})$ for a latent
sample with this partition and distinct values
$\theta_1^*,\ldots,\theta_c^*$. A realization $M$ of the latent matching $\mathcal M$ is a
possibly empty collection of disjoint pairs of block indices,
each identifying two observed species arising from the same
original DP weight. Write $\mathcal U(M)$ for the unpaired block indices. Merging the two blocks in every pair leaves a particular partition into $c-|M|$ blocks associated with the original DP weights. Its Dirichlet-process probability is
\[
 \frac{\Gamma(\alpha)}{\Gamma(\alpha+n)}
 \alpha^{c-|M|}
 \prod_{i\in\mathcal U(M)}\Gamma(n_i)
 \prod_{\{i,j\}\in M}\Gamma(n_i+n_j),
\]
as follows from the Dirichlet EPPF \citep{Pitman1996,Pitman2006}.

Within a pair, either species has probability $1/2$. An unmatched
block can correspond to either species, while two matched blocks
can be assigned to the two species in either order. Each pair
represented in the sample therefore contributes a factor two,
and the $n$ individual allocations contribute $2^{-n}$. Hence, writing
\[
 W_\alpha(M;\boldsymbol n)=
 (2\alpha)^{c-|M|}
 \prod_{i\in\mathcal U(M)}\Gamma(n_i)
 \prod_{\{i,j\}\in M}\Gamma(n_i+n_j),
\]
we obtain
\begin{equation}\label{eq:paired-joint-matching}
 \Pp(\mathcal P_n=\mathcal P,\mathcal M=M)
 =\frac{\Gamma(\alpha)}{2^n\Gamma(\alpha+n)}
 W_\alpha(M;\boldsymbol n).
\end{equation}
Distinct matchings describe disjoint partitions of the sample among the original DP weights. Summing \eqref{eq:paired-joint-matching} over them gives the EPPF in (38), while Bayes' rule gives
\begin{equation}\label{eq:paired-matching-posterior}
 \Pp(\mathcal M=M\mid\thetab_{\zb})
 =\frac{W_\alpha(M;\boldsymbol n)}
        {\sum_{M'}W_\alpha(M';\boldsymbol n)}.
\end{equation}
The observed values contribute no additional information about the matching because all atom locations are independent draws from the same diffuse base measure.

\par\medskip\Needspace{3\baselineskip}\noindent\textit{Prediction given a matching.}\enspace
Let $n_{g_M(i)}$ be the total sample count for the pair containing
species $i$, equal to $n_i$ when $i\in\mathcal U(M)$ and to
$n_i+n_j$ when $\{i,j\}\in M$. Conditional on the matching,
the DP posterior assigns mean total weight $n_g/(\alpha+n)$
to each observed pair and $\alpha/(\alpha+n)$ to the unobserved
pairs. Each pair
shares its weight equally between its two species, giving
\begin{equation}\label{eq:paired-source-predictive}
 \begin{aligned}
 \Pp(\theta_{z_{n+1}}=\theta_i^*\mid\thetab_{\zb},\mathcal M=M)
 &=\frac{n_{g_M(i)}}{2(\alpha+n)},\\
 \Pp(\theta_{z_{n+1}}\notin\{\theta_1^*,\ldots,\theta_c^*\}
          \mid\thetab_{\zb},\mathcal M=M)
 &=\frac{\alpha+\tfrac12\sum_{i\in\mathcal U(M)}n_i}{\alpha+n}.
 \end{aligned}
\end{equation}
The new-species probability combines mass in entirely unobserved pairs with the unobserved partners of unpaired species. Averaging \eqref{eq:paired-source-predictive} under \eqref{eq:paired-matching-posterior} gives prediction conditional on $\thetab_{\zb}$. These probabilities sum to one for every matching and therefore also after averaging. They agree with the EPPF ratios for an increment to an existing block or the addition of a singleton block, completing the proof of Proposition~3.11.

For a kernel mixture, let $\ell_B=\log\int\prod_{i\in B} f(x_i\mid\theta)G_0(\dd\theta)$ whenever this integral is available. The joint posterior on the observation partition and its latent matching is
\begin{equation}\label{eq:paired-mixture-partition-posterior}
 \Pp(\mathcal P_n=\mathcal P,\mathcal M=M\mid\bm x)
 \propto W_\alpha(M;\boldsymbol n)
         \exp\!\left\{\sum_{B\in\mathcal P}\ell_B\right\}.
\end{equation}
Here $\alpha$ is fixed, so the factor $\Gamma(\alpha)/\{2^n\Gamma(\alpha+n)\}$ is common to all terms. This supplies a marginal route to inference, while the finite-index likelihood retains the weights and atoms. In the Gaussian model of Supplement~\ref{sec:gaussian}, the previously derived block likelihoods apply directly, and replacing the PY EPPF by (38) gives exact posterior enumeration for the paired prior.

\par\medskip\Needspace{3\baselineskip}\noindent\textit{Coincidence, occupancy and discovery.}\enspace
Equal splitting multiplies the $r$th power sum of the original DP weights by $2^{1-r}$. The Dirichlet power-sum identity consequently gives
\begin{equation}\label{eq:paired-r-coincidence}
 \Pp(\theta_{z_1}=\cdots=\theta_{z_r})
 =\E\sum_j w_j^r
 =2^{1-r}\frac{\Gamma(\alpha+1)\Gamma(r)}{\Gamma(\alpha+r)},
 \qquad r\geq2.
\end{equation}
At $r=2$, this gives the prior coincidence probability
$1/\{2(\alpha+1)\}$ in (41). The equality $\E(w_1)=\E\sum_jw_j^2$ also follows from size-biased invariance and agrees with $2w_1\sim\Beta(1,\alpha)$.

Let $d_n$ denote the probability that observation $n+1$ belongs to an unobserved species, averaged over the preceding sample. Conditional on the paired weights it is $\sum_jw_j(1-w_j)^n$. A size-biased draw from the original DP weights has a $\Beta(1,\alpha)$ distribution, so
\begin{equation}\label{eq:paired-discovery-supp}
 \begin{aligned}
 d_n
 &=\Pp(\theta_{z_{n+1}}\notin\{\theta_{z_1},\ldots,\theta_{z_n}\})\\
 &=\alpha\int_0^1(1-x/2)^n(1-x)^{\alpha-1}\dd x\\
 &=\frac{\alpha}{2^n}\sum_{h=0}^n\binom nh\frac1{\alpha+h}
 =\E\!\left(\frac{\alpha}{\alpha+B_n}\right),
 \qquad B_n\sim\operatorname{Binomial}(n,1/2).
 \end{aligned}
\end{equation}
Expanding $(1+y)^n$ after setting $y=1-x$ gives the positive
finite sum in (42) and its binomial
expectation representation. Writing $M_{1,n}$ for the number of singleton species, summing the indicators of discovery and singleton occurrence gives
\[
 \E(c_n)=\sum_{m=0}^{n-1}d_m,
 \qquad \E(M_{1,n})=n d_{n-1}.
\]
To derive (43), view each draw as
selecting an original DP weight and then, independently,
choosing one of its two species with equal probability.
The number $B_n$ assigned to the first species of their
respective pairs is therefore $\operatorname{Binomial}(n,1/2)$
and is independent of the DP weights. A DP sample of size
$b$ has expected occupancy
$m_\alpha(b)=\alpha\{\psi(\alpha+b)-\psi(\alpha)\}$,
where $\psi=\Gamma'/\Gamma$ and $m_\alpha(0)=0$
\citep[Section~3.3]{Pitman2006}.

Conditional on $B_n=b$, the two groups have expected
occupancies $m_\alpha(b)$ and $m_\alpha(n-b)$. Although
their counts may be dependent through the shared DP weights,
linearity of expectation allows their contributions to be
added, and binomial symmetry gives
\begin{equation}\label{eq:paired-occupancy-supp}
 \begin{aligned}
 \E(c_n)
 &=\E\{m_\alpha(B_n)+m_\alpha(n-B_n)\}\\
 &=\frac{2\alpha}{2^n}\sum_{b=0}^n\binom nb
       \{\psi(\alpha+b)-\psi(\alpha)\}.
 \end{aligned}
\end{equation}

These are prior summaries averaged over the sample, whereas \eqref{eq:paired-source-predictive} conditions on its observed partition. At $\alpha=1$, \eqref{eq:paired-discovery-supp} gives
$d_n=2\{1-2^{-(n+1)}\}/(n+1)$, while
\eqref{eq:paired-occupancy-supp} reduces to
$\E(c_n)=2\sum_{r=1}^n(1-2^{-r})/r$. After block sizes $(5,1)$ and $(3,3)$, averaging over their possible matchings instead gives new-species probabilities $13/49$ and $19/112$, respectively.

\begin{figure}[tbp]\centering
\includegraphics[width=.84\linewidth]{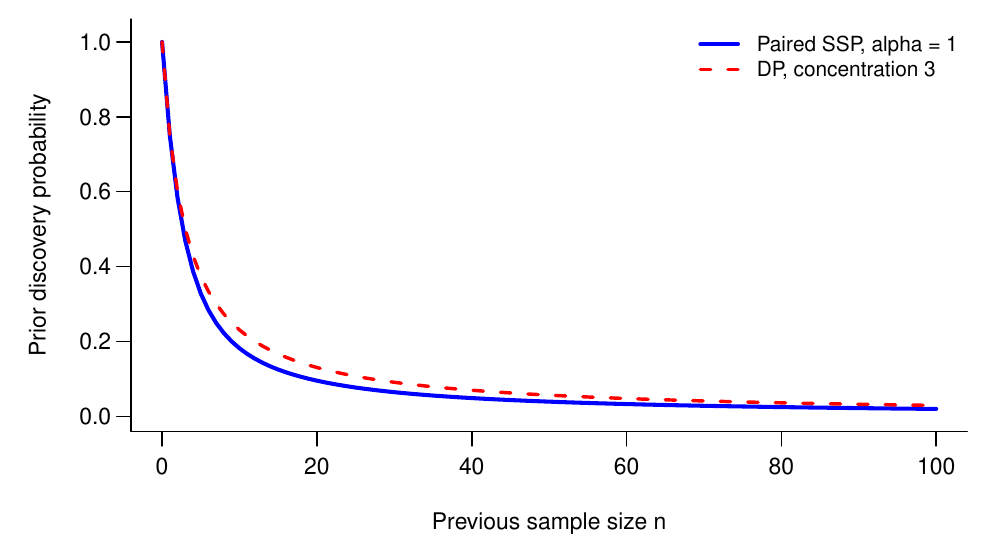}
\caption{Prior discovery probabilities for paired weights with $\alpha=1$ and a DP with total mass $3$. Both have prior tie probability $1/4$, while their discovery curves follow \eqref{eq:paired-discovery-supp} and $3/(3+n)$, respectively.}
\label{fig:paired-discovery}
\end{figure}
Figure~\ref{fig:paired-discovery} compares the prior discovery probabilities with those of a DP having the same tie probability. The different curves demonstrate that eliciting a single recurrence summary leaves room for distinct discovery behavior.

\par\medskip\Needspace{3\baselineskip}\noindent\textit{Finite evaluation.}\enspace
For a subset $A\subseteq\{1,\ldots,c\}$, define $H(\varnothing)=1$. Choosing its least element $i$ and separating matchings according to whether $i$ is unpaired or paired with some $j\in A\setminus\{i\}$ gives
\begin{equation}\label{eq:paired-matching-recursion}
 H(A)=2\alpha\Gamma(n_i)H(A\setminus\{i\})
 +\sum_{j\in A\setminus\{i\}}2\alpha\Gamma(n_i+n_j)
                  H(A\setminus\{i,j\}).
\end{equation}
The recursion gives the EPPF as
$2^{-n}\Gamma(\alpha)H(\{1,\ldots,c\})/\Gamma(\alpha+n)$
using only positive terms, which can be accumulated on the
logarithmic scale without cancellation from alternating
signs. Storing and reusing the value for each subset requires
at most $2^c$ states, hence making the calculation practical when
the number of observed species is modest.

The marginal-index construction in
Section~3.7 also permits direct evaluation
of the index probabilities through one-dimensional integration
and finite sums. Let $b(x)$ and
$\lambda_1(x),\lambda_2(x)$ denote its Bernoulli probability
and Poisson means. For $0<x<1$ and integers $r\geq0$, write
\[
 h_r(x)=\Pp(N_1+2N_2=r\mid X=x)
 =e^{-\lambda_1(x)-\lambda_2(x)}
  \sum_{a=0}^{\lfloor r/2\rfloor}
  \frac{\lambda_1(x)^{r-2a}\lambda_2(x)^a}
       {(r-2a)!a!}.
\]
With $h_{-1}(x)=0$, averaging over the Bernoulli variable
and then over $X$ gives, for $j\geq1$,
\begin{equation}\label{eq:paired-mean-quadrature}
 \overline w_j=\Pp(D=j)
 =\int_0^1
   \bigl\{[1-b(x)]h_{j-1}(x)+b(x)h_{j-2}(x)\bigr\}\dd x.
\end{equation}
The coefficients can be computed successively from
$h_0(x)=e^{-\lambda_1(x)-\lambda_2(x)}$ using the
nonnegative recursion
\[
 r h_r(x)
 =\lambda_1(x)h_{r-1}(x)+2\lambda_2(x)h_{r-2}(x),
 \qquad r\geq1.
\]
The finite convolution in (35)
then yields each $\pi_k$.

The remaining index mass admits an explicit bound.
For any integer $H\geq0$ and
$1<z<\min\{\sqrt{1+1/\alpha},1/q\}$,
\begin{equation}\label{eq:paired-tail-bound}
 \Pp(K>H)\leq C(z)z^{-H-1},\qquad
 C(z)=\frac{z^2(1-q)}
 {[1-\alpha(z^2-1)](1-qz)}.
\end{equation}
To obtain this bound, combine the Poisson factors in
\eqref{eq:paired-pgf} to write
\[
 \exp\{\lambda_1(x)(z-1)+\lambda_2(x)(z^2-1)\}
 =x^{-\alpha(z^2-1)}(2-x)^{-\alpha(z-1)^2}
 \leq x^{-\alpha(z^2-1)}.
\]
Since the Bernoulli factor is at most $z$ and
$K=D+S_0$ with independent $S_0\sim\Geom_0(1-q)$,
\[
 \E(z^K)
 \leq\frac{z^2(1-q)}{1-qz}
       \int_0^1x^{-\alpha(z^2-1)}\dd x
 =C(z).
\]
Finally, $K>H$ implies $z^K\geq z^{H+1}$, so Markov's
inequality gives \eqref{eq:paired-tail-bound}.

The tabulated predictive probabilities and sample summaries
are evaluated by finite sums, so quadrature is needed only
for individual mean weights or index probabilities. For these
quantities, the quadrature error estimates are reported
separately from the analytic bound on the remaining index mass.
The supplied code also verifies EPPF normalization over all
$203$ set partitions at $n=6$, the predictive addition rule,
and agreement between EPPF ratios and conditional predictions
averaged over the posterior distribution of the latent matching.
Throughout these calculations, the original infinite law is
retained without redistributing its remaining mass among the
generated atoms.

\subsection{Posterior work and slice overhead}\label{sec:illustration-kalli}\label{sec:work-illustration}
\begin{figure}[tbp]\centering
\includegraphics[width=.88\linewidth]{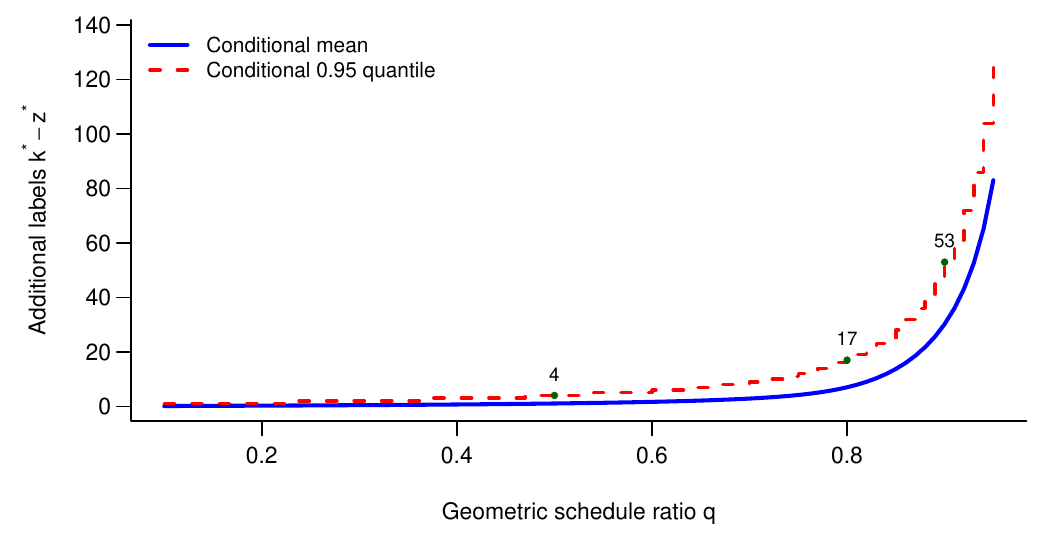}
\caption{Conditional mean and $0.95$ quantile of the
additional labels $k^\star-z^\star$ required for an allocation
update under the geometric schedule $\xi_j=cq^{j-1}$.
The same allocation vector, with $n=120$, $c_n=16$ and
$z^\star=26$, is used throughout.}
\label{fig:budget-optimized}
\end{figure}

We illustrate the effect of the schedule using a fixed
allocation vector $\zb$ with $n=120$, $c_n=16$ and
$z^\star=26$.
The occupied labels $(1,\ldots,12,15,18,22,26)$ have counts
$(30,20,15,12,10,8,6,5,4,3,2,1,1,1,1,1)$, respectively. For $\xi_j=cq^{j-1}$,
Figure~\ref{fig:budget-optimized} shows how the conditional
mean and $0.95$ quantile of $k^\star-z^\star$ vary with $q$,
as calculated from \eqref{eq:kstar-cdf}. At
$q=0.5,0.8,0.9$, the corresponding means are approximately
$1.04,7.00,30.20$, and the quantiles are $4,17,53$,
respectively. Thus varying the schedule changes the number
of additional labels required while preserving the
underlying SSP and its marginal posterior.

For any fixed $q$, \eqref{eq:tdelta} gives a conservative
upper bound on this extension with probability at least
$1-\delta$, without evaluating the full product CDF.
A related result is established by
\citet[Theorem~3.3]{FranzoliniGaffi2026}, who show that
a DP slice sampler with fixed concentration requires
$K_n-H_n=O_{\Pp}(\log n)$ additional instantiated components,
uniformly over the current sample partition. Their sampler
places the $H_n$ occupied components first, draws their
weights from the posterior Dirichlet law, and extends the
representation to $K_n$ components until the remaining mass
falls below the smallest weight-based slice variable.

The two bounds use different reference points.
Equation~\eqref{eq:tdelta} controls $k^\star-z^\star$
for any proper SSP under a deterministic geometric schedule,
whereas $K_n-H_n$ measures excess over occupancy.
In our ordering, the latter comparison also includes
the $z^\star-c_n$ unoccupied labels below $z^\star$,
as expressed by \eqref{eq:overhead-decomposition}.
The reference points coincide when $z^\star=c_n$, although
relabeling to achieve this equality requires consistent
treatment of the weights and schedule. Both bounds concern
representation size per iteration, while mixing rates and
extensions of the DP result to other SSPs with weight-dependent
schedules require separate analysis.

Remember that exact posterior enumeration for the six-observation Gaussian
example is given in Supplement~\ref{sec:gaussian}.

\bibliography{references_sspfin}